%% file: LieGroups.tex
\documentclass[11pt]{amsart}
\pdfoutput=1

\usepackage{amssymb,latexsym}
\usepackage{amsmath,amsfonts,amsthm}
\usepackage{stmaryrd}
\usepackage{enumerate}
\usepackage{enumitem}
\usepackage[T1]{fontenc}
\usepackage[latin1]{inputenc}
\usepackage[english]{babel}

\usepackage{accents}

\usepackage{graphicx}
\usepackage{texdraw}

\usepackage[bookmarks=true,%
    colorlinks=true,
    linkcolor=blue,%
    citecolor=blue,%
    filecolor=blue,%
    menucolor=blue,%
    urlcolor=blue,%
    breaklinks=true]{hyperref}

\usepackage{caption}

\usepackage[all]{xy}
\CompileMatrices
\def\cxymatrix#1{\xy*[c]\xybox{\xymatrix#1}\endxy}

\theoremstyle{plain}
\newtheorem{theorem}{Theorem}[section]
\theoremstyle{definition}
\newtheorem{proposition}[theorem]{Proposition}
\newtheorem{lemma}[theorem]{Lemma}
\newtheorem{definition}[theorem]{Definition}

\newtheorem{remark}[theorem]{Remark}
\newtheorem{corollary}[theorem]{Corollary}
\newtheorem{conjecture}[theorem]{Conjecture}

\newtheorem{example}[theorem]{Example}

\newtheorem{convention}[theorem]{Convention}

\numberwithin{equation}{section}

\newcommand{\Z}{\mathbb Z}
\newcommand{\Q}{\mathbb Q}
\newcommand{\R}{\mathbb R} 
\newcommand{\C}{\mathbb C}
\renewcommand{\H}{\mathbb H}

\newcommand{\T}{\mathcal T}

\newcommand{\A}{\mathcal A}

\newcommand{\M}{\mathcal M}

\newcommand{\g}{\mathfrak g}
\newcommand{\n}{\mathfrak n}
\newcommand{\h}{\mathfrak h}

\newcommand{\Matrix}[4]{\left(\begin{smallmatrix}#1&#2\\#3&#4\end{smallmatrix}\right)}
\newcommand{\Vector}[2]{\left(\begin{smallmatrix}#1\\#2\end{smallmatrix}\right)}

\newcommand{\wzero}{\overline{w_0}}

\DeclareMathOperator{\Spin}{Spin}
\DeclareMathOperator{\Sp}{Sp}
\DeclareMathOperator{\SL}{SL}

\DeclareMathOperator{\diag}{diag}

\DeclareMathOperator{\id}{id}

\DeclareMathOperator{\red}{red}

\DeclareMathOperator{\gen}{gen}

\DeclareMathOperator{\Hom}{Hom}
\DeclareMathOperator{\Conf}{Conf}
\DeclareMathOperator{\rot}{rot}
\DeclareMathOperator{\flip}{flip}

\usepackage{scalerel}

\usepackage{ifpdf}
\ifpdf
\fi

\usepackage[left=2.8cm,right=2.8cm,bottom=3.1cm,top=3cm]{geometry}

\title[Fock-Goncharov coordinates for rank two Lie groups]{Fock-Goncharov coordinates for rank two Lie groups}
\author{Christian K. Zickert}
\address{University of Maryland \\
         Department of Mathematics \\
         College Park, MD 20742-4015, USA \newline
         {\tt \url{http://www2.math.umd.edu/~zickert}}}
\email{zickert@math.umd.edu}

\thanks{The author was supported in part by DMS-13-09088. \\
\newline
2010 {\em Mathematics Classification.} Primary 32G15, 57M27. Secondary 13F60, 57M50.
}

\date{}
\begin{document}
\maketitle
\begin{abstract}
Let $G$ be a simply connected, simple, complex Lie group of rank $2$. We give explicit Fock-Goncharov coordinates for configurations of triples and quadruples of affine flags in $G$. We show that the action on triples by orientation preserving permutations corresponds to explicit quiver mutations, and that the same holds for the flip (changing the diagonal in a quadrilateral). This gives explicit coordinates on higher Teichm\"uller space, and also coordinates for boundary-unipotent representations of $3$-manifold groups. As an application, we compute the (generic) boundary-unipotent representations in $\Sp(4,\C)/\langle -I\rangle$ for the figure-eight knot complement.
\end{abstract}
\section{Introduction}
Let $G$ be a simply connected, semisimple, complex Lie group with adjoint group $G'$. For an oriented, punctured surface $S$ with negative Euler characteristic, Fock and Goncharov~\cite{FockGoncharov} define a pair $(\mathcal A_{G,S},\mathcal X_{G',S})$ of moduli spaces, which can be viewed as an ``algebraic-geometric avatar of Higher Teichm\"uller theory''~\cite[p.~5]{FockGoncharov}. We shall here only consider the space $\mathcal A_{G,S}$. The space $\A_{G,S}$ has a birational atlas with a chart $\A_{G,\T}$ for each ideal triangulation $\T$ together with an ordering of the vertices of each triangle compatible with the orientation of $S$. Each chart is a complex torus, and is constructed by gluing together copies of a configuration space of triples of affine flags in general position via a gluing pattern determined by the triangulation. Fock and Goncharov show that the atlas is \emph{positive}, i.e.~that the transition functions are subtraction free rational functions. This allows one to define the space of positive points of $\mathcal A_{G,S}$. When $G=\SL(2,\C)$ this space is Penner's decorated Teichm\"uller space~\cite{Penner}, and the positive coordinates coming from an ideal triangulation are Penner's $\lambda$-coordinates.

Our main result is that when $G$ is simple of rank 2, the transition functions are given by explicit quiver mutations. For this it is enough to consider a \emph{rotation} (a cyclic permutation of the vertex ordering of a triangle) and a \emph{flip} (change of the diagonal in a quadrilateral). We also give explicit algorithms for the transition functions in higher rank, and we conjecture that they are always given by quiver mutations. When $G=\SL(n,\C)$ explicit quiver mutations were given by Fock and Goncharov~\cite[Sec.~10]{FockGoncharov}.

Garoufalidis, Thurston and Zickert~\cite{GaroufalidisThurstonZickert} (see also \cite{BergeronFalbelGuilloux,DimofteGabellaGoncharov}) used the work of Fock and Goncharov to define coordinates (called \emph{Ptolemy coordinates}) for boundary-unipotent $\SL(n,\C)$-representations of 3-manifold groups. The relations between these coordinates (called \emph{Ptolemy relations}) are exactly the mutation relations found by Fock and Goncharov. The coordinates seem to be very efficient for concrete computations (see \cite{CURVE,FalbelKoseleffRouillier} for a database). Our main results provide similar coordinates for all simply connected, simple, complex Lie groups of rank~2. 

\section{Statement of results}\label{sec:Results}
Let $G$ be a simply connected, simple, complex Lie group of rank $2$, i.e., $G$ is either $A_2=\SL(3,\C)$, $B_2=\Spin(5,\C)\cong\Sp(4,\C)=C_2$ or $G_2$. There is a canonical central element $s_G\in G$, which is either trivial or of order $2$ (see Section~\ref{sec:sG}). It is trivial for $A_2$ and $G_2$, and non-trivial for $B_2$ and $C_2$.

Fix a maximal unipotent subgroup $N$ and let $\A=\A_G=G/N$ denote the affine flag variety of $N$-cosets in $G$. The diagonal left $G$ action on $\A^k$ does not have a geometric quotient, but if we restrict to tuples that are sufficiently generic (see Definition~\ref{def:SufGen}), there is a geometric quotient $\Conf_k^*(\A)$. It is a sub-variety of the algebro-geometric quotient $\A^k/\!/ G$.

To each of the groups $A_2$, $B_2$, $C_2$ and $G_2$ we associate a weighted quiver $Q_G$ (see Definition~\ref{def:Quiver}) of weight $1$, $2$, $2$, and $3$, respectively. We think of the graphs as being immersed in the plane (in fact, in a triangle), but the immersion only serves as a visual representation, providing a convenient labeling scheme, and is not formally part of the data.

\begin{figure}[htb]
\begin{center}
\begin{minipage}[c]{0.24\textwidth}
\scalebox{0.7}{\input{figures_gen/QuiverA2Ex.tex}}
\end{minipage}
\begin{minipage}[c]{0.24\textwidth}
\scalebox{0.7}{\input{figures_gen/QuiverB2Ex.tex}}
\end{minipage}
\begin{minipage}[c]{0.24\textwidth}
\scalebox{0.7}{\input{figures_gen/QuiverC2Ex.tex}}
\end{minipage}
\begin{minipage}[c]{0.24\textwidth}
\scalebox{0.7}{\input{figures_gen/QuiverG2Ex.tex}}
\end{minipage}
\\
\hspace{-1cm}
\begin{minipage}[t]{0.25\textwidth}
\caption{$Q_{A_2}$.}\label{fig:QuiverA2}
\end{minipage}
\begin{minipage}[t]{0.25\textwidth}
\caption{$Q_{B_2}$.}\label{fig:QuiverB2}
\end{minipage}
\begin{minipage}[t]{0.25\textwidth}
\caption{$Q_{C_2}$.}\label{fig:QuiverC2}
\end{minipage}
\begin{minipage}[t]{0.25\textwidth}
\caption{$Q_{G_2}$.}\label{fig:QuiverG2}
\end{minipage}
\end{center}
\end{figure}

Every quiver $Q$ has an associated \emph{seed torus} $T_Q$ (see Definition~\ref{def:SeedTorus}), which is a complex torus with a coordinate for each vertex. Mutation (see Definition~\ref{def:Mutation}) in a vertex $v_k$ of $Q_G$ gives rise to another quiver $\mu_{v_k}(Q_G)$ together with a birational map of seed tori $\mu_{v_k}\colon T_{Q_G}\to T_{\mu_{v_k}(Q_G)}$.
For a sequence $(i_1,\dots,i_k)$ of vertex indices define 
\begin{equation}\mu_{(i_1,\dots,i_k)}=\mu_{v_{i_k}}\mu_{v_{i_{k-1}}}\dots\mu_{v_{i_1}}. 
\end{equation}

\subsection{Rotations}
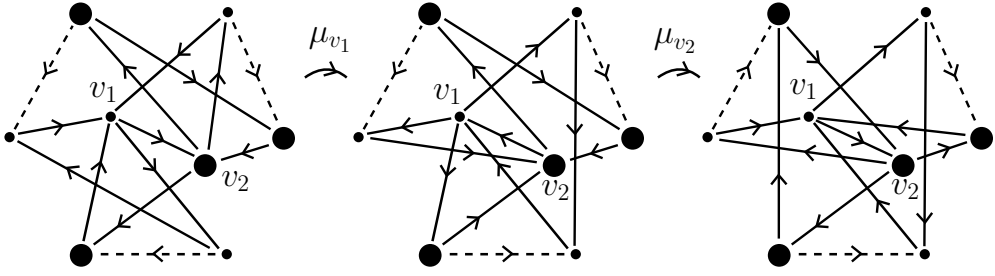
\begin{figure}[htb]
\begin{center}
\scalebox{0.75}{\input{figures_gen/MutationExample.tex}}
\end{center}
\caption{The mutation $\mu_{B_2}^{\rot}$ corresponds to a rotation by 120 degrees (after rearranging the position of $v_1$ and $v_2$).}\label{fig:MutationExample}
\end{figure}

Let
\begin{equation}
\mu_{A_2}^{\rot}=\id,\qquad \mu_{B_2}^{\rot}=\mu_{C_2}^{\rot}=\mu_{(1,2)},\qquad \mu_{G_2}^{\rot}=\mu_{(1,2,3,1,4,2)}.
\end{equation}
The following is a simple verification, which is illustrated in Figure~\ref{fig:MutationExample} for $G=B_2$.
\begin{lemma}\label{lemma:RotLemma} The quiver $\mu_G^{\rot}(Q_G)$ is isomorphic to $Q_G$ via an isomorphism which corresponds to a rotation by 120 degrees.
\end{lemma}

\begin{theorem}\label{thm:Rotation} There is a canonical birational equivalence
\begin{equation}
\M\colon \Conf_3^*(\A_G)\to T_{Q_G}
\end{equation}
such that the map $(g_0N,g_1N,g_2N)\to (g_2N,g_0N,g_1N)$ corresponds to the mutation sequence $\mu_G$ under the isomorphism $T_{\mu_G^{\rot}(Q_G)}\cong T_{Q_G}$ induced by Lemma~\ref{lemma:RotLemma}.
\end{theorem}
\begin{remark}\label{rm:M}
The map $\M$ is the composition of a birational equivalence $\Delta\colon\Conf_3^*(\A_G)\to T_{Q_G}$ given by minor coordinates (see Proposition~\ref{prop:EdgeFaceMinors}) and a monomial map $m\colon T_{Q_G}\to T_{Q_G}$ (see Section~\ref{sec:MonomialMap}).
\end{remark}

\subsection{The flip}
The generic configurations form an incomplete simplicial set with face maps
\begin{equation}\label{eq:SimpFaceMap}
\varepsilon_{i}\colon \Conf_k^*(\A)\to\Conf_{k-1}^*(\A),\qquad (g_0N,\dots,g_{k-1}N)\mapsto (g_0N,\dots,\widehat{g_i N},\dots,g_{k-1}N).
\end{equation}
For $i=0,\dots,k$, let $\kappa_i$ denote the map $\Conf_k^*(\A)\to\Conf_k^*(\A)$ which replaces the coset $g_iN$ by $g_is_GN$ leaving all other cosets fixed.

\subsubsection{Gluing configurations} We now consider configuration spaces $\Conf_3^*(\A)\times_{02}^{s_G}\Conf_3^*(\A)$ and $\Conf_3^*(\A)\times_{13}^{s_G}\Conf_3^*(\A)$ obtained by gluing together copies of $\Conf_3^*(\A)$ together along $\Conf_2^*(\A)$. Each is birationally equivalent to $\Conf_4^*(\A)$ and is defined by the pushout diagram
\begin{equation}\label{eq:Pushout}
\cxymatrix{{&\Conf_3^*(\A)\ar^-{\varepsilon_{j-1}\circ \kappa_{l-1}}[dr]&\\\Conf_4^*(\A)\ar^-{\varepsilon_i}[ur]\ar_-{\varepsilon_j\circ \kappa_k}[dr]\ar@{-->}^-{\Psi_{kl}}[r]&\Conf_3^*(\A)\times_{kl}^{s_G}\Conf_3^*(\A)\ar@{-->}[r]\ar[u]\ar[d]&\Conf_2^*(\A),\\&\Conf_3^*(\A)\ar_-{\varepsilon_i}[ru]&}}
\end{equation}
where $(i,j,k,l)$ is either $(0,2,1,3)$ or $(1,3,0,2)$. The map $\Psi_{kl}$ is illustrated in Figure~\ref{fig:MapConf4}.

\begin{figure}[htb]
\begin{center}
\scalebox{0.8}{\input{figures_gen/MapConf4.tex}}
\caption{Element in $\Conf_4^*(\A)$ and its image in $\Conf_3^*(\A)\times_{02}^{s_G}\Conf_3^*(\A)$ and $\Conf_3^*(\A)\times_{13}^{s_G}\Conf_3^*(\A)$.}\label{fig:MapConf4}
\end{center}
\end{figure}

\subsubsection{Gluing quivers} Similar to the gluing of configurations, we can glue together copies of $Q_G$ to form quivers $Q_G\cup_{02}Q_G$ and $Q_G\cup_{13} Q_G$. The formal construction is described in Section~\ref{sec:GluingQuivers}. We denote the (non-frozen) vertices of the ``right copy'' of $Q_G$ in $Q_G\cup_{02}Q_G$ by $\bar v_i$, and those of the  ``left copy'' by $v_i$. Similarly, we use $v_i$ for the ``top copy'' of $Q_G$ in $Q_G\cup_{13}Q_G$ and $\bar v_i$ for the bottom copy. The two vertices on the common edge are indexed by $0$ and $\infty$ (see Figures~\ref{fig:QC202} and~\ref{fig:QC213}). Let
\begin{equation}
\begin{gathered}
\mu_{A_2}^{\flip}=\mu_{(0,\infty,1,\bar 1)},\qquad \mu_{B_2}^{\flip}=\mu_{C_2}^{\flip}=
\mu_{(0,\infty,1,2,\bar 1,\bar 2,0,1,\bar 1)}\\
\mu_{G_2}^{\flip}=\mu_{(0,\infty,3,2,1,3,4,2,\bar 1,\bar 2,\bar 4,0,3,1,4,\bar 3,\bar 1,\bar 2,0,3,\bar 3,\bar 2,\bar 1,\bar 3)}.
\end{gathered}
\end{equation}

The following is a simple verification, which is illustrated in Figure~\ref{fig:FlipC2} for $G=C_2$.
\begin{lemma}\label{lem:FlipLemma}
There is a canonical isomorphism between $\mu_G^{\flip}(Q_G\cup_{02}Q_G)$ and $Q_G\cup_{13}Q_G$.
\end{lemma}

We may thus identify the seed tori $T_{\mu_G^{\flip}(Q_G\cup_{02}Q_G)}$ and $T_{Q_G\cup_{13}Q_G}$.

\begin{remark}
For verification of Lemmas~\ref{lem:FlipLemma} and \ref{lemma:RotLemma}, the java applet~\cite{KellerApp} by Mark Keller is very useful.
\end{remark}

\begin{figure}[htb]
\begin{center}
\begin{minipage}[c]{0.48\textwidth}
\scalebox{0.75}{\input{figures_gen/QC202.tex}}
\end{minipage}
\begin{minipage}[c]{0.47\textwidth}
\scalebox{0.75}{\input{figures_gen/QC213.tex}}
\end{minipage}
\\
\begin{minipage}[t]{0.48\textwidth}
\caption{The quiver $Q_{C_2}\cup_{02}Q_{C_2}$.}\label{fig:QC202}
\end{minipage}
\begin{minipage}[t]{0.47\textwidth}
\caption{The quiver $Q_{C_2}\cup_{13}Q_{C_2}$.}\label{fig:QC213}
\end{minipage}
\end{center}
\end{figure}


\begin{theorem}\label{thm:Flip}
We have a commutative diagram
\begin{equation}
\cxymatrix{{\Conf_3^*(\A)\times_{02}^{s_G}\Conf_3^*(\A)\ar[r]^-\cong\ar[d]_-{\Psi_{13}\Psi_{02}^{-1}}&T_{Q_G\cup_{02} Q_G}\ar[d]^-{\mu_G^{\flip}}\\\Conf_3^*(\A)\times_{13}^{s_G}\Conf_3^*(\A)\ar[r]^-\cong&T_{Q_G\cup_{13} Q_G,}}}
\end{equation}
where the vertical maps are induced by the map $\M$ in Theorem~\ref{thm:Rotation}.
\end{theorem}

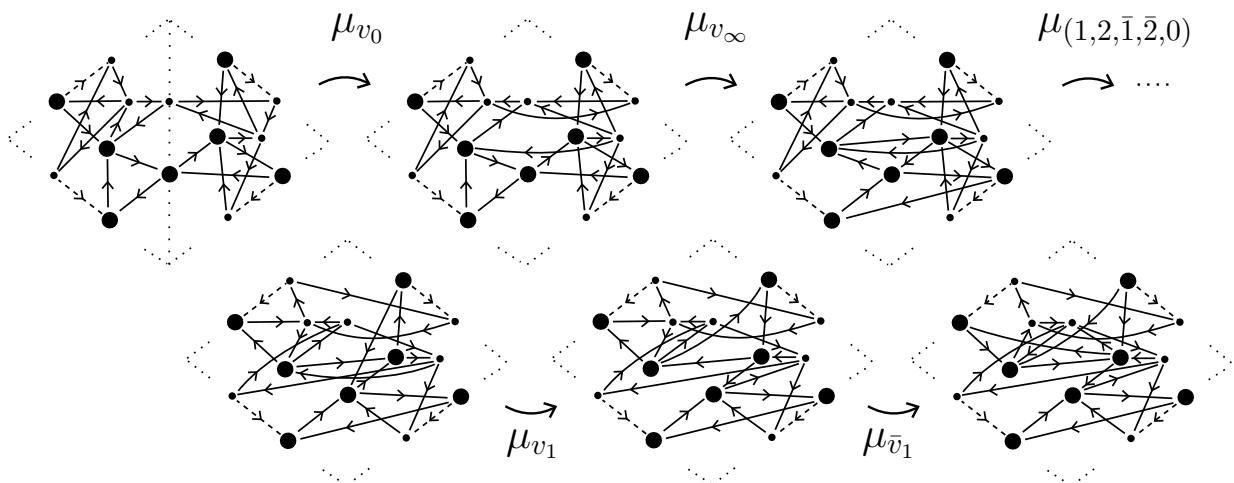
\begin{figure}[htb]
\begin{center}
\hspace{-0.5cm}\scalebox{1.12}{\input{figures_gen/FlipC2.tex}}
\end{center}
\caption{First and last two mutations in $\mu_{C_2}^{\flip}(Q_{C_2})$. The final configuration corresponds to $Q_G\cup_{13} Q_G$ after rearranging the vertices inside the dotted square.}\label{fig:FlipC2}
\end{figure}

\begin{conjecture}\label{conj:GeneralCase} For any semisimple, simply connected, complex Lie group $G$, there exists a quiver $Q_G$ and quiver mutations $\mu_G^{\rot}$ and $\mu_G^{\flip}$ such that Theorems~\ref{thm:Rotation} and \ref{thm:Flip} hold. The map $\M$ should be a composition of minor coordinates and a monomial map (see Remark~\ref{rm:M}).
\end{conjecture}
\begin{remark}
The minor coordinates in Conjecture~\ref{conj:GeneralCase} depend on a choice of reduced word for the longest element in the Weyl group. Hence, the quiver in Conjecture~\ref{conj:GeneralCase} should as well. In rank $2$ there are only two reduced words, and we have selected the one starting with $1$.
\end{remark}

\begin{remark} After finishing this paper I became aware of the very interesting recent preprint~\cite{ClusterStructuresLe} by Ian Le obtaining similar results for the groups $B_n$, $C_n$ and $D_n$ using different methods.
\end{remark}

\subsection{An atlas on $\A_{G,S}$}\label{sec:CoordinatesAGS}
Let $S$ be an oriented punctured surface with negative Euler characteristic. The universal cover of $S$ is an open oriented disk $D$, and the lift of the punctures define a countable, cyclically ordered $\pi_1(S)$-set $F_\infty(S)$ of points on $\partial\overline D$, the \emph{Farey set}~\cite[Sec.~1.3]{FockGoncharov}. Let $\overline F_\infty(S)$ be the double cover of $F_\infty(S)$ induced by the double cover of $\partial\overline D$, and let $\sigma$ denote the non-trivial automorphism. The fundamental group of the punctured tangent bundle is a $\Z$-extension of $\pi_1(S)$, and the quotient by $2\Z$ is a central $\Z/2\Z$-extension $\overline\pi_1(S)$ (see~\cite[Sec.~2.4]{FockGoncharov}). Let $\overline\sigma$ denote the generator. 

The space $\A_{G,S}$ is the moduli space of decorated twisted $G$-local systems on $S$ (\cite[Def~2.4]{FockGoncharov}).
When $s_G$ is trivial we may regard it as the quotient stack of pairs $(\rho,D)$ by the $G$-action $g(\rho,D)=(g\rho g^{-1},gD)$, where $\rho\colon \pi_1(S)\to G$ is boundary-unipotent (loops encircling punctures map to conjugates of $N$), and $D\colon F_\infty(S)\to \A$ is a $\rho$-equivariant map. When $s_G$ is non-trivial, it is the quotient stack of pairs $(\overline\rho,\overline D)$, where $\overline\rho\colon\overline\pi_1(S)\to G$ is a boundary-unipotent representation taking $\overline\sigma$ to $s_G$, and $\overline D\colon\overline F_\infty(S)\to\A$ is $\overline \rho$-equivariant (see~\cite[Sec.~8.6]{FockGoncharov}). 

Given a topological ideal triangulation $\T$ of $S$ we get an atlas on $\A_{G,S}$ as in~\cite[Sec.~8]{FockGoncharov}. The process is illustrated in Figures~\ref{fig:FundPoly} and \ref{fig:FundPolyB2} in the case when $S$ is a twice punctured torus and $G=B_2$. Pick a fundamental polyhedron $P$ for $\T$ in $D$. The triangulation of $S$ induces a triangulation of $P$. Pick an ordering $O$ of the vertices of $P$ agreeing with the cyclic ordering on $F_\infty(S)$. This associates a copy of the quiver $Q_G$ to each triangle, and by gluing these together, we obtain a quiver whose seed torus embeds in $\A_{G,S}$. Note that if two edges in $P$ are identified, the corresponding coordinates are identified as well. By~\cite[Thm.~8.2]{FockGoncharov} this provides a positive atlas with a chart for each triple $(\T,P,O)$. Our main results give explicit formulas for how the coordinates change when changing the triple.
The pair $(\rho,D)$, or $(\overline\rho,\overline D)$, corresponding to a collection of coordinates can be explicitly computed using a \emph{natural cocycle} (see Section~\ref{sec:NaturalCocycle}).



\begin{figure}[htb]
\begin{center}
\begin{minipage}[c]{0.48\textwidth}
\begin{center}
\scalebox{0.75}{\input{figures_gen/FundPoly.tex}}
\end{center}
\hfill
\end{minipage}
\begin{minipage}[c]{0.47\textwidth}
\begin{center}
\scalebox{0.75}{\input{figures_gen/FundPolyB2.tex}}
\end{center}
\end{minipage}
\\
\begin{minipage}[t]{0.55\textwidth}
\caption{Fundamental polyhedron for the twice punctured torus $S$.}\label{fig:FundPoly}
\end{minipage}
\begin{minipage}[t]{0.4\textwidth}
\caption{Coordinates on $A_{B_2,S}$.}\label{fig:FundPolyB2}
\end{minipage}
\end{center}
\end{figure}

\subsection{$3$-manifold groups and Ptolemy varieties}\label{sec:PtolemyIntro}
Let $M$ be a compact 3-manifold with a topological ideal triangulation $\T$.  A representation $\pi_1(M)\to G$ is \emph{boundary-unipotent} if peripheral subgroups map to conjugates of $N$, and a \emph{decoration} of a boundary-unipotent representation is a $\rho$-equivariant assignment of an $N$-coset to each ideal point in the universal cover of $M$. In Section~\ref{sec:PtolemyVariety} we define a variety $P_G(\T)$ by gluing together configurations spaces $\Conf_4^*(\A)$ using a gluing pattern determined by the triangulation. Most of the results in~\cite{GaroufalidisThurstonZickert} on Ptolemy varieties for $\SL(n,\C)$ have natural analogues for $G$. As in~\cite[(9.26)]{GaroufalidisThurstonZickert} there is a natural one-to-one correspondence
\begin{equation}\label{eq:OneToOnePG}
\xymatrix{\left\{\txt{Points in\\$P_{G}(\T)$}\right\}\ar@{<->}[r]^-{1:1}&\left\{\txt{Generically decorated, boundary-unipotent\\$\pi_1(M)\to G$}\right\}\Big/G}
\end{equation}
and our main results yield an explicit formula for this map.

By a result of Kostant~\cite{Kostant} there is a canonical homomorphism $\SL(2,\C)\to G$, which preserves unipotent elements and takes $s_{\SL(2,\C)}$ to $s_G$. If $M=\H^3/\Gamma$ is a cusped hyperbolic 3-manifold, there is thus a canonical boundary-unipotent representation $\rho_G\colon\pi_1(M)\to G\big/\langle s_G\rangle$. As explained e.g.~in~\cite{ZickertDuke,GaroufalidisThurstonZickert}, $\rho_G$ need not have a boundary-unipotent lift to $G$, and the obstruction to the existence of such a lift is a class in $H^2(M,\partial M;\Z/2\Z)$. For each $\sigma\in H^2(M,\partial M;\Z/2\Z)$ there is variety $P_G^\sigma(\T)$, and the analogue of~\eqref{eq:OneToOnePG} is (c.f.~\cite[(9.31)]{GaroufalidisThurstonZickert})
\begin{equation}\label{eq:OneToOnePGsigma}
\xymatrix{\left\{\txt{Points in\\$P_G^\sigma(\T)$}\right\}\ar[r]^-{z:1}&\left\{\txt{Generically decorated, boundary-unipotent\\
$\pi_1(M)\to G\big/\langle s_G\rangle$ with obstruction class $\sigma$}\right\}\Big/G,}
\end{equation}
where $z$ is the order of the group $Z^1(M,\partial M;\Z/2\Z)$ of $\Z/2\Z$ valued $1$-cocycles (with cell structure induced by $\T$).

As in~\cite[Sec.~4.1]{GaroufalidisThurstonZickert}, there is a natural action of $H^c$ on $P_G(\T)$ and $P_G^\sigma(\T)$, where $H$ is the maximal torus in $G$ and $c$ is the number of boundary components of $M$. The action is monomial and the quotients are denoted by $P_G(\T)_{\red}$ and $P_G^\sigma(\T)_{\red}$. The maps~\eqref{eq:OneToOnePG} and \eqref{eq:OneToOnePGsigma} induce maps 
\begin{equation}
\xymatrixcolsep{1pc}\xymatrix{P_{G}(\T)_{\red}\ar[r]&\left\{\txt{Boundary-unipotent\\$\pi_1(M)\to G$}\right\}\Big/G},\quad
\xymatrixcolsep{1pc}\xymatrix{P_{G}^\sigma(\T)_{\red}\ar[r]&\left\{\txt{Boundary-unipotent\\$\pi_1(M)\to G$,\\obstruction class $\sigma$}\right\}\Big/G}
\end{equation}
which are generically $1:1$ and $\vert H^1(M,\partial M;\Z/2\Z)\vert:1$, respectively.

\subsection{Computations for the figure-eight knot complement}
Let $G=\Sp(4,\C)$, and let $M$ be the figure-eight knot complement with the standard ideal triangulation $\T$ with $2$ simplices.
The fundamental group of $M$ has a presentation
\begin{equation}\label{eq:TwoBridge}
\pi_1(M)=\langle x_1,x_2\bigm\vert x_1w=wx_2,w=x_2x_1^{-1}x_2^{-1}x_1\rangle.
\end{equation}
In Section~\ref{sec:Fig8} we show that $P_G(\T)_{\red}$ is empty and that $P_G^{\sigma}(\T)_{\red}$ consists of two zero-dimensional components of degree 2 and 6, respectively.
The component of degree $2$ is defined over $\Q(\sqrt{-3})$, and the corresponding representation in $\Sp(4,\C)/\langle-I\rangle$ takes $x_1$ and $x_2$ to
\begin{equation}
\scalebox{1}{$
\begin{bmatrix}
 1 & -\frac{9(1+\sqrt{-3})}{8} & \frac{3(-1+\sqrt{-3})}{4} & 1+\sqrt{-3} \\
 0 & 1 & -1-\sqrt{-3} & -\frac{16}{9} \\
 0 & 0 & 1 & 0 \\
 0 & 0 & \frac{9\left(1+\sqrt{-3}\right)}{8}  & 1 \\
\end{bmatrix}
$},\enspace
\scalebox{0.89}{$
\begin{bmatrix}
 0 & 0 & \frac{3\left(-1+\sqrt{-3}\right)}{4}  & 0 \\
 0 & 0 & -2+2\sqrt{-3} & \frac{4\left(-1+\sqrt{-3}\right)}{9}  \\
 \frac{1}{3}+\frac{\sqrt{-3}}{3} & -\frac{3\left(1+\sqrt{-3}\right)}{2}  & 8 & \frac{16}{3} \\
 0 & \frac{9}{16} \left(1+\sqrt{-3}\right) & -9 & -4
\end{bmatrix}
$},
\end{equation}
respectively. The component of degree $6$ is defined over $\Q(\omega)$, where 
\begin{equation}\label{eq:Omega}
\omega^6-\omega^5+3\omega^4-5\omega^3+8\omega^2-6\omega+8=0,
\end{equation}
and the corresponding representation is given by
\begin{equation}
x_1\mapsto
\begin{bmatrix}
 1 & a_2 & b_1 & b_2 \\
 0 & 1 & b_3 & b_4 \\
 0 & 0 & 1 & 0 \\
 0 & 0 & c_3 & 1
\end{bmatrix},
\qquad
x_2\mapsto
\begin{bmatrix}
 0 & 0 & b_1' & 0 \\
 0 & 0 & b_3' & b_4' \\
 c_1' & c_2' & d_1' & d_2' \\
 0 & c_4' & d_3' & d_4'
\end{bmatrix},
\end{equation}
where
\begin{equation}
\scalebox{1}{$
\begin{gathered}
a_2=-c_3=\frac{\omega ^5}{16}+\frac{7 \omega ^3}{16}-\frac{5 \omega ^2}{8}-\frac{5 \omega }{8}-\frac{3}{2},\quad b_1=-\frac{\omega ^5}{8}+\frac{\omega ^3}{8}+\frac{\omega ^2}{4}-\frac{3 \omega }{4}-1,\\
b_2=-b_3=-\frac{\omega ^5}{32}+\frac{\omega ^4}{16}+\frac{3 \omega ^3}{32}+\frac{\omega ^2}{16}-\frac{\omega
   }{16}-\frac{3}{4},\quad b_4=2\frac{b_2}{a_2},\\
c_1'=-b_1'^{-1}=\frac{3 \omega ^5}{32}-\frac{3 \omega ^4}{16}+\frac{7 \omega ^3}{32}-\frac{11 \omega ^2}{16}+\frac{11 \omega
   }{16}-\frac{1}{4},\\
c_2'=-\frac{\omega ^5}{4}+\frac{\omega ^4}{2}-\frac{3 \omega ^3}{4}+\omega ^2-\frac{5 \omega }{2}+3,\quad
c_4'=-b_2'^{-1}=-\frac{\omega ^5}{4}+\frac{\omega ^4}{2}-\frac{3 \omega ^3}{4}+\omega ^2-\frac{5 \omega }{2}+3,\\
b_3'=\frac{3 \omega ^5}{8}-\frac{\omega ^4}{2}+\frac{3 \omega ^3}{8}-\frac{3 \omega ^2}{2}+\frac{9 \omega
   }{4}-\frac{1}{2},\quad
b_4'=\frac{b_3'}{b_1'c_2'},\\
d_1'=\frac{\omega ^5}{8}+\frac{3 \omega ^4}{4}-\frac{3 \omega ^3}{8}+\frac{3 \omega ^2}{4}-\frac{7 \omega }{4}+7,
\quad d_2'=\frac{b_4'd_3'-b_3'd_4'}{b_1'},\\
d_3'=-\omega ^5+\omega ^4-\frac{3 \omega ^3}{2}+\frac{7 \omega ^2}{2}-6 \omega -1,\quad 
d_4'=-\frac{\omega ^5}{8}-\frac{3 \omega ^4}{4}+\frac{3 \omega ^3}{8}-\frac{3 \omega ^2}{4}+\frac{7 \omega }{4}-3.
\end{gathered}$}
\end{equation}
These representations all lift to representations in $\Sp(4,\C)$, but no lift is boundary-unipotent.

\begin{remark}
We stress that the notion of genericity depends on the triangulation. There may be more representations than those detected by the Ptolemy variety.
A triangulation independent Ptolemy variety detecting all irreducible representations is defined for $G=\SL(2,\C)$ in~\cite{IndependentPtolemy}.
\end{remark}

\section{Quivers, seed tori, and mutations}\label{sec:Quivers}
The following definition of a (weighted) quiver serves our needs. The definition is a special case of the notion of a \emph{seed} as defined by Fock and Goncharov~\cite{FGClusterEnsembles}; see Remark~\ref{rm:Seed}. For closely related notions see e.g.~\cite{Marsh,Keller}. 

\begin{definition}\label{def:Quiver} Let $m\geq 1$ be an integer. A \emph{quiver} (of weight $m$) is a directed graph without $2$-cycles together with a partition of the vertices and edges into two types; \emph{fat vertices} (of weight $m$) or not, and \emph{half-edges or not}, respectively. In the case when $m=1$ we do not distinguish between vertices. All edges joining two vertices are required to have the same type, and the multiplicity of a half-edge must be odd. An \emph{isomorphism} of quivers is an isomorphism of graphs preserving the types of edges and vertices.  
\end{definition}

\begin{example}The graphs in Figures~\ref{fig:QuiverA2}, \ref{fig:QuiverB2}, \ref{fig:QuiverC2}, and \ref{fig:QuiverG2} define quivers $Q_{A_2}$, $Q_{B_2}$, $Q_{C_2}$ and $Q_{G_2}$, and we declare the weights to be $1$, $2$, $2$, and $3$, respectively.\end{example}

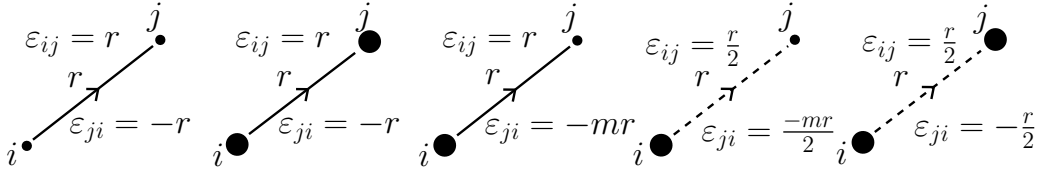
\begin{figure}[htb]
\begin{center}
\scalebox{0.75}{\input{figures_gen/Epsilonijs.tex}}
\caption{Definition of $\varepsilon_{ij}$ when $i$ and $j$ are joined by an edge of multiplicity $r$.}\label{fig:Epsilonijs}
\end{center}
\end{figure}

For a quiver $Q$ let $V_Q$ denote the set of vertices. When $V_Q=\{v_i\}_{i\in I}$, we shall denote a vertex either by $v_i$ or simply by $i$. For vertices $i$, and $j$, let $\sigma_{ij}$ denote the number of directed edges from $i$ to $j$ counting a half-edge as $1/2$, and counting an edge from $j$ to $i$ negative. A quiver determines a pair of functions
\begin{equation}\label{eq:dande}
\begin{aligned}
d_Q\colon V_Q&\to\{1,m\},\quad i\mapsto d_i,&\qquad \varepsilon_Q\colon V_Q\times V_Q&\to\frac{1}{2}\Z,\quad (i,j)\mapsto\varepsilon_{ij}\\
d_i&=\begin{cases}m&\text{ if $i$ is fat}\\1& \text{ otherwise}\end{cases},&\qquad \varepsilon_{ij}&=\frac{d_j}{\gcd(d_i,d_j)}\sigma_{ij}.
\end{aligned}
\end{equation}
The $\varepsilon_{ij}$ are illustrated in Figure~\ref{fig:Epsilonijs}.

\begin{lemma} For any set $V$ and functions $d\colon V\to\{1,m\}$ and $\varepsilon\colon V\times V\to\frac{1}{2}\Z$ such that $\varepsilon_{ij}/d_j=-\varepsilon_{ji}/d_i\in\Q$ for all $(i,j)\in V\times V$, there is a unique quiver $Q$ with vertex set $V$ satisfying that $d_Q=d$ and $\varepsilon_Q=\varepsilon$.
\end{lemma}
\begin{proof}
By~\eqref{eq:dande}, $d$ determines which vertices are fat, and $\varepsilon$ determines the multiplicity of an edge. An edge is a half-edge if and only if either $\varepsilon_{ij}$ or $\varepsilon_{ji}$ is a half-integer.
\end{proof}

\begin{definition} A vertex of a quiver is \emph{frozen} if it lies on a half-edge. The set of frozen vertices is denoted by $V^0_Q$.
\end{definition}
For the quivers $Q_G$ we index the six frozen vertices by pairs $ij$ as shown in Figure~\ref{fig:ExampleB2}.

\begin{remark}\label{rm:Seed} The tuple $(V_Q,V_Q^0,\varepsilon_Q,d_Q)$ is a \emph{seed} in the sense of~\cite[Def.~1.6]{FGClusterEnsembles}.
\end{remark}

\subsection{Quiver mutations and seed tori}
\begin{figure}[htb]
\begin{center}
\scalebox{0.75}{\input{figures_gen/MutationRules.tex}}
\end{center}
\caption{Mutation of the graph. The integer $\alpha_{ijk}$ is $m$ if $d_i=d_j\neq d_k$ and $1$ otherwise.}\label{fig:MutationRules}
\end{figure}
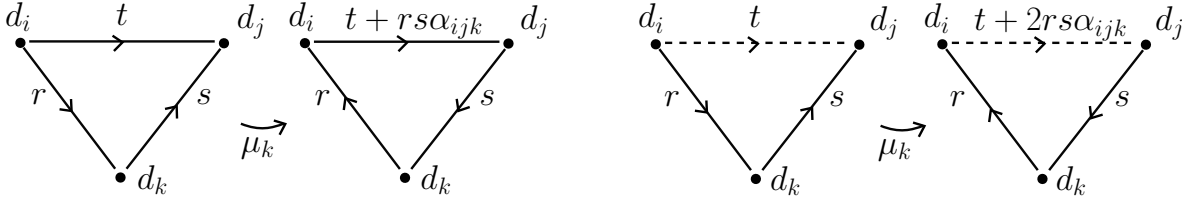

A process called~\emph{mutation} transforms one quiver to another. We follow Fock and Goncharov~\cite{FGClusterEnsembles}. 
\begin{definition}\label{def:Mutation} Let $Q$ be a quiver and $k$ a non-frozen vertex. Let $\mu_k(Q)$ be the unique quiver with $V_{\mu_k(Q)}=V_Q$, $d_{\mu_k(Q)}=d_Q$, and
\begin{equation}\label{eq:EpsilonMutation}
\varepsilon_{\mu_k(Q)}(i,j)=
\begin{cases}-\varepsilon_{ij}&\text{if } k\in\{i,j\}\\
\varepsilon_{ij}&\text{if } \varepsilon_{ik}\varepsilon_{kj}\leq 0,\enspace k\notin\{i,j\}\\
\varepsilon_{ij}+\vert \varepsilon_{ik}\vert\varepsilon_{kj}&\text{if } \varepsilon_{ik}\varepsilon_{kj}>0,\enspace k\notin\{i,j\}.
\end{cases}
\end{equation}
We say that $\mu_k(Q)$ is obtained from $Q$ by a \emph{mutation} at $k$.
\end{definition}
Note that mutation is an involution, i.e., $\mu_k(\mu_k(Q))=Q$. The formula~\eqref{eq:EpsilonMutation} implies that a mutation transforms the graph as shown in Figure~\ref{fig:MutationRules}. We refer to the Figures~\ref{fig:MutationExample}, \ref{fig:FlipC2}, and \ref{fig:ExampleB2} for examples.

\begin{definition}\label{def:SeedTorus}
The \emph{seed torus} associated to a quiver $Q$ is the complex torus
\begin{equation}
T_Q=\Hom_\Z(\Lambda_Q,\C^*),
\end{equation}
where $\Lambda_Q$ is the free abelian group generated by $V_Q$. The natural identification of  $T_Q$ with $(\C^*)^{\vert V_Q\vert}$, endows $T_Q$ with a coordinate system $\{a_i\}_{i\in V_Q}$.
\end{definition}

A mutation induces a birational map of seed tori
\begin{equation}
\mu_k\colon T_Q\to T_{\mu_k(Q)}
\end{equation}
\begin{equation}
\mu_k^*(a_k')=\frac{1}{a_k}\big(\prod_{j|\varepsilon_{kj}>0}a_j^{\varepsilon_{kj}}+\prod_{j|\varepsilon_{kj}<0}a_j^{-\varepsilon_{kj}}\big),\qquad \mu_k^*(a_i')=a_i, \text{ for }i\neq k
\end{equation}
Since mutations are only allowed at non-frozen vertices, the coordinates of the frozen vertices always stay fixed.
 
\begin{example}\label{ex:MutB2Example}
For the mutation shown in Figure~\ref{fig:ExampleB2} we have 
\begin{equation}
\begin{gathered}
a_{ij}=a_{ij}'=a_{ij}'',\qquad a_1'=\frac{1}{a_1}\big(a_{01}a_{02}^2a_{12}+a_{20}a_2^2\big),\qquad a_2'=a_2\\
a_1''=a_1',\qquad a_2''=\frac{1}{a_2'}\big(a_1'a_{10}'+a_{01}'a_{21}'a_{02}'\big)=\frac{a_{01} a_{02}^2 a_{10 }a_{12} + a_{10} a_2^2 a_{20} + a_{01} a_{02} a_1 a_{21}}{a_1 a_2}.
\end{gathered}
\end{equation}

\begin{figure}[htb]
\begin{center}
\scalebox{0.75}{\input{figures_gen/ExampleB2.tex}}
\end{center}
\caption{Coordinates on $T_{Q_{B_2}}$, $T_{\mu_{v_1}(Q_{B_2})}$, and $T_{\mu_{v_2}\mu_{v_1}(Q_{B_2})}$.}\label{fig:ExampleB2}
\end{figure}
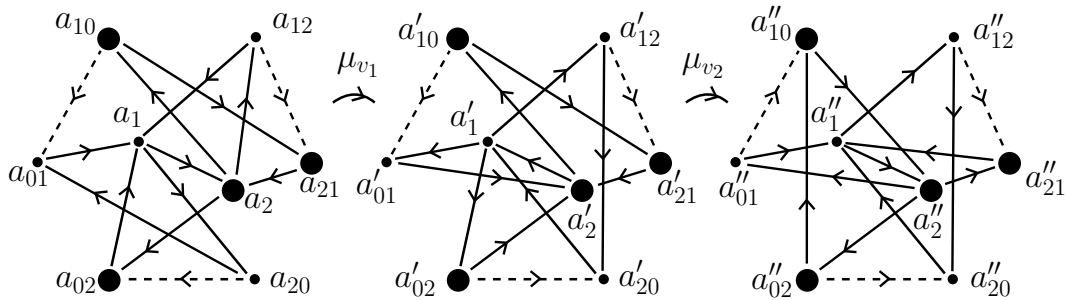
\end{example}

\subsection{Gluing quivers along frozen vertices}\label{sec:GluingQuivers}
For a subset $S$ of $V_Q$, let $Q_S$ denote the largest subgraph of $Q$ with vertex set $S$. It inherits a quiver structure from $Q$.
\begin{definition}
Let $Q$ and $Q'$ be quivers and let $W\subset V^0_{Q}$ and $W'\subset V^0_{Q'}$ be subsets of frozen vertices, and $\phi\colon Q_W\to Q_{W'}$ an isomorphism. The quiver $Q\cup_\phi Q'$ is the quiver obtained by gluing together $Q$ and $Q'$ via $\phi$, eliminating $1$ and $2$-cycles, and declaring that the gluing of two half-edges is a full edge (not a half-edge).
\end{definition}

Note that the vertices in $Q\cup_\phi Q'$ corresponding to $W$ and $W'$ are no longer frozen, and are thus open for mutation. Also note that the seed torus for $Q\cup_\phi Q'$ is the fiber product $T_Q\times_{T_{Q_W}}T_{Q'}$. 

\subsubsection{The quivers $Q_G\cup_{02}Q_G$ and $Q_G\cup_{13}Q_G$}
We now give a formal definition of the quivers $Q_G\cup_{02}Q_G$ and $Q_G\cup_{13}Q_G$ introduced in Section~\ref{sec:Results}. Recall that the frozen vertices of $Q_G$ are indexed by pairs $ij$ with $i,j\in\{0,1,2\}$. Let $\sigma_G\in S_2$ be the trivial permutation when $G=A_2$ and the non-trivial permutation otherwise (this is the permutation of the fundamental weights given by the longest element in the Weyl group, see Section~\ref{sec:WeylGroup}). Denote the non-frozen vertices of one copy of $Q_G$ by $\bar v_i$ and let
\begin{equation}
\phi_{02}\colon\{v_{01},v_{10}\}\to\{\bar v_{02},\bar v_{20}\},\qquad \phi_{13}\colon\{v_{02},v_{20}\}\to\{\bar v_{12},\bar v_{21}\},
\end{equation}
be such that $\phi_{02}$ takes the pair $(v_{01},v_{01})$ to $\sigma_G(\bar v_{02},\bar v_{20})$ and $\phi_{13}$ takes $(v_{02},v_{20})$ to $\sigma_G(\bar v_{12},\bar v_{21})$. We can now define
$Q_G\cup_{kl}Q_G$ to be $Q_G\cup_{\phi_{kl}}Q_G$.
We denote the images of $v_{01}$ and $v_{10}$ in $Q_G\cup_{02}Q_G$ by $v_\infty$ and $v_0$, respectively. Similarly, we denote the images $v_{02}$ and $v_{20}$ in $Q_G\cup_{13}Q_G$ by $v_\infty$ and $v_0$. The frozen vertices of $Q_G\cup_{kl} Q_G$ are indexed according to the edges of a quadrilateral (see Figures~\ref{fig:QC202} and \ref{fig:QC213}).

\subsection{Explicit formulas for mutations}\label{sec:ExplicitMutations}
An isomorphism of quivers induces an isomorphism of seed tori. In particular, by Lemma~\ref{lemma:RotLemma}, we may identify $T_{\mu^{\rot}(Q_G)}$ with $T_{Q_G}$, and the identification is such that $a_{ij}'=a_{i-1,j-1}$ (indices modulo $3$). 
\subsubsection{Formulas for $\mu_G^{\rot}$} Following Example~\ref{ex:MutB2Example} the explicit formulas for the non-frozen coordinates for $B_2$ and $C_2$ are given by
\begin{equation}\label{eq:MurotABC}
\begin{aligned}
&B_2:& a_1'&=\frac{1}{a_1}\big(a_{01}a_{02}^2a_{12}+a_{20}a_2^2\big),&\quad a_2'&=\frac{a_{01} a_{02}^2 a_{10 }a_{12} + a_{10} a_2^2 a_{20} + a_{01} a_{02} a_1 a_{21}}{a_1 a_2}\\
&C_2:& a_1'&=\frac{a_{01} a_{02} a_{12} + a_2 a_{20}}{a_1},&\quad a_2'&=
\frac{a_{10} (a_{01} a_{02} a_{12} + a_2 a_{20})^2 + a_{01}^2 a_{02} a_1^2 a_{21}}{a_1^2 a_2}.
\end{aligned}
\end{equation}
For $G_2$ the closed formula is rather lengthy so we instead introduce a ``dummy variable'' for each intermediate mutation.
\begin{equation}
\begin{aligned}
&G_2: &z_1{a_1}&=a_2 a_{20} + a_{01} a_{02} a_3,&z_2{a_2}&=a_{01}^3 a_{02}^2 a_4 + z_1^3, & a_1'{a_3}&=a_{20} a_4 + a_{12} z_1,\\
&&a_2'{a_4}&=a_{10} a_1'^3 + a_{21} z_2,&a_3'{z_1}&=a_{01}^2 a_{02} a_1' + z_2,&a_4'{z_2}&=a_{01}^3 a_{02} a_2' + a_{10} a_3'^3.\\
\end{aligned}
\end{equation}


\subsubsection{Formulas for $\mu_G^{\flip}$} As for $\mu_{G_2}^{\rot}$ we express the formulas for the non-frozen coordinates via dummy variables $z_i$. For $A_2$ the relations are \emph{Ptolemy relations}, i.e.~of the form $ef=ab+cd$.


\begin{equation}
\begin{aligned}
&A_2:& \bar a_1'{a_0}&=a_{01} a_1 + a_{03} \bar a_1,& a_1'{a_{\infty}}&=a_1 a_{21} + \bar a_1 a_{23},\\
&&a_0'{a_1}&= a_{30} a_1' + a_{32} \bar a_1',& a_{\infty}'{\bar a_1}&= a_{10} a_1' + a_{12} \bar a_1'.
\end{aligned}
\end{equation}


\begin{equation}
\begin{aligned}
&B_2:& z_1a_0&=a_2 \bar a_2 + \bar a_1 a_{32},& \bar a_1' a_{\infty}&= a_{01} a_1 + \bar a_1 a_{30}, &z_2a_1&= a_{03}^2 \bar a_1 a_{23} + a_2^2 \bar a_1',\\
&& \bar a_2'a_2& = a_{03} z_1 + z_2,&
 z_3\bar a_1&=z_1^2 + a_{12} z_2,&a_2'\bar a_2&= a_{10} a_{12} a_{32} + a_{21} z_1, \\
 &&a_{\infty}'{z_1}&= a_2' \bar a_2' + a_{10} z_3,& a_0' z_2&= a_{23} \bar a_2'^2 + \bar a_1' z_3,&a_1'z_3&= a_0' a_2'^2 + a_{12} a_{23} a_{\infty}'^2.
\end{aligned}
\end{equation}


\begin{equation}\label{eq:C2Flipzs}
\begin{aligned}
&C_2:& z_1{a_0} &= a_2 \bar a_2 + \bar a_1^2 a_{32},& \bar a_1' {a_{\infty}}&= a_{01} a_1 + \bar a_1 a_{30}, &z_2{a_1}& = a_{03} \bar a_1 a_{23} + a_2 \bar a_1',\\
&& \bar a_2'{a_2}&= a_{03} z_1 + z_2^2,&z_3{\bar a_1}& = z_1 + a_{12} z_2,& a_2' {\bar a_2}&= a_{10} a_{12}^2 a_{32} + a_{21} z_1,\\
&&a_{\infty}' {z_1}&= a_2' \bar a_2' + a_{10} z_3^2,& a_0'{z_2} &= a_{23} \bar a_2' + \bar a_1' z_3,& 
 a_1' {z_3}&=a_0' a_2' + a_{12} a_{23} a_{\infty}'.
\end{aligned}
\end{equation}


\begin{equation}
\begin{aligned}
&G_2:& 
z_1{a_0}&= a_{32} \bar a_3^3 + a_4 \bar a_4,& \bar a_1'a_{\infty} &= a_{01} a_1 + a_{30} \bar a_3,&z_2{a_3}&= a_2 a_{23} + a_1 a_4,\\
&&z_3{a_2}&= a_{03} a_4^2 + z_2^3,&z_4{a_1}& =a_{03} a_{23} \bar a_3 + \bar a_1' z_2,&
 z_5{z_2}&= \bar a_3 z_3 + a_4 z_4,\\
 &&z_6{a_4}&= z_1 z_3 + z_5^3, &\bar a_2'{z_3} &= z_4^3 + a_{03} z_6,&
 z_7{\bar a_3}&= z_1 + \bar a_1 z_5,\\
 &&z_8 {\bar a_4}&= \bar a_1^3 a_{32} + \bar a_2 z_1,&
a_4' {\bar a_2}&= a_{10} a_{12}^3 a_{32} + a_{21} z_8,&
 z_9{z_1}&= z_7^3 + z_6 z_8,\\
&&z_{10}{z_5}&=z_6 + z_4 z_7,& \bar a_3' {z_4}&= a_{23} \bar a_2' + \bar a_1' z_{10},&
\bar a_4' {z_6}&= z_{10}^3 + \bar a_2' z_9,\\
&&z_{11}{\bar a_1}&= a_{12} z_7 + z_8, &
z_{12}{z_7}&= z_{10} z_{11} + z_9, &z_{13}{z_8}&=a_{10} z_{11}^3 + a_4' z_9,\\
&&a_{\infty}'{z_9} &= a_{10} z_{12}^3 + \bar a_4' z_{13}, &
a_0'{z_{10}}&=a_{23} \bar a_4' + \bar a_3' z_{12},&
z_{14}{z_{11}}&= a_4' z_{12} + a_{12} z_{13},\\
&&a_2'{z_{13}}&= a_4'^2 a_{\infty}' + z_{14}^3,&
 a_1'{z_{12}}&= a_{12} a_{23} a_{\infty}' + a_0' z_{14},&
a_3'{z_{14}} &= a_{23} a_2' + a_1' a_4'.
\end{aligned}
\end{equation}



\section{Preliminaries on Lie groups}\label{sec:LieGroups}
Let $G$ be a simply connected, semisimple, complex Lie group of rank $r$ with Lie algebra $\g$. 
It is well known that $G$ is the $\C$ points of a linear algebraic group over $\Z$, and is thus an affine variety.

\subsection{Basic notions}
Fix a Cartan subalgebra $\h$ of $\g$, and a set $\Pi=\{\alpha_1,\dots,\alpha_r\}\subset \h^*$ of simple roots. This gives rise to a root space decomposition 
\begin{equation}
\g=\n_-\oplus\h\oplus\n,\qquad \n_-=\bigoplus_{\alpha\in\Delta_-}\g_\alpha, \qquad \n_+=\bigoplus_{\alpha\in\Delta_+}\g_\alpha,
\end{equation} 
where $\Delta_-$ and $\Delta_+$ denote the sets of negative, respectively, positive roots, and $\g_\alpha$ denotes the root space for a root $\alpha$.
Let $N_-$, $H$, and $N$ denote the Lie subgroups of $G$ with Lie algebras $\n_-$, $\h$, and $\n$, respectively.
Fix Serre generators $e_i\in \g_{\alpha_i}$, $f_i\in\g_{-\alpha_i}$, and $h_i\in\h$ of $\g$, and let
\begin{equation}
x_i(t)=\exp(te_i)\in N,\qquad y_i(t)=\exp(tf_i)\in N_-,\qquad t\in\C.
\end{equation}

\subsubsection{Fundamental weights and the Cartan matrix} Let $\langle,\rangle$ denote the symmetric bilinear form on $\h^*$ dual to the Killing form $B$ on $\h$. For each root $\alpha$, let $H_\alpha\in\h$ be the unique element satisfying that $\alpha(H)=B(H,H_\alpha)$, and let
\begin{equation}
\alpha^\vee=\frac{2}{\langle\alpha,\alpha\rangle}\alpha,\qquad h_\alpha=\frac{2}{\langle\alpha,\alpha\rangle}H_\alpha.
\end{equation}
The element $h_{\alpha_i}$ is the Serre generator $h_i$.
The set of $\gamma\in\h^*$ with $\gamma(h_i)\in\Z$ for all $i$ form a lattice $P$ generated by the \emph{fundamental weights}, which are the elements $\omega_1,\dots,\omega_r\in \h^*$ satisfying that $\omega_i(h_{\alpha_j})=\delta_{ij}$, or equivalently, that $\langle\omega_i,\alpha_j^\vee\rangle=\delta_{ij}$. The \emph{Cartan matrix} is the matrix $A$ with entries $A_{ij}=\langle\alpha_i^\vee,\alpha_j\rangle$.

\subsubsection{Coordinates on H} For every weight $\omega$, there is a \emph{character} $\chi_\omega\colon H\to\C^*$, and for every root $\alpha$ there is a \emph{cocharacter} $\chi^*_\alpha\colon \C^*\to H$. These are defined by
\begin{equation}
\chi_\omega(\exp(h))=e^{\omega(h)},\qquad \chi^*_\alpha(e^t)=\exp(h_\alpha t),
\end{equation} 
and satisfy
\begin{equation}
\chi_\omega\circ\chi^*_\alpha(t)=t^{\langle \omega,\alpha^\vee\rangle},\qquad t\in\C^*.
\end{equation}
It follows that we have an isomorphism
\begin{equation}\label{eq:CoordsOnH}
H\cong (\C^*)^r,\qquad h\mapsto\big(\chi_{\omega_1}(h),\dots,\chi_{\omega_r}(h)\big),\qquad \chi^*_{\alpha_1}(h_1)\cdots\chi^*_{\alpha_r}(h_r)\mapsfrom (h_1,\dots,h_r)
\end{equation} 
We may thus identify $H$ with $(\C^*)^r$. We sometimes denote $\chi^*_{\alpha_i}(t)$ by $h_i^t$.

\subsubsection{The Weyl group and reduced words}\label{sec:WeylGroup}
The \emph{Weyl group} $W$ is the group generated by the \emph{simple root reflections} $s_i$ given by
\begin{equation}\label{eq:RootReflection}
s_i(\gamma)=\gamma-\langle\gamma,\alpha_i^\vee\rangle\alpha_i,\qquad \gamma\in \h^*.
\end{equation}
The Weyl group is isomorphic to $N_G(H)/H$, and there is a section (see \cite[Sec.~1.4]{FominZelevinsky})
\begin{equation}\label{eq:WeylSplit}
\begin{gathered}
W\to N_G(H),\qquad w\mapsto \overline w\\
\overline{s_{i_1}\!\cdots s_{i_k}}\mapsto \overline{s_{i_1}}\cdots \overline{s_{i_k}},\qquad \overline {s_i}=x_i(-1)y_i(1)x_i(-1).
\end{gathered}
\end{equation}
The Weyl group is a Coxeter group and there is a unique longest element $w_0$. The Weyl group acts on $\h^*$ and permutes the simple roots and fundamental weights. It also acts on $H$ via~\eqref{eq:WeylSplit}, i.e.~$w(h)=\overline w h\overline w^{-1}$.  The action by $w_0$ is such that $w_0(\omega_i)=-\omega_{\sigma_G(i)}$ for a permutation $\sigma_G\in S_r$. In particular, if $h=(h_1,\dots,h_r)\in H$ we have
\begin{equation}\label{eq:PermuteCoordinates}
w_0(h)=(h_{\sigma_G(1)}^{-1},\dots,h_{\sigma_G(r)}^{-1}).
\end{equation}
Using the explicit root data given in Section~\ref{sec:Rank2} one checks that $\sigma_G\in S_2$ is trivial for $B_2$, $C_2$ and $G_2$ and non-trivial for $A_2$.

A \emph{reduced word} for $w\in W$ is a tuple $\mathbf i=(i_1,\dots,i_m)$, with $m$ minimal, such that
\begin{equation}
w_0=s_{i_1}\cdots s_{i_m}.
\end{equation} 
In all of the following we shall fix a reduced word $\mathbf i=(i_0,\dots,i_m)$ for $w_0$. The length $m$ is equal to the number of positive roots.


\subsection{The element $s_G$}\label{sec:sG}
Consider the element
\begin{equation}
s_G=\prod_{\alpha\in\Delta_+}\chi_\alpha^*(-1)\in H.
\end{equation}
As shown in~\cite[Sec.~2.3]{FockGoncharov}, $s_G$ is central in $G$ and has order dividing $2$, and $\overline{w_0}^{-1}=\overline{w_0}s_G$. Clearly, all coordinates of $s_G$ are either $1$ or $-1$.
\subsection{Chamber weights and generalized minors}
References for this section include \cite{BerensteinZelevinsky,FominZelevinsky,FockGoncharov}. We adopt the notation of~\cite{BerensteinZelevinsky,FominZelevinsky}, and warn the reader that the symbols $w$ and $\omega$ look very similar.

Let $G_0$ be the Zariski open subset of elements $g\in G$ admitting a (necessarily unique) factorization
\begin{equation}\label{eq:yhx}
g=[g]_-[g]_0[g_+], \qquad [g]_-\in N_-,\quad  [g]_0\in H, \quad [g]_+\in N.
\end{equation}
When writing $g=yhx$, we shall often implicitly assume that $y\in N_-$, $h\in H$, and $x\in N$. The factors $y$, $h$, and $x$ are regular functions of $g\in G_0$.
\begin{definition} A \emph{chamber weight} is an element $\gamma\in\h^*$ in the Weyl orbit of a fundamental weight, i.e.~$\gamma=w\omega_i$ for some $i\in\{1,\dots,r\}$ and $w\in W$. 
\end{definition}
\begin{definition} For a chamber weight $\gamma=w\omega_i$, the \emph{(generalized) minor} associated to $\gamma$ is the regular function $\Delta^\gamma\colon G\to\C$ whose restriction to $\overline wG_0$ is given by
\begin{equation}
\Delta^\gamma(g)=\chi_{\omega_i}([\overline w^{-1}g]_0)\in\C.
\end{equation}
\end{definition}
\begin{remark}\label{rm:SLnMinors} For $G=A_r$, $W=S_r$ and the Chamber weight for $\sigma \omega_i$ is the $i\times i$ minor with rows $\sigma(1),\dots,\sigma(i)$ and columns $1,\dots,i$.
\end{remark}
Recall that we have fixed a reduced word $\mathbf i=(i_1,\dots i_m)$ for $w_0$. 
\begin{definition} An \emph{$\mathbf i$-chamber weight} is a chamber weight of the form
\begin{equation}
w_k\omega_i, \quad w_k=s_{i_m}s_{i_{m-1}}\cdots s_{i_k}, \qquad i\in\{1,\dots,r\},\quad k\in\{1,\dots,m+1\}.
\end{equation}
The corresponding minor is called an \emph{$\mathbf i$-minor}.
\end{definition}
\begin{proposition}[{\cite[Prop.~2.9]{BerensteinZelevinsky}}]\label{prop:DistinctWeights}
There are $m+r$ distinct $\mathbf i$-chamber weights, the fundamental weights $\omega_i$ and the weights
\begin{equation}\label{eq:gammak}
\gamma_k=w_k\omega_{i_k}=s_{i_m}s_{i_{m-1}}\cdots s_{i_k}\omega_{i_k},\qquad k\in\{1,\dots,m\}.
\end{equation}
Moreover, all chamber weights $w_0\omega_i$ are $\mathbf i$-chamber weights.
\end{proposition}

\begin{definition} We call the minors $\Delta^{w_0\omega_i}$ and $\Delta^{\omega_i}$ \emph{edge minors} and the $m-r$ remaining $\mathbf i$-minors \emph{face minors}.\end{definition}

Note that the edge minors $\Delta^{w_i}$ are the coordinates on $H$ given in~\eqref{eq:CoordsOnH}.
\subsection{Some biregular isomorphisms}
The \emph{transpose map} (see e.g.~\cite{FominZelevinsky,BerensteinZelevinsky}) is the unique biregular antiautomorphism $\Psi\colon G\to G$ satisfying 
\begin{equation}
\Psi(x_i(t))=y_i(t),\qquad \Psi(h)=h,\qquad \Psi(y_i(t))=x_i(t),\qquad t\in\C, h\in H. 
\end{equation}
One has (see \cite[p.~55]{FockGoncharov})
\begin{equation}
\Psi(\wzero)=\wzero^{-1}=\wzero s_G.
\end{equation}

The varieties $N\cap G_0\wzero$ and $N_-\cap\wzero G_0$ will be of special significance. One easily checks that $\Psi$ restricts to a biregular isomorphism between them.
Fomin and Zelevinsky~\cite{FominZelevinsky} define a biregular isomorphism
\begin{equation}
\pi_-\colon N\cap G_0\overline{w_0}\to N_-\cap\overline{w_0}G_0,\qquad 
x\mapsto \overline{w_0}^{-1}[x\overline{w_0}^{-1}]_+\overline{w_0},\qquad [\overline{w_0}y]_+\mapsfrom y.
\end{equation}
Similarly, one has a biregular isomorphism (also considered in~\cite{Lusztig,FockGoncharov})
\begin{equation}\label{eq:PhiDef}
\Phi\colon N\cap G_0\wzero\to N_-\cap \wzero G_0,\qquad
x\mapsto [x\wzero]_-,\qquad \wzero [\wzero^{-1} y]_-\wzero^{-1}\mapsfrom y.
\end{equation}
Note that $\Phi$ is determined by the (equivalent) properties
\begin{equation}\label{eq:PhiProp}
x\wzero N=\Phi(x)[x\wzero]_0N,\qquad xN_-=\Phi(x)[x\wzero]_0\wzero s_GN_-,\quad x\in N\cap G_0\wzero,
\end{equation}
which allow one to write a coset $x\wzero hN$ as $ykN$ and vice versa (and similarly for $N_-$ cosets).
Each of the isomorphisms respects conjugation by elements $h\in H$, i.e., we have
\begin{equation}\label{eq:Conjugation}
\Psi(hxh^{-1})=h^{-1}\Psi(x)h,\qquad \Phi(hxh^{-1})=h\Phi(x)h^{-1}, \qquad \pi_-(hxh^{-1})=h\pi_-(x)h^{-1}.
\end{equation}

\subsection{Factorization coordinates}
Consider the map
\begin{equation}
x_{\mathbf i}\colon\C^m\to N,\qquad (t_1,\dots,t_m)\mapsto x_{i_1}(t_1)\cdots x_{i_m}(t_m).
\end{equation}
Theorem~\ref{thm:MinorsAndFactorization} below summarizes~\cite[Thm~2.19]{FominZelevinsky} and \cite[Thms.~1.4, 4.3]{BerensteinZelevinsky}.
\begin{theorem}\label{thm:MinorsAndFactorization} Let $x=x_{\mathbf i}(t_1,\dots,t_m)\in N\cap G_0\overline{w_0}$, and let $y=\pi_-(x)$. The $t_i$ and the $\mathbf i$-minors of $y$ are related by the monomial expressions
\begin{equation}\label{eq:MinorsFact} 
t_k=\frac{1}{\Delta^{w_k\omega_{i_k}}(y)\Delta^{w_{k+1}\omega_{i_k}}(y)}\prod_{j\neq i_k}\Delta^{w_k\omega_j}(y)^{-A_{j,i_k}},\qquad \Delta^{\gamma_k}(y)=\prod_{l\geq k} t_l^{\langle \gamma_k,(\alpha^{\mathbf i}_l)^\vee\rangle},
\end{equation}
where $\alpha_l^{\mathbf i}=w_{l+1}(\alpha_{i_l})$.
\end{theorem}

\begin{remark}\label{rm:WhichMinor}Every minor occurring is equal to either $\Delta^{\omega_i}$, or some $\Delta^{\gamma_k}$. This can be seen using that $s_i\omega_j=\omega_j$ for $i\neq j$ (see e.g.~\cite[(2.5)]{BerensteinZelevinsky}). For example, if $\mathbf i=(1,2,1,2,1,2)$, then $w_3\omega_2=s_2s_1s_2s_1\omega_2=s_2s_1s_2\omega_2=\gamma_4$.
\end{remark}

\begin{corollary} The variety $x_{\mathbf i}((\C^*)^m)\cap G_0\wzero$ is isomorphic to the Zariski open subset of $N_-\cap \wzero G_0$ of points where the $\mathbf i$-minors are non-zero.\qed
\end{corollary}

\begin{corollary}\label{cor:Nminus} The map $N_-\to \C^m$ taking $y$ to $(\Delta^{\gamma_1}(y),\dots,\Delta^{\gamma_m}(y))$ is a birational equivalence.\qed
\end{corollary}




\section{Configuration spaces of tuples}
Let $G$ be as in Section~\ref{sec:LieGroups}, i.e.~semisimple of rank $r$. Most of the results of this section can be found in Fock-Goncharov~\cite[Sec.~8]{FockGoncharov}. Since our notation differs slightly from that of Fock and Goncharov, we give complete proofs.
\begin{definition}\label{def:SufGen} A tuple $(g_0N,\dots,g_{k-1}N)\in\A^k$ is \emph{sufficiently generic} if 
\begin{equation}\label{eq:SuffGen}
g_i^{-1}g_j\in \overline{w_0}G_0, \qquad i\neq j\in\{0,\dots,k-1\}, 
\end{equation}
a condition, which is open, and independent of the choice of coset representatives. The subvariety of $\A^k$ of sufficiently generic tuples is denoted by $\A^{k,*}$, and the quotient of $\A^{k,*}$ by the diagonal left $G$ action is denoted by $\Conf^*_k(\A)$.
\end{definition}
It is convenient to view a tuple $(g_0N,\dots,g_{k-1}N)$ as an ordered $(k-1)$-simplex $\Delta^{k-1}$ together with a labeling of the $i$th vertex by $g_iN$. 
\subsection{The variety structure on $\Conf^*_k(\A)$}
For $k>2$, let $\mathcal W_k$ be the Zariski open subset of $(B_-\cap\wzero G_0)^{k-2}$ consisting of points $(a_2,\dots,a_{k-1})$ with $a_i^{-1}a_j\in \wzero G_0$ for $i\neq j$. Let $\mathcal W_2$ be a singleton.
\begin{proposition}\label{prop:Akstar} For $k>1$ we have an isomorphism of varieties
\begin{equation}\label{eq:Akstar}
G\times H\times \mathcal W_k\to \A^{k,*},\qquad  (g,h,a_2,\dots, a_j)\mapsto g(N,\wzero hN,a_2N,\dots,a_{k-1}N) 
\end{equation}
\end{proposition}
\begin{proof}
Let $\alpha=(g_0N,\dots,g_{k-1}N)\in \A^{k,*}$ with $g_i\in G$ fixed coset representatives. Since $g_i^{-1}g_j\in\wzero G_0$, we have factorizations $g_0^{-1}g_i=\wzero y_ih_ix_i$. Let
\begin{equation}
a_i=\wzero y_1^{-1}y_ih_i[\overline{w_0}y_1^{-1}y_ih_i]_+^{-1}\in B_-\cap\overline G_0,\qquad i=2,\dots,k-1.
\end{equation} 
The $a_i$ and $h_i$ are independent of the coset representatives $g_i$ and are regular functions of $\alpha$. Letting $g=g_0\wzero y_1\wzero^{-1}$, one has $\alpha=g(N,\wzero h_1N,a_2N,\dots,a_{k-1}N)$. This proves the result.
\end{proof}
\begin{corollary}
The quotient $\Conf_k^*(\A)=\A^{k,*}/G$ is a variety isomorphic to $H\times \mathcal W_k$.\qed
\end{corollary}
\begin{example}For $k=2$ and $3$, we have
\begin{equation}\label{eq:Conf2and3}
\begin{gathered}
H\cong \Conf_2^*(\A),\qquad H\times B_-\cap\wzero G_0\cong\Conf_3^*(\A).\\
h\mapsto (N,\wzero hN),\qquad (h,a)\mapsto(N,\wzero hN,aN).
\end{gathered}
\end{equation}
Note that $(g_0N,g_1N)\in \Conf_2^*(\A)$ corresponds to $[\wzero^{-1}g_0^{-1}g_1]_0\in H$.
\end{example}
\begin{definition}The representative of $\alpha\in\Conf_k(\A)$ of the form $(N,\wzero hN,a_1N,\dots ,a_{k-2}N)$ is called the \emph{canonical representative}.\end{definition}

\subsection{Edge coordinates}
We have regular maps
\begin{equation}\label{eq:piij}
\pi_{ij}\colon\Conf_k^*(\A)\to H, \qquad (g_0N,\dots,g_{k-1}N)\mapsto [\wzero^{-1}g_i^{-1}g_j]_0, \quad i\neq j.
\end{equation}
Note that under the isomorphism $H\cong\Conf_2^*(\A)$, $\pi_{ij}$ takes $(g_0N,\dots,g_{k-1}N)$ to $(g_iN,g_jN)$. Since $H\cong(\C^*)^r$, a configuration thus gives rise to $r$ coordinates for each edge (see Figure~\ref{fig:EdgeCoordinates}) given by the edge minors $\Delta^{w_i}$. The following simple, but important, result illustrates the significance of the element $s_G$. 
\begin{lemma}\label{lem:ijVSji} Let $\alpha\in\Conf_k(\A)$. If $\pi_{ij}(\alpha)=h$ then $\pi_{ji}(\alpha)=w_0(h^{-1})s_G$. 
\end{lemma}
\begin{proof}
If $(g_iN,g_jN)=(N,\wzero hN)$, then $(g_jN,g_iN)=(\wzero hN,N)=(N,\wzero kN)$, where $k\in H$ equals $[\wzero^{-1}(\wzero h)^{-1}]_0=w_0(h^{-1})s_G$. This proves the result.
\end{proof}
By~\eqref{eq:PermuteCoordinates} this shows that when changing the orientation of an edge, the edge coordinates are permuted and multiplied by a sign (see Figure~\ref{fig:OrientationChange}).
\begin{lemma}\label{lemma:CanonicalRep} Let $\alpha\in\Conf_3(\A)$ and let $h_1=\pi_{01}(\alpha)$, $h_2=\pi_{12}(\alpha)$, and $h_3=\pi_{20}(\alpha)$. The canonical representative of $\alpha$ equals
\begin{equation}
(N,\wzero h_1N,uw_0(h_1)h_2s_GN),
\end{equation}
where $u$ is an element in $N_-$ satisfying that $[\wzero^{-1}u]_0=(w_0(h_3h_1)h_2)^{-1}$.
\end{lemma}
\begin{proof} The canonical representative has the form $(N,\wzero hN,ukN)$ for some $h,k\in H$, $u\in N_-$. By \eqref{eq:piij}, we have $h_1=[\wzero^{-1}\wzero h]_0=h$ and $h_2=[\wzero^{-1}(\wzero h)^{-1}uk]_0=w_0(h^{-1})ks_G$,
which together imply that $k=w_0(h_1)h_2s_G$, proving the first statement.
For the second statement, Lemma~\ref{lem:ijVSji} implies that
$w_0(h_3^{-1})s_G=\pi_{02}(\alpha)=[\wzero^{-1}uw_0(h_1)h_2s_G]_0$,
and it follows that $[\wzero^{-1}u]_0=(w_0(h_3h_1)h_2)^{-1}$ as desired.
\end{proof}

\begin{figure}[htb]
\begin{center}
\begin{minipage}[c]{0.47\textwidth}
\scalebox{0.8}{\input{figures_gen/EdgeCoordinates.tex}}
\end{minipage}
\hfill
\begin{minipage}[c]{0.4\textwidth}
\scalebox{0.9}{\input{figures_gen/OrientationChange.tex}}
\end{minipage}
\\
\begin{minipage}[t]{0.42\textwidth}
\caption{Edge coordinates.}\label{fig:EdgeCoordinates}
\end{minipage}
\begin{minipage}[t]{0.45\textwidth}
\caption{Changing  the orientation of an edge.}\label{fig:OrientationChange}
\end{minipage}
\end{center}
\end{figure}

\subsection{Face coordinates}\label{sec:MinorCoordinates}

Consider the maps
\begin{equation}\label{eq:FiberMaps}
H^3\to H,\quad (h_1,h_2,h_3)\to (w_0(h_3h_1)h_2)^{-1},\qquad N_-\cap\wzero G_0\to H,\quad u\mapsto [\wzero^{-1} u]_0.
\end{equation}
The following is a restatement of Lemma~\ref{lemma:CanonicalRep}. 
\begin{lemma}
We have an isomorphism of varieties
\begin{equation}\label{eq:H3u}
\Conf_3(\A)\to H^3\times_H N_-\cap\wzero G_0, \qquad \alpha\mapsto \big(\pi_{01}(\alpha),\pi_{12}(\alpha),\pi_{20}(\alpha),\pi_{N_-}(\alpha)\big),
\end{equation}
where $\pi_{N_-}$ is the map $\Conf_3(\A)\to N_-\cap\wzero G_0$ induced by the projection $B_-=N_-H\to N_-$, and $\times_H$ denotes the fiber product with respect to the maps~\eqref{eq:FiberMaps}.\qed
\end{lemma}


By Proposition~\ref{prop:DistinctWeights} there exist $j_1<\dots<j_{m-r}\in\{1,\dots,m\}$ such that the face minors are $\Delta^{\gamma_{j_1}},\dots,\Delta^{\gamma_{j_{m-r}}}$. We let $\Delta^{\accentset{\circ}{\gamma}_k}$ denote the face minor $\Delta^{\gamma_{j_k}}$.

\begin{proposition}\label{prop:EdgeFaceMinors} The edge and face minors define a birational equivalence
\begin{equation}
\begin{gathered}
\Delta\colon\Conf_3^*(\A)\cong H^3\times_H N_-\cap\wzero G_0\to (\C^*)^{3r}\times (\C^*)^{m-r},\\
(h_1,h_2,h_3,u)\mapsto \big(\{\Delta^{w_i}(h_1)\}_{i=1}^r,\{\Delta^{w_i}(h_2)\}_{i=1}^r,\{\Delta^{w_i}(h_3)\}_{i=1}^r,\{\Delta^{\accentset{\circ}{\gamma}_i}(u)\}_{i=1}^{m-r}\big).
\end{gathered}
\end{equation}
\end{proposition}
\begin{proof}
By definition, the edge minors $\Delta^{w_0w_i}$ of $u$ are the coordinates of $[\wzero^{-1}u]_0$, which by Lemma~\ref{lemma:CanonicalRep} are rational functions of $h_1$, $h_2$ and $h_3$. The result now follows from Corollary~\ref{cor:Nminus}.
\end{proof}

By Lemma~\ref{lem:ijVSji}, $\Delta$ also defines a birational equivalence (also denoted by $\Delta$)
 \begin{equation}
 \Delta\colon\Conf_3^*(\A)\times_{kl}^{s_G}\Conf_3^*(\A)\to (\C^*)^{5r}\times (\C^*)^{2(m-r)}.
 \end{equation}
In Section~\ref{sec:Rank2} we shall identify the codomains with seed tori when $G$ is $A_2$, $B_2$, $C_2$ or $G_2$.
 




\subsection{Comparison with the Ptolemy coordinates}\label{sec:PtolemyComparison}
Using the standard root datum (as in~\cite{Knapp}) for $G=\SL(n,\C)$, the group $H$ is the diagonal matrices, and $N$ is the upper triangular matrices with $1$ on the diagonal. The map $\chi_{\omega_i}\colon H\to \C^*$ takes ${\diag(a_1,\dots,a_n)}$ to $a_1a_2\cdots a_i$, and the element $\wzero$ is the counter diagonal matrix whose $(n+1-i,i)$ entry is $(-1)^{i-1}$.

Given a triple $(g_0N,g_1N,g_2N)$ of $N$-cosets in $\SL(n,\C)$, there is a Ptolemy coordinate $c_t$ for each triple $t=(t_0,t_1,t_2)$ of non-negative integers summing to $n$ defined by 
\begin{equation}
c_t=\det(\{g_0\}_{t_0},\{g_1\}_{t_1},\{g_2\}_{t_2}),
\end{equation}
where $\{g\}_{k}$ denotes the first $k$ column vectors of a matrix $g$ (see~\cite{GaroufalidisThurstonZickert,GaroufalidisGoernerZickert}). The Ptolemy coordinate $c_t$ of $(N,\wzero h_1N,uw_0(h_1)h_2s_GN)$ is up to a sign equal to the product (undefined terms are $1$) of $\chi_{\omega_{t_1}}(h_1)$, $\chi_{\omega_{t_2}}(w_0(h_1)h_2)$ and the $(t_2\times t_2)$-minor  of $u$ given by the rows $t_0+1,\dots, t_0+t_1$ and the columns $1,\dots,t_1$. Using the ``standard'' word $\mathbf i=t_{n-1}\cdots t_2t_1$, where $t_i=s_1s_2\cdots s_i$, it follows from Remark~\ref{rm:SLnMinors} that these minors are the $\mathbf i$-minors. To summarize, the Ptolemy coordinates are up to a sign and a monomial transformation equal to the minor coordinates for $\mathbf i$.

\subsection{The action on $\Conf_3^*(\A)$ by rotations}
\begin{proposition}\label{prop:ShiftFormula} The rotation map $\rot\colon\Conf_3^*(\A)\to\Conf_3^*(\A)$ taking $(g_0N,g_1N,g_2N)$ to $(g_2N,g_0N,g_1N)$ corresponds to the map
\begin{equation}
(h_1,h_2,h_3,u)\mapsto\big(h_3,h_1,h_2,h_2^{-1}w_0(h_1)^{-1}(\Phi\circ\Psi)^2(u)w_0(h_1)h_2\big).
\end{equation}
under the isomorphism~\eqref{eq:H3u}.
\end{proposition}
The proof uses the following technical lemmas.
\begin{lemma}\label{lem:Tech1} For any $u\in N_-\cap\wzero G_0$, we have $[\Psi(u)\wzero]_0=[\wzero^{-1}u]_0$.
\end{lemma}
\begin{proof}
Let $\wzero^{-1}u=yhx$. We then have $[\wzero^{-1}u]_0=h$, and $\Psi(u)=\Psi(x)h\Psi(y)\Psi(\wzero)$, from which it follows that $[\Psi(u)\wzero]_0=h=[\wzero^{-1}u]_0$.
\end{proof}

\begin{lemma}\label{lem:Tech3} For any $u\in N_-\cap\wzero G_0$, we have $[\wzero^{-1}\Phi\Psi(u)]_0=[\wzero^{-1}u]_0^{-1}$.
\end{lemma}
\begin{proof}
Let $\Psi(u)\wzero=\Phi\Psi(u)hx$ for $h\in H$, $x\in N$. Then $[\wzero^{-1}u]_0=[\Psi(u)\wzero]_0=h$ and 
\begin{equation}
\wzero^{-1}\Phi\Psi(u)=(\wzero^{-1}\Psi(u)\wzero)x^{-1}h^{-1}=(\wzero^{-1}\Psi(u)\wzero)h^{-1}(hx^{-1}h^{-1}).
\end{equation}
Hence, $[\wzero^{-1}\Phi\Psi(u)]_0=h^{-1}$, and the result follows.
\end{proof}

\begin{lemma}\label{lem:Tech2} For any $u\in N_-\cap\wzero G_0$, we have
\begin{equation}
u^{-1}N=\Psi\Phi\Psi(u)^{-1}[\wzero^{-1}u]_0^{-1}s_G\wzero N.
\end{equation}
\end{lemma}
\begin{proof} 
Suppose $\Psi(u)N_-=yh\wzero N_-$. We then have,
\begin{equation}
u^{-1}N=\big(\Psi\big(\Psi(u)N_-\big)\big)^{-1}=\big(\Psi\big(yh\wzero N_-\big)\big)^{-1}=\Psi(y)^{-1}h^{-1}\wzero N.
\end{equation}
By~\eqref{eq:PhiProp}, $y=\Phi\Psi(u)$ and $h=[\Psi(u)\wzero]_0s_G=[\wzero^{-1}u]_0s_G$. This proves the result.
\end{proof}

\begin{proof}[Proof of Proposition~\ref{prop:ShiftFormula}]
Let $\alpha=\big(N,\wzero h_1N,uw_0(h_1)h_2s_GN\big)\in\Conf_3(\A)$. One has
\begin{equation}
\begin{aligned}
\alpha&=\big(N,\wzero h_1N,uw_0(h_1)h_2s_GN\big)\\
&=\big(u^{-1}N\,,\,\wzero h_1N\,,\,w_0(h_1)h_2s_GN\big)\\
&=\big(\Psi\Phi\Psi(u)^{-1}[\wzero^{-1} u]_0^{-1}s_G\wzero N,\wzero h_1N,w_0(h_1)h_2s_GN\big)\\
&=\big([\wzero^{-1}u]_0^{-1}s_G\wzero N,\Psi\Phi\Psi(u)\wzero h_1N,w_0(h_1)h_2s_GN\big)\\
&=\big(h_2^{-1}w_0(h_1)^{-1}[\wzero^{-1}u]_0^{-1}\wzero N,h_2^{-1}w_0(h_1)^{-1}(\Phi\Psi)^2(u)[\Psi\Phi\Psi(u)\wzero]_0h_1s_GN,N\big)\\
&=\big(\wzero h_3N,h_2^{-1}w_0(h_1)^{-1}(\Phi\Psi)^2(u)w_0(h_1)h_2w_0(h_3)h_1s_GN,N\big).
\end{aligned}
\end{equation}
The third equality follows from Lemma~\ref{lem:Tech2}, the fifth from \eqref{eq:PhiProp}, and the last from Lemma~\ref{lemma:CanonicalRep}, which together with Lemmas~\ref{lem:Tech1} and \ref{lem:Tech3} imply that 
\begin{equation}
[\Psi\Phi\Psi(u)\wzero]_0=[\wzero^{-1}\Phi\Psi(u)]_0=[\wzero^{-1}u]_0^{-1}=w_0(h_3h_1)h_2. 
\end{equation}
This concludes the proof.
\end{proof}


\subsection{$\Conf_4^*(\A)$ and the flip}\label{sec:FlipSection}

For $\alpha=(g_0N,g_1N,g_2N,g_3N)\in\Conf_4^*(\A)$ let 
\begin{equation}
\begin{aligned}
\alpha_{012}&=(g_0s_GN,g_1N,g_2N),&\quad \alpha_{023}&=(g_0N,g_2N,g_3N),\\
\alpha_{123}&=(g_1N,g_2N,g_3N),&\quad\alpha_{013}&=(g_0N,g_1s_GN,g_2N),
\end{aligned}
\end{equation}
so that $\Psi_{02}(\alpha)=(\alpha_{012},\alpha_{023})$ and $\Psi_{13}(\alpha)=(\alpha_{123},\alpha_{013})$. We wish to relate the canonical rep\-resen\-tatives of  $\Psi_{02}(\alpha)$ to those of $\Psi_{13}(\alpha)$. Let $\alpha_{120}=\rot^{-1}(\alpha_{012})$ and $\alpha_{130}=\rot^{-1}(\alpha_{013})$. We then have
\begin{equation}\label{eq:RotPsi}
\Psi_{02}(\alpha)=(\rot(\alpha_{120}),\alpha_{023}),\qquad \Psi_{13}(\alpha)=(\alpha_{123},\rot(\alpha_{130})).
\end{equation}
Hence, by Proposition~\ref{prop:ShiftFormula} it is enough to relate the canonical representatives of $\alpha_{120}$ and $\alpha_{023}$ to those of $\alpha_{123}$ and $\alpha_{130}$. 

Each $\alpha\in\Conf_4^*(\A)$ has a unique representative of the form
$(N,yk_1N,\wzero k_2N,\Phi^{-1}(v)k_3N)$. Letting $h_{ij}=\pi_{ij}(\alpha)\in H$ it follows from Lemma~\ref{lemma:CanonicalRep} that this representative is given by 
\begin{equation}
\alpha=\big(N,w_0(h_{02})w_0(h_{12})^{-1}u^{-1}N,\wzero h_{02}N,\Phi^{-1}(v)\wzero h_{03}N\big).
\end{equation}
In particular, we have
\begin{equation}
\begin{aligned}
\alpha_{120}&=\big(w_0(h_{02})w_0(h_{12})^{-1}u^{-1}N,\wzero h_{02}N,s_GN)=(N,\wzero h_{12}N,uw_0(h_{12})w_0(h_{02}^{-1})s_GN)\\
\alpha_{023}&=(N,\wzero h_{02}N,\Phi^{-1}(v)\wzero h_{03}N)=(N,\wzero h_{02}N,vw_0(h_{02})h_{23}s_GN).
\end{aligned}
\end{equation}

Each element in $G_0$ also admits a factorization $xyh$ with $x\in N_+$, $y\in N_-$ and $h\in H$. In other words, the identity induces an isomorphism
\begin{equation}\label{eq:iotaDef}
\iota\colon N_-\times H\times N\to N\times N_-\times H.
\end{equation}

\begin{proposition}\label{prop:Flip} Let $k=w_0(h_{12})w_0(h_{02})^{-1}\in H$, and let $c,d\in N_-$, and $l\in H$ be elements satisfying that
$\iota(u,k,\Phi^{-1}(v))=(\Phi^{-1}(c),d,l)$. Then $l=h_{31}^{-1}h_{30}$, and we have 
\begin{equation}
\alpha_{123}=(N,\wzero h_{12}N,cw_0(h_{12})h_{23}s_GN),\enspace
\alpha_{130}=(N,\wzero w_0(h_{31}^{-1})N,dh_{31}^{-1}h_{30}s_GN).
\end{equation}
\end{proposition}
\begin{proof} By left multiplication by $\Phi^{-1}(c)^{-1}uk=dl\Phi^{-1}(v)^{-1}$, we have
\begin{equation}
\begin{aligned}
\alpha&=\big(N,w_0(h_{02})w_0(h_{12})^{-1}u^{-1}N,\wzero h_{02}N,\Phi^{-1}(v)\wzero h_{03}N\big)\\
&=\big(dlN,N,\Phi^{-1}(c)^{-1}\wzero h_{12}N,\wzero w_0(l)h_{03}N\big).
\end{aligned}
\end{equation}
This shows that $h_{13}=w_0(l)h_{03}$, yielding the formula for $l$. The formulas for $\alpha_{123}$ and $\alpha_{130}$ now follow from their definition.
\end{proof}

For the groups $A_2$, $B_2$, $C_2$ and $G_2$, Theorem~\ref{thm:Flip} states that after a monomial transformation, the minor coordinates of $\alpha_{012}$ and $\alpha_{023}$ are related to those of $\alpha_{123}$ and $\alpha_{013}$ by quiver mutations. The example below shows the much simpler case $A_1=\SL(2,\C)$. The case of $\SL(n,\C)$ is treated in~\cite[Sec.~10]{FockGoncharov}.

\begin{example}
For $G=\SL(2,\C)$, $s_G=-I$. There are no face coordinates, and the edge coordinates $\pi_{ij}$ are the Ptolemy coordinates $c_{ij}=\det(g_i\Vector{1}{0},g_j\Vector{1}{0})$. Figure~\ref{fig:MapConf4SL2} shows the corresponding coordinates in $\Conf_3^*(\A)\times_{kl}^{s_G}\Conf_3^*(\A)$. The Ptolemy coordinates satisfy the \emph{Ptolemy relation} $c_{03}c_{12}+c_{01}c_{23}=c_{02}c_{13}$, which is equivalent to $c_{02}(-c_{13})=c_{23}(-c_{01})+c_{12}(-c_{03})$, the mutation relation arising from a mutation at the middle vertex of the quiver shown on the right in Figure~\ref{fig:MapConf4SL2}.
\begin{figure}[htb]
\begin{center}
\scalebox{0.8}{\input{figures_gen/MapConf4SL2.tex}}
\caption{Ptolemy coordinates of a tuple and edge coordinates of its images in $\Conf_3^*(\A)\times_{02}^{s_G}\Conf_3^*(\A)$ and $\Conf_3^*(\A)\times_{13}^{s_G}\Conf_3^*(\A)$.}\label{fig:MapConf4SL2}
\end{center}
\end{figure}
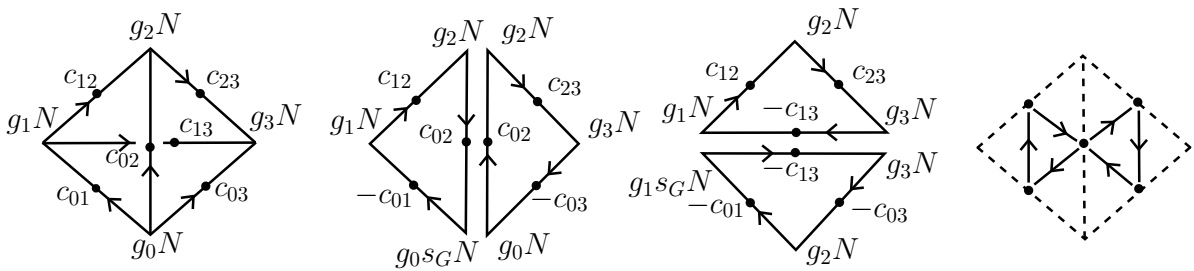
\end{example}

\section{The natural cocycle}\label{sec:NaturalCocycle}
We now show that there is an explicit one-to-one correspondence between $\Conf_k^*(\A)$ and certain $G$-valued 1-cocycles on a truncated simplex labeling long edges by elements in $\wzero H$ and short edges by elements in $N\cap G_0\wzero$. This result allows us to explicitly recover a representation from its coordinates. 

For a CW complex $X$ let $V(X)$ denote the set of vertices of $X$ and $E(X)$ the set of oriented edges. A \emph{$G$-cocycle} on $X$ is a function $\tau\colon E(X)\to G$ such that $\tau(\varepsilon_1)\tau(\varepsilon_2)\dots\tau(\varepsilon_l)=1$, whenever $\varepsilon_1\cdot\varepsilon_2\cdots\varepsilon_l$ is a contractible loop. The \emph{coboundary} of a \emph{$0$-cochain} $\eta\colon V(X)\to G$ is the $G$-cocycle taking an edge from vertex $v$ to vertex $w$ to $\eta(v)^{-1}\eta(w)$.

Let $\Delta^n$ denote a standard $n$-simplex, and let $\overline\Delta^n$ denote the corresponding truncated simplex. Let $v_{ij}$ denote the vertex of $\overline\Delta^n$ near vertex $i$ of $\Delta^n$ on the edge between $i$ and $j$ of $\Delta^n$. Each edge of $\overline\Delta^n$ is either \emph{long} (from $v_{ij}$ to $v_{ji}$) or \emph{short} (from $v_{ij}$ to $v_{ik}$).

\begin{definition}
A $G$-cocycle $\tau$ on $\overline\Delta^n$ is a \emph{natural cocycle} if $\tau(\varepsilon)\in N\cap G_0\wzero$ when $\varepsilon$ is a short edge, and $\tau(\varepsilon)\in \wzero H$, when $\varepsilon$ is a long edge.
\end{definition}

\begin{convention} Given a natural cocycle on $\overline\Delta^n$, we denote the labeling of the short edge from $v_{ij}$ to $v_{ik}$ by $\beta^i_{jk}$,  and the labeling of the long edge from $v_{ij}$ to $v_{ji}$ by $\alpha_{ij}$ (see Figures~\ref{fig:NaturalCocyclePart} and \ref{fig:NaturalCocycle}).
\end{convention}

\begin{definition}\label{def:CocycleForAlpha} Let $\alpha=(g_0N,\dots,g_{k-1}N)\in\Conf_k^*(\A)$. The \emph{natural cocycle} associated to $\alpha$ is the coboundary of the $0$-cochain $\eta_\alpha$ taking $v_{ij}$ to $g$ if $(g_iN,g_jN)=g(N,\wzero hN)$ with $h\in H$.
\end{definition}

Note that the set of natural cocycles is a variety, and that the map taking a configuration to its natural cocycle is an isomorphism. We wish to give an explicit formula for the edges. 
\begin{figure}[htb]
\begin{center}
\begin{minipage}[c]{0.47\textwidth}
\scalebox{0.75}{\input{figures_gen/NaturalCocycle3d.tex}}
\end{minipage}
\hfill
\begin{minipage}[c]{0.47\textwidth}
\scalebox{0.75}{\input{figures_gen/NaturalCocycle.tex}}
\end{minipage}
\\
\begin{minipage}[t]{0.48\textwidth}
\caption{Natural cocycle on a $3$-simplex.}\label{fig:NaturalCocyclePart}
\end{minipage}
\begin{minipage}[t]{0.47\textwidth}
\caption{Natural cocycle on a 2-simplex.}\label{fig:NaturalCocycle}
\end{minipage}
\end{center}
\end{figure}

\begin{lemma}\label{lem:CocycleLemma} Let $h,k\in H$. The natural cocycle for $(N,\wzero hN,uks_GN)$ has 
\begin{equation}
\beta^2_{01}=k^{-1}\Psi\Phi\Psi(u)k.
\end{equation}
\end{lemma}
\begin{proof}
We have $\beta^2_{01}=\eta(v_{20})^{-1}\eta(v_{21})$, where $\eta$ is the $0$-cochain from Definition~\ref{def:CocycleForAlpha}. Let $y=k^{-1}uk$. By Lemma~\ref{lem:Tech2}, we have $y^{-1}k^{-1}s_GN=\Psi\Phi\Psi(y)^{-1}[\wzero y]_0^{-1}\wzero k^{-1}N$. Hence,
\begin{equation}
(uks_GN,N)=uks_G(N,y^{-1}k^{-1}s_GN)=uks_G\Psi\Phi\Psi(y)^{-1}(N,[\wzero y]_0^{-1}\wzero k^{-1}N),
\end{equation}
from which it follows that $\eta(v_{20})=uks_G\Psi\Phi\Psi(y)^{-1}$. Similarly,
\begin{equation}
(uks_GN,\wzero hN)=uks_G(N,\wzero hw_0(k)s_GN),
\end{equation}
so $\eta(v_{21})=uks_G$. It follows that $\beta^2_{01}=\eta(v_{20})^{-1}\eta(v_{21})=\Psi\Phi\Psi(y)=k^{-1}\Psi\Phi\Psi(u)k$.
\end{proof}

\begin{proposition}\label{prop:NaturalCocycleFormula} Let $\alpha=(h_1,h_2,h_3,u)\in H^3\times_H N_-\cap\wzero G_0=\Conf_3^*(\A)$, and let 
\begin{equation}
u_0=u,\enspace u_1=h_2^{-1}w_0(h_1)^{-1}(\Phi\Psi)^2(u_0)w_0(h_1)h_2,\enspace u_2=h_1^{-1}w_0(h_3)^{-1}(\Phi\Psi)^2(u_1)w_0(h_3)h_1.
\end{equation}
The natural cocycle for $\alpha$ is given by
\begin{equation}
\begin{aligned}
\beta^2_{01}&=h_2^{-1}w_0(h_1)^{-1}\Psi\Phi\Psi(u_0)w_0(h_1)h_2,&\qquad \beta^1_{20}&=h_1^{-1}w_0(h_3)^{-1}\Psi\Phi\Psi(u_1)w_0(h_3)h_1\\
\beta^0_{12}&=h_3^{-1}w_0(h_2)^{-1}\Psi\Phi\Psi(u_2)w_0(h_2)h_3,&\qquad \alpha_{01}&=\wzero h_1,\quad \alpha_{12}=\wzero h_2,\quad \alpha_{20}=\wzero h_3.
\end{aligned}
\end{equation}
\end{proposition}
\begin{proof}
The formula for the long edges $\alpha_{ij}$ is an immediate consequence of the definition, and the formula for the short edges $\beta^i_{jk}$ follow from Lemma~\ref{lem:CocycleLemma} and Proposition~\ref{prop:ShiftFormula}.
\end{proof}

\section{Deriving explicit formulas}\label{sec:ExplicitFormulas}
We now derive formulas for the rotation (Proposition~\ref{prop:ShiftFormula}) and the flip (Proposition~\ref{prop:Flip}) in terms of the minor coordinates. Explicit computations are given for the rank two groups in Section~\ref{sec:Rank2}. Let $N_-^{\neq}$ denote the open subset of $N_-\cap\wzero G_0$ with non-vanishing $\mathbf i$-minors, and let
\begin{equation}
N^{\mathbf i}=x_{\mathbf i}((\C^*)^m)\cap G_0\wzero, \qquad N_-^{\bar{\mathbf i}}=y_{\bar{\mathbf i}}((\C^*)^m)\cap\wzero G_0,
\end{equation}
where $y_{\bar{\mathbf i}}(s_1,\dots,s_m)=y_{i_m}(s_1)y_{i_{m-1}}\dots y_{i_1}(s_m)$. Note that the factorization of elements in $N_-^{\bar{\mathbf i}}$ is with respect to the opposite word $\bar{\mathbf i}=s_{i_m}\cdots s_{i_2}s_{i_1}$. The factorization coordinates on $N^{\mathbf i}$, $N_-^{\bar{\mathbf i}}$ and the $\mathbf i$-minors on $N_-^{\neq}$ define canonical birational equivalences of each of these spaces with $(\C^*)^m$. 
\subsection{Rotations}\label{sec:ExplicitRotations}
By Proposition~\ref{prop:ShiftFormula} we need a formula for $(\Phi\Psi)^2$ and a formula for how the minor coordinates change under conjugation. 
We begin with the latter.

\begin{lemma}\label{lem:ConjMinors} For any $w\in W$,
$\Delta^{w\omega_i}(k^{-1}uk)=\chi_{\omega_i}(w^{-1}(k^{-1})k)\Delta^{w\omega_i}(u)$.
\end{lemma}
\begin{proof}
For $u\in \overline w G_0$, one easily checks that for $[\overline w^{-1} k^{-1}uk]_0=w^{-1}(k^{-1})k[\overline w^{-1}u]_0$. This proves the result.
\end{proof}

To obtain a formula for $(\Phi\Psi)^2$ first observe that
\begin{equation}
(\Phi\Psi)^n=(\Psi\Phi\Psi)^{-1}\circ(\Psi\Phi)^n\circ\Psi\Phi\Psi,\qquad n\in\Z.
\end{equation}
The basic observation below allows us to apply  Theorem~\ref{thm:MinorsAndFactorization} to explictly compute $\Psi\Phi\Psi$.

\begin{lemma}\label{lem:PsiPhiPsi} For any $u\in N_-\cap \wzero G_0$, we have $\pi_-(\Psi\Phi\Psi(u))=u$.
\end{lemma}
\begin{proof}
Let $\wzero^{-1}u=yhx$. Then $x=\pi_-^{-1}(u)$, and $u=\wzero yhx$. Hence, $\Psi(u)=\Psi(x)h\Psi(y)\wzero s_G$, so $\Phi\Psi(u)=\Psi(x)$. This proves the result.
\end{proof}
\begin{corollary}\label{cor:PsiPhiPsi} The map $\Psi\Phi\Psi$ extends to a biregular isomorphism $N_-^{\neq}\to N^{\mathbf i}$ given explicitly by~\eqref{eq:MinorsFact}.\qed.
\end{corollary}

%
%

\begin{remark}
By Proposition~\ref{prop:NaturalCocycleFormula}, this provides an explicit formula for the natural cocycle of $\alpha\in\Conf_3(\A)$ whenever the minor coordinates of $\alpha$, $\rot(\alpha)$ and $\rot^2(\alpha)$ are non-zero. For $G=\SL(n,\C)$ and  the ``standard word'' (see Section~\ref{sec:PtolemyComparison}) this formula agrees with the one given in~\cite{GaroufalidisThurstonZickert} via \emph{diamond coordinates}.
\end{remark}

\subsection{The flip}
For all $u\in N_-\cap\wzero G_0$, we have 
\begin{equation}
u=\Psi(\Psi\Phi)^{-1}\Psi\Phi\Psi(u),\qquad \Phi^{-1}(u)=(\Psi\Phi)^{-2}\Psi\Phi\Psi(u).
\end{equation}
This motivates the definition of birational equivalences
\begin{equation}
\begin{gathered}
\Gamma_1\colon N_-^{\neq}\to N_-^{\bar{\mathbf i}}, \qquad\Gamma_2\colon N_-^{\neq}\to N^{\mathbf i}\\
u\mapsto \Psi(\Psi\Phi)^{-1}\Psi\Phi\Psi(u),\qquad u\mapsto(\Psi\Phi)^{-2}\Psi\Phi\Psi(u).
\end{gathered}
\end{equation}

Let $f_{(c,d,l)}$ denote the composition
\begin{equation}\label{eq:fabl}
\xymatrixcolsep{1.1pc}\cxymatrix{{N_-^{\neq}\times H\times N_-^{\neq}\ar[rrr]^-{(\Gamma_1,\id,\Gamma_2)}&&&N_-^{\bar{\mathbf i}}\times H\times N^{\mathbf i}\ar[r]^-{\iota}&N^{\mathbf i}\times N_-^{\bar{\mathbf i}}\times H\ar[rrr]^-{(\Gamma_2^{-1},\Gamma_1^{-1},\id)}&&&N_-^{\neq}\times N_-^{\neq}\times H,}}
\end{equation}
where $\iota$ is the map~\eqref{eq:iotaDef}.
Note that if $u$, $v$, $k$, $c$, $d$, and $l$ are as in Proposition~\ref{prop:Flip}, then 
$(c,d,l)=f_{(c,d,l)}(u,k,v)$. In particular, the flip is given explicitly in terms of $\iota$ and the maps $\Psi\Phi\Psi$, $\Psi\Phi$ and their inverses.


\subsection{Formulas for $\Psi\Phi$ and $\iota$}\label{sec:PsiPhiIota}
The maps $\Psi\Phi$ and $\iota$ can be computed explicitly using the following elementary properties (see e.g.~\cite{Lusztig}):
\begin{gather}
x_i(s)y_j(t)=y_j(t)x_i(s),\quad i\neq j\label{eq:ijcommute}\\
x_i(s)y_i(t)=y_i(\frac{t}{1+st})h_i^{1+st}x_i(\frac{s}{1+st}),\quad y_i(s)x_i(t)=x_i(\frac{t}{1+st})h_i^{\frac{1}{1+st}}y_i(\frac{s}{1+st}),\label{eq:ii}\\
h_i^sy_j(t)=y_j(ts^{-A_{ij}})h_i^s,\quad h_i^{s}x_j(t)=x_j(ts^{A_{ij}})h_i^s\label{eq:ijconj}\\
x_j(t)\overline{s_j}\,\overline{w}B=y_j(1/t)\overline{w}B, \qquad w=s_{j_1}\cdots s_{j_k},\quad s_jw\text{ reduced}.\label{eq:xwyw}
\end{gather}
\begin{example}
We compute $\Psi\Phi$ for the group $A_2$ using the word $\mathbf i=(1,2,1)$. The Cartan matrix is $\Matrix{2}{-1}{-1}{2}$ and we have
\begin{equation}
\begin{aligned}
x_1(a)x_2(b)x_1(c)\overline{s_1s_2s_1}B&=x_1(a)x_2(b)y_1(1/c)\overline{s_2s_1}B\\
&=x_1(a)y_1(1/c)y_2(1/b)\overline{s_1}B\\
&=y_1(\frac{1}{a+c})h_1^{1+a/c}x_1(\frac{ac}{a+c})y_2(1/b)\overline{s_1}B\\
&=y_1(\frac{1}{a+c})y_2(\frac{1}{b}(1+a/c))y_1(\frac{c}{a(a+c)})B
\end{aligned}
\end{equation}
proving that $\Psi\Phi(a,b,c)=(\frac{c}{a(a+c)},\frac{a+c}{bc},\frac{1}{a+c})$.
\end{example}
\begin{example} This toy example illustrates how to compute $\iota$. Assume that $A_{12}=-1$.
\begin{equation}
\begin{aligned}
y_2(a)y_1(b)x_1(c)x_2(d)&=y_2(a)x_1(\frac{c}{1+bc})h_1^{\frac{1}{1+bc}}y_1(\frac{b}{1+bc})x_2(d)\\&=x_1(\frac{c}{1+bc})y_2(a)x_2(d(1+bc))y_1(b(1+bc))h_1^{\frac{1}{1+bc}}\\
&=x_1(\frac{c}{1+bc})x_2(\frac{d(1+bc)}{1+ad(1+bc)})y_2(a(1+ad(1+bc)))\\&\phantom{{}={}}y_1(\frac{b (1 + b c)}{1 + ad(1 + bc)})h_1^{\frac{1}{1+bc}}h_2^{\frac{1}{1+ad(1+bc)}}.
\end{aligned}
\end{equation}
\end{example}

\section{Groups of rank 2}\label{sec:Rank2}
We now compute the functions in Section~\ref{sec:ExplicitFormulas} explicitly for the groups $A_2$, $B_2$, $C_2$ and $G_2$. There are two reduced words: $(1,2,1)$ and $(2,1,2)$ for $A_2$, $(1,2,1,2)$ and $(2,1,2,1)$ for $B_2$, and $(1,2,1,2,1,2)$ and $(2,1,2,1,2,1)$ for $G_2$. We shall always use the word starting with $1$.

We use the root data from Knapp~\cite[Appendix C]{Knapp}. We identify $\h^*$ with $\R^2$ for $B_2$ and $C_2$, and with $\{v\in\R^3|\langle v,e_1+e_2+e_3\rangle=0\}$ for $A_2$ and $G_2$. The $e_i$ are the standard basis vectors, and $\langle,\rangle$ is the standard inner product.
\begin{equation}\label{eq:ExplicitRootData}
\begin{aligned}
&A_2:&\alpha_1=e_1-e_2,\quad&\alpha_2=e_2-e_3, &\omega_1&=e_1,&\quad\omega_2&=e_1+e_2,\\ 
&B_2:&\alpha_1=e_1-e_2,\quad&\alpha_2=e_2, &\omega_1&=e_1,&\quad\omega_2&=\frac{1}{2}(e_1+e_2),\\
&C_2:&\alpha_1=e_1-e_2,\quad&\alpha_2=2e_2, &\omega_1&=e_1,&\quad\omega_2&=e_1+e_2,\\ 
&G_2:& \alpha_1=e_1-e_2,\quad&\alpha_2=-2e_1+e_2+e_3,\quad &\omega_1&=-e_2+e_3,&\quad\omega_2&=-e_1-e_2+2e_3, 
\end{aligned}
\end{equation}
Using this, one easily verifies that $w_0(w_i)=-w_i$ for $B_2$, $C_2$ and $G_2$, and that $w_0(w_i)=-w_{3-i}$ for $A_2$, proving that $\sigma_G$ is trivial for $B_2$, $C_2$ and $G_2$, and non-trivial for $A_2$. The Cartan matrices are $A_2=\Matrix{2}{-1}{-1}{2}$, $B_2=\Matrix{2}{-1}{-2}{2}$, $C_2=\Matrix{2}{-2}{-1}{2}$, and $G_2=\Matrix{2}{-3}{-1}{2}$,
and one has
\begin{equation}
s_{A_2}=(1,1),\qquad s_{B_2}=(1,-1),\qquad s_{C_2}=(-1,1),\qquad s_{G_2}=(1,1)
\end{equation}
under the identification~\eqref{eq:CoordsOnH} of $H$ with $(\C^*)^2$.
\subsubsection{The map $\Psi\Phi\Psi$} 
By Corollary~\ref{cor:PsiPhiPsi}, the map $\Psi\Phi\Psi\colon N_-^{\neq}\to N^{\mathbf i}$ is given explicitly by~\eqref{eq:MinorsFact}. Displayed below are $\Psi\Phi\Psi(u_1,\dots,u_m)$ and $(\Psi\Phi\Psi)^{-1}(t_1,\dots,t_m)$, where the $u_i$ are the minor coordinates $\Delta^{\gamma_i}$ on $N_-^{\neq}$ (see Corollary~\ref{cor:Nminus}), and the $t_i$ are the factorization coordinates on $N^{\mathbf i}$.



\begin{equation}\label{eq:PsiPhiPsiFormula}
\begin{aligned}
&A_2:&&\!\!\!(\frac{u_2}{u_1u_3},\frac{u_3}{u_2},\frac{1}{u_3}),&(\frac{1}{t_1t_2},\frac{1}{t_2t_3},\frac{1}{t_3}) \\
&B_2:& &\!\!\!(\frac{u_2^2}{u_1u_3},\frac{u_3}{u_2u_4},\frac{u_4^2}{u_3},\frac{1}{u_4}),
&(\frac{1}{t_1t_2^2t_3},\frac{1}{t_2t_3t_4},\frac{1}{t_3t_4^2},\frac{1}{t_4})\\
&C_2:&&\!\!\!(\frac{u_2}{u_1u_3},\frac{u_3^2}{u_2u_4},\frac{u_4}{u_3},\frac{1}{u_4}),
&(\frac{1}{t_1t_2t_3},\frac{1}{t_2t_3^2t_4},\frac{1}{t_3t_4},\frac{1}{t_4})\\
&G_2:& &\!\!\!(\frac{u_2}{u_1u_3},\frac{u_3^2}{u_2u_4},\frac{u_4}{u_3u_5},\frac{u_5^3}{u_4u_6},\frac{u_6}{u_5},\frac{1}{u_6}),&
(\frac{1}{t_1t_2t_3^2t_4t_5},\frac{1}{t_2t_3^3t_4^2t_5^3t_6},\frac{1}{t_3t_4t_5^2t_6},\frac{1}{t_4t_5^3t_6^2},\frac{1}{t_5t_6},\frac{1}{t_6})
\end{aligned}
\end{equation}

\subsubsection{Formula for $\Psi\Phi$} Using the algorithm in Section~\ref{sec:ExplicitRotations} we obtain (the displayed formulas are $\Psi\Phi(t_1,\dots,t_m)$ and $(\Psi\Phi)^{-1}(s_1,\dots,s_m)$)
\begin{equation}\label{eq:PsiPhiFormula}
\begin{aligned}
&A_2:& (\frac{t_3}{t_1^2+t_1t_3},\frac{t_1+t_3}{t_2t_3},\frac{1}{t_1+t_3}), \qquad (\frac{1}{s_1+s_3},\frac{s_1+s_3}{s_1s_2},\frac{s_1}{s_1s_3+s_3^2})\\
&B_2:& (\frac{t_3t_4^2}{t_1\alpha_1},\frac{\alpha_1}{t_2t_3t_4\alpha_2},\frac{\alpha_2^2}{\alpha_1},\frac{1}{\alpha_2}),\qquad(\frac{1}{\beta_1},\frac{\beta_1}{\beta_2},\frac{\beta_2^2}{s_1s_2^2s_3\beta_1},\frac{s_1s_2}{s_4\beta_2}),\\
&&\alpha_1=t_3t_4^2 + t_1(t_2 + t_4)^2,\quad \alpha_2=t_2+t_4,
\quad \beta_1=s_1+s_3,\quad \beta_2=s_1 s_2 + (s_1 + s_3)s_4\\
&C_2:& (\frac{t_3t_4}{t_1\alpha_1},\frac{\alpha_1^2}{t_2t_3^2t_4\alpha_2},\frac{\alpha_2}{\alpha_1},\frac{1}{\alpha_2}),\qquad(\frac{1}{\beta_1},\frac{\beta_1^2}{\beta_2},\frac{\beta_2}{s_1s_2s_3\beta_1},\frac{s_1^2s_2}{s_4\beta_2}),\\
&&\alpha_1=t_3t_4 + t_1(t_2 + t_4),\quad \alpha_2=t_2+t_4,
\quad \beta_1=s_1+s_3,\quad \beta_2=s_1^2s_2+(s_1+s_3)^2s_4,
\end{aligned}
\end{equation}
for $A_2$, $B_2$ and $C_2$, while for $G_2$, we have
\begin{equation}
\begin{aligned}
\Psi\Phi(t)= (\frac{t_3t_4t_5^2t_6}{t_1\alpha_1},\frac{\alpha_1^3}{t_2t_3^3t_4^2t_5^3t_6\alpha_2},\frac{\alpha_2}{\alpha_1\alpha_3},\frac{\alpha_3^3}{\alpha_2\alpha_4},\frac{\alpha_4}{\alpha_3},\frac{1}{\alpha_4}),\\
(\Psi\Phi)^{-1}(s)=
(\frac{1}{\beta_1},\frac{\beta_1^3}{\beta_2},\frac{\beta_2}{\beta_1\beta_3},\frac{\beta_3^3}{\beta_2\beta_4},\frac{\beta_4}{s_1s_2s_3^2s_4s_5\beta_3},\frac{s_1^3s_2^2s_3^3s_4}{s_6\beta_4}),\\
\alpha_1=t_4 (t_1 t_2 t_3^2 + t_1 t_5^2 t_6 + t_3 t_5^2 t_6)+t_1 t_2 t_6 (t_3 + t_5)^2,\\
\alpha_2=t_4 (t_2 t_3^3 t_4 + 2 t_2 t_3^3 t_6 + 3 t_2 t_3^2 t_5 t_6 + t_5^3 t_6^2)+t_2t_6^2(t_3 + t_5)^3,\\
\alpha_3=t_1 t_2 + t_1 t_4 + t_3 t_4 + t_1 t_6 + t_3 t_6 + t_5 t_6,\quad\alpha_4=t_2 + t_4 + t_6,\\
\beta_1=s_1 + s_3 + s_5,\quad \beta_2=s_6(s_1 + s_3 + s_5)^3 +s_4(s_1 + s_3)^3+s_1^3s_2,\\
\beta_3=s_1^2s_2(s_3+s_5)+(s_1+s_3)^2s_3s_4,\\
\beta_4=s_1^2 s_2 s_4 (s_1 s_2 s_3^3 +3 s_1 s_3 s_5^2 s_6 + 3 s_3^2 s_5^2 s_6 + 
   2 s_1 s_5^3 s_6 + 3 s_3 s_5^3 s_6)+\\s_1^3s_2^2s_6(s_3 + s_5)^3+s_4^2s_5^3s_6(s_1 + s_3)^3.
\end{aligned}
\end{equation}
%

\subsubsection{Formula for $k^{-1}uk$} For $u\in N_-^{\neq}$ let $u_i=\Delta^{\gamma_i}(u)$ be the $i$th coordinate, and let $k_1$ and $k_2$ denote the coordinates of $k\in H$. Note that $\Delta^{\gamma_1}$ and $\Delta^{\gamma_2}$ are always edge minors, so we shall only need formulas for $(k^{-1}uk)_i$ for $2<i\leq m$. These can be computed using Lemma~\ref{lem:ConjMinors} using the fact that $\overline w\chi_{\alpha}^*(t)\overline w^{-1}=\chi_{w(\alpha)}^*(t)$. We obtain
\begin{equation}\label{eq:ConjRank2}
\begin{aligned}
&A_2:& (k^{-1}uk)_3&=k_1^2/k_2u_3\\
&B_2:& (k^{-1}uk)_3&=k_2^2u_3,\qquad&(k^{-1}uk)_4&=k_2^2/k_1u_4\\
&C_2:& (k^{-1}uk)_3&=k_2u_3, \qquad&(k^{-1}uk)_4&=k_2^2/k_1^2u_4\\
&G_2:& (k^{-1}uk)_3&=k_2u_3,\qquad&(k^{-1}uk)_4&=k_2^3/k_1^3u_4,\\
&&(k^{-1}uk)_5&=k_2/k_1u_5,\qquad&(k^{-1}uk)_6&=k_2^2/k_1^3u_6.
\end{aligned}
\end{equation}

\subsection{A monomial transformation}\label{sec:MonomialMap}
Define a monomial transformation $m_G\colon T_{Q_G}\to T_{Q_G}$ as follows:
\begin{equation}\label{eq:MonomialFormula}
\begin{aligned}
&&m_G^*(a_{ij})&=a_{ij},& m_{A_2}^*(a_1)&=a_1\frac{a_{01}a_{12}}{a_{10}}&&\\
m_{B_2}^*(a_1)&=a_1a_{12},& m_{B_2}^*(a_2)&=a_2\frac{a_{21}}{a_{10}},& m_{C_2}^*(a_1)&=a_1a_{12},& m_{C_2}^*(a_2)&=a_2\frac{a_{21}a_{01}^2}{a_{10}}\\
m_{G_2}^*(a_1)&=a_1a_{12},& m_{G_2}^*(a_2)&=a_2\frac{a_{21}a_{01}^3}{a_{10}},& m_{G_2}^*(a_3)&=a_3a_{01}a_{12},& m_{G_2}^*(a_4)&=a_4\frac{a_{21}a_{01}^3}{a_{10}}&
\end{aligned}
\end{equation}

We identify the codomain of the map $\Delta$ in Proposition~\ref{prop:EdgeFaceMinors} with the seed torus of $Q_G$ by identifying the edge coordinates with the frozen coordinates and the face coordinates with the non-frozen coordinates.
We can now define the map $\M\colon\Conf_3^*(\A)\to T_{Q_G}$ in Theorem~\ref{thm:Rotation} to be the composition of $\Delta$ with $m_G$. This is illustrated in Figure~\ref{fig:QuiverG2Minors} for $G=G_2$. Similarly, one identifies the codomain of $\Delta\colon\Conf_3^*(\A)\times^{s_G}_{kl}\Conf_3^*(\A)\to (\C^*)^{5}\times(\C^*)^{2(m-2)}$ with the seed torus $T_{Q_G\cup_{kl} Q_G}$ for $kl=02$ or $13$.
\begin{figure}[htb]
\begin{center}
\hspace{-1cm}\scalebox{0.72}{\input{figures_gen/QuiverG2Minors.tex}}
\caption{The map $\M\colon\Conf_3^*(\A)\to T_{Q_{G_2}}$.}\label{fig:QuiverG2Minors}
\end{center}
\end{figure}

\subsection{Proof of Theorem~\ref{thm:Rotation}}\label{sec:RotProof}
We wish to prove that 
\begin{equation}
\mu_G^{\rot}\M(h_{01},h_{12},h_{20},u)=\M\big(h_{20},h_{01},h_{12},(w_0(h_{01})h_{12})^{-1}(\Phi\Psi)^2(u)w_0(h_{01})h_{12}\big)
\end{equation}
for all $\alpha=(h_{01},h_{12},h_{20},u)\in N_-\cap \wzero G_0\cong\Conf_3^*(\A)$.
This is simply a matter of applying the explicit formulas above and comparing with the formula for $\mu_G^{\rot}$ in Section~\ref{sec:ExplicitMutations}. Clearly the frozen coordinates correspond, so we only need to check the non-frozen coordinates. We do this for $C_2$ and leave the other groups to the reader. Let $a_{ij}$ and $a_{ji}$ denote the coordinates of $h_{ij}$. These are the frozen coordinates.
By~\eqref{eq:MurotABC} the non-frozen coordinates of $\mu_{C_2}^{\rot}\M(\alpha)$ are given by
\begin{equation}\label{eq:RotC2Recap}
a_1'=\frac{a_{01} a_{02} a_{12} + a_2 a_{20}}{a_1},\quad a_2'=
\frac{a_{10} (a_{01} a_{02} a_{12} + a_2 a_{20})^2 + a_{01}^2 a_{02} a_1^2 a_{21}}{a_1^2 a_2}.\end{equation}
By Lemma~\ref{lemma:CanonicalRep}
we have
\begin{equation}\label{eq:u1u2}
(u_1,u_2)=[\wzero^{-1}u]_0=(w_0(h_{20}h_{01})h_{12})^{-1}=(\frac{a_{20}a_{01}}{a_{12}},\frac{a_{02}a_{10}}{a_{21}})
\end{equation}
and by~\eqref{eq:MonomialFormula}, $a_1=u_3a_{12}$ and $a_2=u_4\frac{a_{21}a_{01}^2}{a_{10}}$. Plugging these into \eqref{eq:RotC2Recap} and using~\eqref{eq:u1u2}
we obtain
\begin{equation}
a_1'=\frac{a_{01}a_{21}(u_2+u_1u_4)}{a_{10}u_3},\qquad a_2'=\frac{a_{21}(u_2^2+u_1^2u_4^2+u_2(u_3^2+2u_1u_4))}{u_3^2u_4}.
\end{equation}

We now compare this to the coordinates of $\M(\rot(\alpha))$.
Using~\eqref{eq:PsiPhiPsiFormula} and \eqref{eq:PsiPhiFormula} we obtain
\begin{equation}
(\Phi\Psi)^2(u_1,u_2,u_3,u_4)=\Big(u_1, u_2,\frac{u_2 + 
u_1 u_4}{u_3},\frac{u_2 (u_2^2 + u_1^2 u_4^2 + 
u_2 (u_3^2 + 2 u_1 u_4))}{u_1^2 u_3^2 u_4}\Big).
\end{equation}
Hence, by~\eqref{eq:ConjRank2} and \eqref{eq:MonomialFormula}, the non-frozen coordinates of $\M(\rot(\alpha))$ are 
\begin{equation}
k_2\frac{u_2 + 
u_1 u_4}{u_3}a_{01},\qquad \frac{k_2^2}{k_1^2}\frac{u_2 (u_2^2 + u_1^2 u_4^2 + 
u_2 (u_3^2 + 2 u_1 u_4))}{u_1^2 u_3^2 u_4}\frac{a_{10}a_{20}^2}{a_{02}},
\end{equation}
where  $k_1=a_{01}^{-1}a_{12}$ and $k_2=a_{10}^{-1}a_{21}$ are the coordinates of $k=w_0(h_{01})h_{12}$. Using~\eqref{eq:u1u2} it follows that these equal $a_1'$ and $a_2'$, respectively. This proves the result.

\subsection{Proof of Theorem~\ref{thm:Flip}}
Let $\alpha_{ijk}$ be as in Section~\ref{sec:FlipSection}. By Theorem~\ref{thm:Rotation} and \eqref{eq:RotPsi} we must prove that
\begin{equation}\label{eq:ProofFlip}
\mu_G^{\flip}\big(\mu_G^{\rot}\M(\alpha_{120}),\M(\alpha_{023})\big)=\big(\M(\alpha_{123}),\mu_G^{\rot}\M(\alpha_{130})\big).
\end{equation}
As in Section~\ref{sec:FlipSection} we may assume that 
\begin{equation}\label{eq:alphaijkformula}
\begin{aligned}
\alpha_{120}&=(N,\wzero h_{12}N,uw_0(h_{12})w_0(h_{02}^{-1})s_GN),&\alpha_{023}&=(N,\wzero h_{02}N,vw_0(h_{02})h_{23}s_GN)\\
\alpha_{123}&=(N,\wzero h_{12}N,cw_0(h_{12})h_{23}s_GN),&
\alpha_{130}&=(N,\wzero w_0(h_{31}^{-1})N,dh_{31}^{-1}h_{30}s_GN),
\end{aligned}
\end{equation}
where $(c,d,l)=f_{(c,d,l)}(u,k,v)$. As in the proof of Theorem~\ref{thm:Rotation} this is simply a matter of computing both sides of \eqref{eq:ProofFlip} using the explicit formulas for $\Psi\Phi$, $\Psi\Phi\Psi$ and their inverses, and the algorithm for computing $\iota$. We give a detailed proof only for $G=C_2$. Clearly, the frozen coordinates agree, so we only need to consider the non-frozen coordinates. 
Let $a_{ij}$, $a_i$, $\bar a_i$, $a_0$ and $a_\infty$ denote the coordinates in $T_{Q_{C_2}\cup_{02} Q_{C_2}}$ of $\big(\mu_{C_2}^{\rot}(\M(\alpha_{120})\big),\alpha_{023})$. 
Note that the coordinates of the elements $h_{12}$, $h_{02}$, and $h_{23}$ involved in \eqref{eq:alphaijkformula} are $(a_{12},a_{21})$, $(a_\infty,a_0)$ and $(a_{23},a_{32})$, respectively. As in Section~\ref{sec:RotProof} we have 
\begin{equation}
\begin{aligned}
\bar a_1&=\frac{a_{12}a_{0}(u_2+u_1u_4)}{a_{21}u_3},&\quad \bar a_2&=a_2'=\frac{a_{0}(u_2^2+u_1^2u_4^2+u_2(u_3^2+2u_1u_4))}{u_3^2u_4}\\
a_1&=v_3a_{23},&\quad a_2&=v_4\frac{a_{32}a_\infty^2}{a_0},
\end{aligned}
\end{equation}
and as in~\eqref{eq:u1u2} we have
\begin{equation}\label{eq:a30a03}
a_{30} =\frac{v_1a_{23}}{a_\infty},\quad a_{03}=\frac{v_2a_{32}}{a_0},\quad a_{01}=\frac{u_1a_\infty}{a_{12}},\quad a_{10}=\frac{u_2a_0}{a_{21}}.
\end{equation}
Using~\eqref{eq:C2Flipzs} we obtain that the face coordinates of 
$\mu_G^{\flip}\big(\mu_G^{\rot}\M(\alpha_{120}),\M(\alpha_{023})\big)$ are


\begin{equation}\label{eq:PrimeCoordinates}
\begin{gathered}
a'_1=a_{23} \Big(u_3+\frac{a_{21}v_3}{a_0}+\frac{a_0 a_{12}^2 u_4 (v_2 + v_3^2 + v_1 v_4)}{a_{21} a_\infty^2 v_3 v_4}\Big),\qquad a'_2=\frac{a_{12}^2 a_{32} u_4}{a_{21}} + \frac{a_{21} a_{32} a_\infty^2 v_4}{a_0^2},\\
\bar a'_1=\frac{a_0 a_{12} a_{23} (u_2 + u_1 u_4) v_1}{a_{21} a_\infty^2 u_3}+ \frac{a_{23} u_1 v_3}{a_{12}},\\
\bar a_2'=a_{32}\Big(\frac{a_0 a_{12}^2 (u_2 + u_1 u_4)^2 v_2}{a_{21}^2 a_\infty^2 u_3^2 v_4}+\frac{(u_2 (u_2 + u_3^2) + 2 u_1 u_2 u_4 + u_1^2 u_4^2) v_2}{a_0 u_3^2 u_4}+\\\frac{(a_{21} a_\infty^2 u_1 u_3 v_3 v_4 + a_0 a_{12}^2 (u_2 + u_1 u_4) (v_2 + v_1 v_4))^2}{a_0 a_{12}^2 a_{21}^2 a_\infty^2 u_3^2 v_3^2 v_4}\Big),
\end{gathered}
\end{equation}
and the non-frozen edge coordinates $a_\infty'$ and $a_0'$ are given by
\begin{equation}
\begin{gathered}
a'_\infty=a_{32}\Big(\frac{v_2}{a_0^2}+\frac{a_{21}a_\infty^2 (u_2 + u_3^2) v_4}{a_0^2 a_{12}^2 u_4}+\frac{2 u_3 (v_2 + v_1 v_4)}{a_0 v_3}+\frac{u_2}{a_{21}}+\\
\frac{a_{12}^2 u_4 (v_2 (v_2 + v_3^2) + 2 v_1 v_2 v_4 + v_1^2 v_4^2)}{a_{21} a_\infty^2 v_3^2 v_4}\Big),\\
a'_0=a_{23}\Big(\frac{a_\infty^2 u_1 + a_{12}^2 v_1}{a_{12} a_\infty^2}+\frac{a_{21} (u_2 + u_3^2 + u_1 u_4) v_3}{a_0 a_{12} u_3 u_4}+\frac{a_0 a_{12} (u_2 + u_1 u_4) (v_2 + v_3^2 + v_1 v_4)}{a_{21} a_\infty^2 u_3v_3 v_4}\Big).
\end{gathered}
\end{equation}

We need to prove the following.

\begin{enumerate}[label=(\roman*)]
\item The non-frozen coordinates of $\M(\alpha_{123})$ are $a_1'$ and $a_2'$.\label{enum:One}
\item The non-frozen coordinates of $\mu_G^{\rot}(\M(\alpha_{130}))$ are $\bar a_1'$ and $\bar a_2'$.\label{enum:Two}
\item The coordinates of $h_{31}=h_{30}l^{-1}$ are $(a_0',a_\infty')$.\label{enum:Three}
\end{enumerate}

To compute $(c,d,l)$ we need a formula for $\iota$. Letting $\iota_I=\iota(-,I,-)$, we have $\iota(u,k,v)=\iota_I(u,kvk^{-1})k$, so we only need a formula for $\iota_I$ (here $I\in H$ is the identity). Applying the algorithm in Section~\ref{sec:PsiPhiIota} we obtain that if  $(x',y',h')=\iota_I(y,x)$, then
\begin{equation}
\begin{aligned}
x'=(\frac{x_1}{\alpha_1},\frac{x_2\alpha_1^2}{\alpha_2},\frac{x_3\alpha_2}{\alpha_1\alpha_3},\frac{x_4\alpha_3^2}{\alpha_2\alpha_4}),\quad 
y'=(\frac{y_1\alpha_4}{\alpha_5},\frac{y_2\alpha_5\alpha_3}{\alpha_4\alpha_6},\frac{y_3\alpha_6^2\alpha_4^2}{\alpha_5\alpha_3^2},\frac{y_4\alpha_3}{\alpha_4\alpha_6}),\quad
h'=(\frac{1}{\alpha_3},\frac{1}{\alpha_4}),
\end{aligned}
\end{equation}
where the $\alpha_i$ are given by
\begin{equation}
\begin{gathered}
\alpha_1=1 + x_1 (y_2 + y_4),\qquad 
\alpha_2=1 + x_2 (y_3 (1 + x_1 y_4)^2 + y_1 (1 + x_1 (y_2 + y_4))^2),\\
\alpha_3=1 + x_3 (y_2 + x_2 y_2 y_3 + y_4) + x_1 (y_2 + y_4 + x_2 x_3 y_2 y_3 y_4),\\
\alpha_4=1 + x_4 \big(y_3 (1 + (x_1 + x_3) y_4)^2 + y_1 (1 + (x_1 + x_3) (y_2 + y_4))^2\big) + 
 x_2 \big(y_3 (1 + x_1 y_4)^2 + \\
    y_1 (1 + x_3^2 x_4 y_2^2 y_3 + 2 x_1 (y_2 + y_4) + x_1^2 (y_2 + y_4)^2)\big)\\
 \alpha_5=1 + x_2 y_3 (1 + x_1 y_4)^2 + x_4 y_3 (1 + (x_1 + x_3) y_4)^2,\qquad
\alpha_6=1 + (x_1 + x_3) y_4.
\end{gathered}
\end{equation}
Using this, together with the explicit formulas for $\Psi\Phi\Psi$, $\Psi\Phi$ and their inverses given in~\eqref{eq:PsiPhiPsiFormula} and \eqref{eq:PsiPhiFormula}, we obtain that $c$, $d$ and $l$ are given by

\begin{equation}
\begin{gathered}
l_1^{-1}=\frac{a_\infty^2 u_1 + a_{12}^2 v_1}{a_{12} a_\infty v_1}   +   \frac{a_{21} a_\infty (u_2 + u_3^2 + u_1 u_4) v_3}{a_0 a_{12} u_3 u_4 v_1}    +   \frac{a_0 a_{12} (u_2 + u_1 u_4) (v_2 + v_3^2 + v_1 v_4)}{a_{21} a_\infty u_3 v_1 v_3 v_4}\\
 l_2^{-1}=\frac{a_{21}}{a_0}+\frac{a_{21}a_\infty^2 (u_2 + u_3^2) v_4}{a_0 a_{12}^2 u_4 v_2}  +   \frac{2 u_3 (v_2 + v_1 v_4)}{v_2 v_3}  +  \frac{a_0u_2}{a_{21}v_2} +\frac{a_0 
   a_{12}^2 u_4 (v_2^2 + v_2v_3^2 + 2 v_1 v_2 v_4 + 
      v_1^2 v_4^2)}{a_{21} a_\infty^2 v_2 v_3^2 v_4},\\
c_1=\frac{a_{12}v_1}{a_\infty l_1},\qquad c_2=\frac{a_{21}v_2}{a_0l_2},\quad
c_3=u_3+\frac{a_{21}v_3}{a_0}+\frac{a_0 a_{12}^2 u_4 (v_2 + v_3^2 + v_1 v_4)}{a_{21} a_\infty^2 v_3 v_4},\quad 
 c_4= u_4 + \frac{a_{21}^2 a_{\infty}^2 v_4}{a_0^2 a_{12}^2},\\
 d_1=\frac{a_\infty u_1}{a_{12}l_1},\qquad d_2=\frac{a_0u_2}{a_{21}l_2},\quad 
 d_3=u_3+\frac{a_{21}a_\infty^2 (u_2 + u_3^2) v_3}{a_0 a_{12}^2 u_4 v_1}+\frac{a_0 u_2 (v_2 + v_3^2 + v_1 v_4)}{a_{21} v_1 v_3 v_4},\\
 d_4=\frac{l_1^2}{l_2}\big(\frac{2 u_3 v_3}{v_1}+\frac{a_{21} a_\infty^2 (u_2 + u_3^2) v_3^2}{a_0 a_{12}^2 u_4 v_1^2}+\frac{a_0 (a_\infty^2 u_2 (v_2 + v_3^2) + a_{12}^2 u_4 v_1^2 v_4)}{a_{21} a_\infty^2 v_1^2 v_4}\big).
 \end{gathered}
 \end{equation}
 By~\eqref{eq:MonomialFormula} the non-frozen coordinates of $\M(\alpha_{123})$ are $c_3a_{23}=a_1'$ and $c_4\frac{a_{32}a_\infty^2}{a_0}=a_2'$ proving~\ref{enum:One}.
Also, by~\eqref{eq:a30a03},
\begin{equation}
h_{31}=h_{30}l^{-1}=(\frac{a_{30}}{l_1},\frac{a_{03}}{l_2})=(\frac{v_1a_{23}}{l_1a_\infty},\frac{v_2a_{32}}{a_0l_2})=(a_0',a_\infty'),
 \end{equation}
 proving \ref{enum:Two}. To prove \ref{enum:Three}, let $b_i$ denote the coordinates of $\M(\alpha_{130})$. Since $h_{13}=w_0(l)h_{30}$, we have
 \begin{equation}
\begin{aligned}
b_{01}&=l_1^{-1}a_{30},& b_{10}&=l_2^{-1}a_{03},&b_{12}&=a_{30},& b_{21}&=a_{03},\\
b_{20}&=a_{01},&b_{02}&=a_{10},&b_1&=d_3b_{12},& b_2&=d_4\frac{b_{21}b_{01}^2}{b_{10}}.
\end{aligned}
 \end{equation}
 Plugging this into the formula \eqref{eq:MurotABC} for $\mu_{C_2}^{\rot}$ we obtain expressions for $b_1'$ and $b_2'$ that equal those of $\bar a_1'$ and $\bar a_2'$ in~\eqref{eq:PrimeCoordinates}. This concludes the proof.

\begin{remark}\label{rm:CoordinateRemark}
Given coordinates on $A_{G,S}$ as described in Section~\ref{sec:CoordinatesAGS} we get a natural cocycle on each triangle of $P$. When $s_G$ is trivial, these glue together to form a natural cocycle on $S$, hence a pair $(\rho,D)$. When $s_G$ is non-trivial, the labelings of identified edges differ by $s_G$ (see Figure~\ref{fig:OrientationChange}), and we instead get a pair $(\overline\rho,\overline D)$. 
\end{remark}

\section{Representations of $3$-manifold groups}
Let $M$ be a compact, oriented $3$-manifold with boundary. Recall that $\rho\colon\pi_1(M)\to G$ is \emph{boundary-unipotent} if peripheral subgroups of $G$ map to conjugates of $N$, and that a \emph{decoration} of such $\rho$ is a  $\rho$-equivariant assignment of a coset $gN\in\A$ to each ideal point (boundary-component) of the universal cover of $M$. A \emph{decorated representation} is a pair $(\rho,D)$, where $\rho$ is boundary-unipotent and $D$ is a decoration of $\rho$. Note that $G$ acts on the set of decorated representations by $g(\rho,D)=(g\rho g^{-1},gD)$. For more details on decorations, we refer to~\cite{GaroufalidisThurstonZickert}. Unless otherwise stated $G$ denotes one of the groups $A_2$, $B_2$, $C_2$ or $G_2$.

\subsection{Generic configurations}
We shall consider a notion of genericity for configurations, which is slightly finer than that of sufficiently generic (Definition~\ref{def:SufGen}). 
\begin{definition}
An element $\alpha\in\Conf_3^*(\A)$ is \emph{generic} if the minor coordinates of $\alpha$, $\rot(\alpha)$, and $\rot^2(\alpha)$ are all non-zero.
\end{definition}
The set $\Conf_3^{\gen}$ of generic configurations in $\Conf_3^*(\A)$ is isomorphic to a Zariski open subset of $T_{Q_G}$. Note that $\mu_G^{\rot}$ is an isomorphism (not just a birational equivalence) on this subset.

\begin{definition} The set $\Conf_4^{\gen}(\A)$ of \emph{generic} configurations in $\Conf_4^*(\A)$ is the largest Zariski open subset $U$ of $\Conf_4^*(\A)$ such that $\Psi_{kl}(U)\in\Conf_3^{\gen}(\A)\times_{kl}\Conf_3^{\gen}(\A)$ for $kl=02$ or $13$, and such that $\mu_G^{\flip}$ defines an isomorphism from $\Psi_{02}(U)$ to $\Psi_{13}(U)$.
\end{definition}

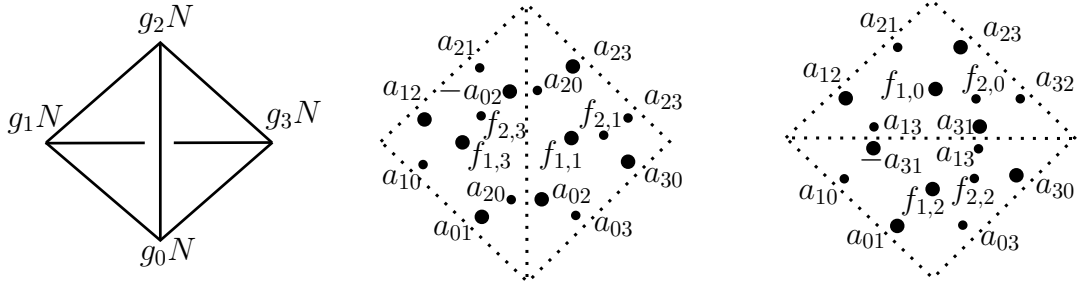
\begin{figure}[htb]
\begin{center}
\scalebox{0.85}{\input{figures_gen/C2Coordinates.tex}}
\caption{Coordinates on the faces of a simplex.}\label{fig:FaceCoordsC2}
\end{center}
\end{figure}

The formulas in Section~\ref{sec:ExplicitMutations} provide explicit defining equations for the variety $\Conf_4^{\gen}(\A)$, and Proposition~\ref{prop:NaturalCocycleFormula} provides an explicit formula for the natural cocycle.
\begin{example}
Let $G=C_2$. Given $\alpha\in\Conf_4^*(\A)$, the simplicial boundary map $\varepsilon_i$ in \eqref{eq:SimpFaceMap} induces configurations on each of the faces with coordinates given by the map $\M$. We denote the coordinates on the $i$th face by $f_{1,i}$ and $f_{2,i}$, and the coordinates on the edges by $a_{ij}$ (see Figure~\ref{fig:FaceCoordsC2}). It now follows from~\eqref{eq:C2Flipzs} that the coordinates satisfy


\begin{equation}\label{eq:C2Simplex}
\begin{aligned}
a_{20} z_1&=a_{32} f_{1,3}^2 + f_{2,1}f_{2,3}, &a_{02} f_{1,2}&=a_{01} f_{1,1}-
 a_{30} f_{1,3},&  f_{1,1}z_2&=-f_{2,1} f_{1,2} + a_{03} a_{23} f_{a,3} ,\\
 f_{2,1} f_{2,2}& = a_{03} z_1 + z_2^2,& f_{1,3} z_3&=z_1 + a_{12} z_2 , &
f_{2,0}f_{2,3}&=a_{10} a_{12}^2 a_{32}  + a_{21} z_1,\\
a_{13} z_1&=f_{2,0}f_{2,2}  + a_{10} z_3^2, &
a_{31} z_2&=a_{23} f_{2,2}  - f_{1,2}z_3, & f_{1,0}z_3&=a_{12} a_{13} a_{23} + f_{2,0}a_{31} 
\end{aligned}
\end{equation}
The ideal generated by these relations defines the Zariski closure of $\Conf_4^{\gen}(\A)$ for $C_2$.
\end{example}

\subsection{Generic decorations and the Ptolemy variety}\label{sec:PtolemyVariety}
Let $\T$ be a topological ideal triangulation of $M$. We assume for simplicity that the triangulation is \emph{ordered}, i.e.~that we have fixed a vertex ordering of each simplex, which is respected by the face pairings. For more on ordered triangulations, see e.g.~\cite{GaroufalidisGoernerZickert}.

A decorated representation $(\rho,D)$ associates a quadruple of affine flags to each 3-simplex of $\T$. We refer to the collection of such as the \emph{associated configurations}.
\begin{definition}
A decoration of a boundary-unipotent representation $\rho$ is \emph{generic} if the associated configurations are in $\Conf_4^{\gen}(\A)$.
\end{definition}
\begin{remark}Note that this notion depends on $\T$.
\end{remark}

The triangulation $\T$ defines a category $J$ with an object for each $k$-simplex and a morphism for each inclusion of a $k$-simplex in an $l$-simplex. For $k=1,2$, let $\Conf_k^{\gen}(\A)=\Conf_k^*(\A)$.
\begin{definition} The \emph{Ptolemy variety} $P_G(\T)$ is the limit of the functor from $J$ to affine varieties taking a $k$-cell to $\Conf_k^{\gen}(\A)$, and an inclusion onto the $i$th face to the face map $\varepsilon_i$ in \eqref{eq:SimpFaceMap}.
\end{definition}
Informally, the Ptolemy variety is the variety built from copies of $\Conf_4^{\gen}(\A)$ by gluing them together using the gluing pattern determined by the triangulation, i.e.~if two faces are identified, the corresponding configuration spaces are identified as well.
Tautologically, we have a one-to-one correspondence between points in the Ptolemy variety and generically decorated boundary-unipotent representations, i.e.~\eqref{eq:OneToOnePG} holds. The natural cocycle provides an explicit formula for this correspondence.

\begin{remark}
One could also consider a Ptolemy variety by gluing together copies of $\Conf_4^*(\A)$ instead of $\Conf_4^{\gen}$. However, we don't have explicit defining equations for this variety.
\end{remark}

\subsection{Obstruction classes}
As mentioned earlier, there are interesting boundary-unipotent representations in $G/\langle s_G\rangle$ that don't have boundary-unipotent lifts to $G$. The obstruction is a class in $H^2(M,\partial M;\Z/2\Z)=H^2(\widehat M;\Z/2\Z)$, where $\widehat M$ is the space obtained from $M$ by collapsing each boundary component to a point. The theory of obstruction classes developed in \cite{GaroufalidisThurstonZickert} for $\SL(n,\C)$ (see~\cite{PtolemyField} for a summary when $n=2$) has a natural analogue for $G$. 
The theory is an elementary generalization of the $\SL(n,\C)$ case, so we only sketch it.


Fix an ordered triangulation $\T$. This determines a $\Delta$-complex structure (as in Hatcher~\cite{Hatcher}) on $\widehat M$. Let $C^*(\widehat M;\Z/2\Z)$ be the simplicial complex of $\Z/ 2\Z$-valued cochains, and let $\sigma\in C^2(\widehat M;\Z/2\Z)$ be a cocycle. The restriction $\sigma_s$ of $\sigma$ to a 3-simplex $\Delta_s$ of $\T$ is a coboundary, so we may represent $\sigma$ by a collection $\eta_s\in C^1(\Delta_s;\Z/2\Z)$ such that $\delta(\eta_s)=\sigma_s$. Note that if a face $f$ of $\Delta_s$ is identified with a face $f'$ of $\Delta_{s'}$, then $\tau_{f,f'}=\eta_{s\vert f}(\eta_{s'\vert f'})^{-1}$ is a cocycle on $\Delta^2$, a standard simplex canonically identified with $f$ and $f'$. Since every such is a coboundary, it is either trivial, or there exists a unique $j=j_{f,f'}\in\{0,1,2\}$ such that $\tau_{f,f'}$ is the coboundary of the $0$-cochain on $\Delta^2$ taking the $j$th vertex to $-1\in\Z/2\Z$ (see Figure~\ref{fig:ObstructionIdentification}). 
We can now define the Ptolemy variety $P^{\sigma}(\T)$ to be the variety obtained by gluing together copies of $\Conf_4^{\gen}(\A)$ in such a way that if two faces $f$ and $f'$ are identified, the corresponding copies of $\Conf_3^{\gen}(\A)$ are identified, not by the identity, but via the map $\kappa_j=\kappa_{j_{f,f'}}$, replacing $g_jN$ by $g_js_GN$ (see Figure~\ref{fig:ObstructionIdentification}). Note that the effect of $\kappa_j$ on the natural cocycle is to leave all three short edges and the long edge opposite $j$ fixed, and to multiply the two long edges extending to $j$ by $s_G$.

One now checks that up to a canonical isomorphism $P^{\sigma}_G(\T)$ only depends on the cohomology class of $\sigma$, and that the set $Z^1(\widehat M;\Z/2\Z)$ of $1$-cocycles acts on $P^\sigma_G(\T)$ with orbits corresponding to decorated boundary-unipotent $G/\langle s_G\rangle$-representations. This proves~\eqref{eq:OneToOnePGsigma}.

\begin{figure}[htb]
\begin{center}
\begin{minipage}[c]{0.48\textwidth}
\scalebox{0.75}{\input{figures_gen/ObstructionIdentification.tex}}
\end{minipage}
\begin{minipage}[c]{0.47\textwidth}
\scalebox{0.75}{\input{figures_gen/Kappa.tex}}
\end{minipage}
\\
\begin{minipage}[t]{0.48\textwidth}
\caption{$\eta_s$, $\eta_{s'}$ and $\eta_{s|f}(\eta_{s'|f'})^{-1}$.}\label{fig:ObstructionIdentification}
\end{minipage}
\begin{minipage}[t]{0.47\textwidth}
\caption{Effect of $\kappa_2$ on the coordinates, $G=C_2$.}\label{fig:Kappa}
\end{minipage}
\end{center}
\end{figure}

\subsection{The reduced Ptolemy variety}
 
If $M$ has a single boundary component, the action of $H$ on $\Conf_k^{\gen}(\A)$ where $h\in H$ acts by replacing each coset $g_iN$ by $g_ihN$ descends to an action on $P_G(\T)$. More generally, if $M$ has $c$ boundary-components, we get an action by $H^c$. This action is also defined for $P_G^{\sigma}(\T)$. We refer to the quotients as \emph{reduced Ptolemy varieties}, and we denote the quotients by $P_G(\T)_{\red}$ and $P_G^{\sigma}(\T)_{\red}$, respectively.

\subsection{Explicit computations for the figure eight not complement}\label{sec:Fig8}
\begin{figure}[htb]
\begin{center}
\scalebox{1}{\input{figures_gen/Fig8Trig.tex}}
\caption{Ordered triangulation of $M$.}\label{fig:Fig8Trig}
\end{center}
\end{figure}
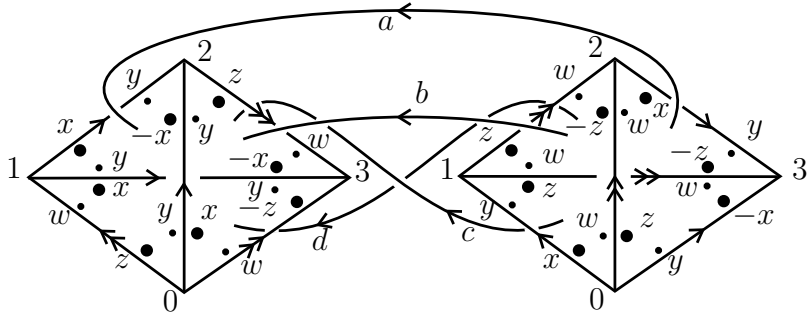

Let $M$ be the figure eight knot complement, and let $\T$ be the standard ideal triangulation of $M$ with two ideal simplices. Figure~\ref{fig:Fig8Trig} shows this triangulation together with the edge coordinates for $G=C_2$ (the face coordinates are not shown).

Using the explicit relations in \eqref{eq:C2Simplex} we obtain that the Zariski closure of $P_{C_2}(\T)$ is given by
\begin{equation}
\begin{aligned}
&f_{2,1}f_{2,3} + f_{1,3}^2w - z_{1,0}y, \quad&
&f_{1,2}x - f_{1,1}z - f_{1,3}z,\quad&
&z_{2,0}f_{1,1} + f_{2,1}f_{1,2} - f_{1,3}wz,\\
&z_{2,0}^2 - f_{2,1}f_{2,2} + z_{1,0}w,&
&z_{1,0} - z_{3,0}f_{1,3} + z_{2,0}x,\quad&
&f_{2,0}f_{2,3} - w^2x^2 - z_{1,0}y,\\
&f_{2,0}f_{2,2} + z_{3,0}^2w - z_{1,0}y,&
&z_{3,0}f_{1,2} - z_{2,0}x - f_{2,2}z,&
&z_{3,0}f_{1,0} + f_{2,0}x - xyz,\\
&f_{2,1}f_{2,3} - z_{1,1}w + f_{1,1}^2y,\quad &
&xf_{1,1} + f_{1,3}x - f_{1,0}z,&
&z_{2,1}f_{1,3} + f_{1,0}f_{2,3} - f_{1,1}xy,\\
&z_{2,1}^2 - f_{2,0}f_{2,3} + z_{1,1}y,&
&z_{1,1} - z_{3,1}f_{1,1} + z_{2,1}z,&
&f_{2,1}f_{2,2} - z_{1,1}w - y^2z^2,\\
&f_{2,0}f_{2,2} - z_{1,1}w + z_{3,1}^2y,\quad&
&z_{3,1}f_{1,0} - f_{2,0}x - z_{2,1}z,&
&z_{3,1}f_{1,2} + f_{2,2}z - wxz.
\end{aligned}
\end{equation}
A computation using Magma~\cite{Magma} shows that there are no solutions where all coordinates are non-zero, and where all rotations of all faces are well defined. Hence, $P_{C_2}(\T)$ is empty.

A simple computations shows that $H^2(\widehat M;\Z/2\Z)=\Z/2\Z$, and that the generator $\sigma$ is represented by the cocycle taking the faces paired by $b$ and $c$ to $-1$.
The Zariski closure of the Ptolemy variety $P_G^\sigma(\T)$ is given by
\begin{equation}
\begin{aligned}
&f_{2,1}f_{2,3} + f_{1,3}^2w - z_{1,0}y,\quad&
&f_{1,2}x - f_{1,1}z - f_{1,3}z,\quad&
&z_{2,0}f_{1,1} + f_{2,1}f_{1,2} + f_{1,3}wz,\\
&z_{2,0}^2 - f_{2,1}f_{2,2} + z_{1,0}w,\quad&
&z_{1,0} - z_{3,0}f_{1,3} + z_{2,0}x,\quad&
&f_{2,0}f_{2,3} - w^2x^2 - z_{1,0}y,\quad\\
&f_{2,0}f_{2,2} + z_{3,0}^2w - z_{1,0}y,\quad&
&z_{3,0}f_{1,2} - z_{2,0}x + f_{2,2}z,\quad&
&z_{3,0}f_{1,0} + f_{2,0}x + xyz,\\
&f_{2,1}f_{2,3} - z_{1,1}w + f_{1,1}^2y,&
&xf_{1,1} + f_{1,3}x - f_{1,0}z,&
&z_{2,1}f_{1,3} - f_{1,0}f_{2,3} + f_{1,1}xy,\\
&z_{2,1}^2 - f_{2,0}f_{2,3} + z_{1,1}y,&
&z_{1,1} + z_{3,1}f_{1,1} + z_{2,1}z,&
&f_{2,1}f_{2,2} - z_{1,1}w - y^2z^2,\\
&f_{2,0}f_{2,2} - z_{1,1}w + z_{3,1}^2y,&
&z_{3,1}f_{1,0} + f_{2,0}x + z_{2,1}z,&
&z_{3,1}f_{1,2} + f_{2,2}z - wxz.
\end{aligned}
\end{equation}
One easily checks that the action by an element $(k_1,k_2)\in H$ multiplies the coordinates $f_{1,i}$ and $f_{2,i}$ by $k_1^2k_2$ and $k_1^2k_2^2$, respectively, so we may add the additional relations $f_{1,0}=1$ and $f_{2,0}=1$ to obtain the reduced Ptolemy variety $P_G^{\sigma}(\T)_{\red}$. A Magma computation shows that there are two zero-dimensional components in $P_G^{\sigma}(\T)_{\red}$. One is defined over $\Q(\sqrt{-3})$ and given by
\begin{equation}
\begin{gathered}
f_{1,0}=f_{2,0}=f_{1,2}=-f_{2,3}=1, \quad f_{1,1}=\frac{1}{2}(-1+\sqrt{-3}),\quad f_{2,1}=-f_{1,2}=-f_{1,3}=\frac{1}{2}(1+\sqrt{-3}),\\
x=\frac{1}{3}(1+\sqrt{-3}),\quad y=\frac{3}{8}(-1+\sqrt{-3}),\quad z=-\frac{1}{3}(1+\sqrt{-3}),\quad w=\frac{3}{4}.
\end{gathered}
\end{equation}
The other component is defined over $\Q(\omega)$, with $\omega$ defined in~\eqref{eq:Omega}, and is given by
\begin{equation}
\begin{gathered}
f_{1,0}=f_{2,0}=1,\quad f_{1,1}=-\frac{3 \omega ^5}{16}+\frac{3 \omega ^4}{8}-\frac{7 \omega ^3}{16}+\frac{7 \omega ^2}{8}-\frac{15 \omega
   }{8}+\frac{3}{2},\\
f_{2,1}=\frac{\omega ^4}{2}-\frac{\omega ^3}{2}+\omega ^2-2 \omega +3,\quad 
f_{1,2}=-\frac{3 \omega ^5}{16}-\frac{\omega ^4}{8}+\frac{\omega ^3}{16}+\frac{3 \omega ^2}{8}-\frac{3 \omega
   }{8}-\frac{3}{2},\\
f_{2,2}=-\frac{\omega ^4}{2}+\frac{\omega ^3}{2}-\omega ^2+2 \omega -3,\quad 
f_{1,3}=\frac{\omega ^5}{16}-\frac{\omega ^4}{8}+\frac{5 \omega ^3}{16}-\frac{\omega ^2}{8}+\frac{5 \omega
   }{8}-\frac{1}{2},\\
x=\frac{3 \omega ^5}{32}-\frac{3 \omega ^4}{16}+\frac{7 \omega ^3}{32}-\frac{11 \omega ^2}{16}+\frac{11 \omega
   }{16}-\frac{1}{4},\quad 
y=-\frac{\omega ^5}{4}+\frac{3 \omega ^4}{8}-\frac{5 \omega ^3}{8}+\frac{3 \omega ^2}{2}-\frac{5 \omega }{4}+1,\\
z=\frac{\omega ^5}{64}+\frac{3 \omega ^4}{32}-\frac{3 \omega ^3}{64}+\frac{3 \omega ^2}{32}-\frac{15 \omega
   }{32}+\frac{9}{8},\quad
w=-\frac{5 \omega ^5}{16}+\frac{\omega ^4}{2}-\frac{11 \omega ^3}{16}+\frac{9 \omega ^2}{8}-\frac{15 \omega
   }{8}+\frac{3}{2}.
\end{gathered}
\end{equation}

\begin{remark}
The reduced Ptolemy variety $P_{B_2}^{\sigma}(\T)_{\red}$ also has two components of degree 2 and 6 defined over $\Q(\sqrt{-3})$ and $\Q(\omega)$, respectively. This is, of course, not surprising since $B_2$ and $C_2$ are isomorphic. We have not been able to explicitly compute the Ptolemy variety for $G_2$.
\end{remark}

\subsection{Recovering the representations}
One can explicitly recover the representation corresponding to a point in the Ptolemy variety using the natural cocycle. 
As described in~\cite[Sec.~4.1]{ZickertEnhancedPtolemy} the fundamental group of the figure eight knot complement has a presentation of the form
\begin{equation}\label{eq:FacePairing}
\langle a,b,c\bigm\vert ca^{-1}bc^{-1}a, ab^{-1}c^{-1}b\rangle,
\end{equation}
where $a$, $b$, and $c$ are the face pairings in Figure~\ref{fig:Fig8Trig}. This presentation is isomorphic to the presentation~\eqref{eq:TwoBridge} via the map taking $x_1$ to $c$ and $x_2$ to $ab^{-1}$.
Let $\alpha_{ij,s}$ and $\beta^i_{jk,s}$ denote the labelings of the natural cocycle associated to simplex $s$. As in~\cite[Sec.~3.5.1]{ZickertEnhancedPtolemy}, the representation is given by
\begin{equation}
a=(\beta^2_{31,0}\alpha_{23,0}\beta^3_{12,0})^{-1},\qquad b=(\beta^3_{01,0})^{-1}(\beta^2_{30,1})^{-1}\alpha_{23,1},\qquad c=\beta^3_{12,0}\beta^3_{12,1}.
\end{equation}
The formulas differ slightly from those of~\cite{ZickertEnhancedPtolemy} due to the fact that we are using an ordered triangulation.
Using the Serre generators for $\mathfrak{sp}(4,\C)$ given in Knapp~\cite{Knapp}, we obtain
\begin{equation}
x_1(t)=\begin{bmatrix}1&t&0&0\\0&1&0&0\\0&0&1&0\\0&0&-t&1\end{bmatrix},\qquad x_2(t)=\begin{bmatrix}1&0&0&0\\0&1&0&t\\0&0&1&0\\0&0&0&1\end{bmatrix},
\end{equation}
and that $h_1^t=\diag(t,t^{-1},t^{-1},t)$ and $h_2^t=\diag(1,t,1,t^{-1})$. Also, $s_{C_2}=-I$, and $\wzero=\Matrix{0}{I}{-I}{0}$.
Using this, we can now recover the natural cocycle explicitly from the coordinates, and we obtain the formulas in Section~\ref{sec:PtolemyIntro}.

\subsection*{Acknowledgment}
The author wishes to thank Matthias Goerner, Stavros Garoufalidis, and Dylan Thurston for helpful comments.

\bibliographystyle{plain}
\bibliography{BibFile}

\end{document}

%% file: figures_gen/QuiverA2Ex.tex
\begingroup
 \setlength{\unitlength}{0.8pt}
 \begin{picture}(176.71072,160.19676)
 \put(0,0){\includegraphics{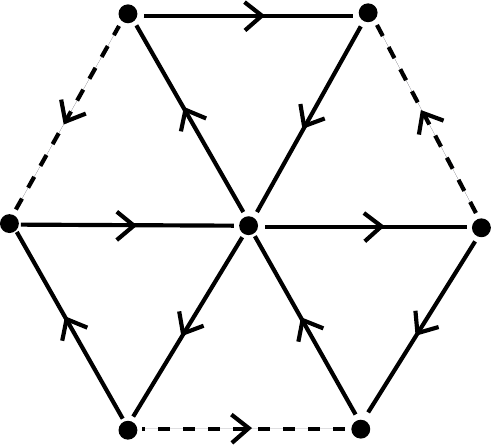}}

\definecolor{inkcol1}{rgb}{0.0,0.0,0.0}
   \put(59.065029,98.422286){\rotatebox{360.0}{\makebox(0,0)[tl]{\strut{}{
    \begin{minipage}[h]{127.655896pt}
\textcolor{inkcol1}{\LARGE{$v_1$}}\\
\end{minipage}}}}}%

 \end{picture}
\endgroup

%% file: figures_gen/QuiverB2Ex.tex
\begingroup
 \setlength{\unitlength}{0.8pt}
 \begin{picture}(182.5,166.07143)
 \put(0,0){\includegraphics{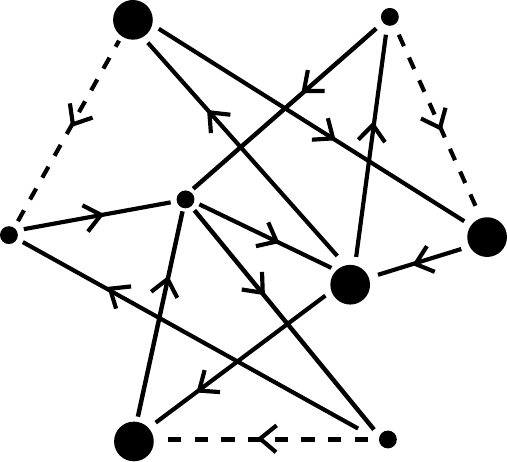}}

\definecolor{inkcol1}{rgb}{0.0,0.0,0.0}
   \put(50.285717,116.887356){\rotatebox{360.0}{\makebox(0,0)[tl]{\strut{}{
    \begin{minipage}[h]{127.655896pt}
\textcolor{inkcol1}{\LARGE{$v_1$}}\\
\end{minipage}}}}}%

\definecolor{inkcol1}{rgb}{0.0,0.0,0.0}
   \put(126.207877,56.143496){\rotatebox{360.0}{\makebox(0,0)[tl]{\strut{}{
    \begin{minipage}[h]{127.655896pt}
\textcolor{inkcol1}{\LARGE{$v_2$}}\\
\end{minipage}}}}}%

 \end{picture}
\endgroup

%% file: figures_gen/QuiverC2Ex.tex
\begingroup
 \setlength{\unitlength}{0.8pt}
 \begin{picture}(178.57143,164.64285)
 \put(0,0){\includegraphics{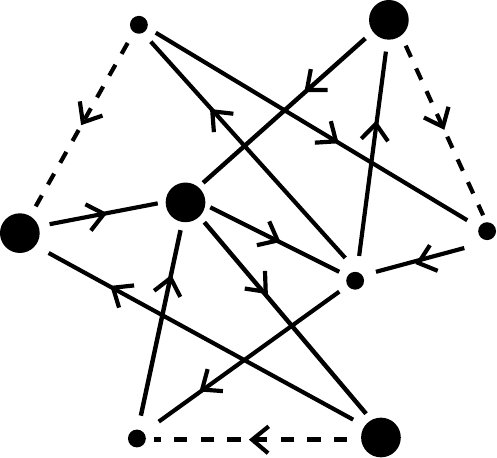}}

\definecolor{inkcol1}{rgb}{0.0,0.0,0.0}
   \put(46.207877,114.357776){\rotatebox{360.0}{\makebox(0,0)[tl]{\strut{}{
    \begin{minipage}[h]{127.655896pt}
\textcolor{inkcol1}{\LARGE{$v_1$}}\\
\end{minipage}}}}}%

\definecolor{inkcol1}{rgb}{0.0,0.0,0.0}
   \put(126.922167,61.143496){\rotatebox{360.0}{\makebox(0,0)[tl]{\strut{}{
    \begin{minipage}[h]{127.655896pt}
\textcolor{inkcol1}{\LARGE{$v_2$}}\\
\end{minipage}}}}}%

 \end{picture}
\endgroup

%% file: figures_gen/QuiverG2Ex.tex
\begingroup
 \setlength{\unitlength}{0.8pt}
 \begin{picture}(178.57141,164.64285)
 \put(0,0){\includegraphics{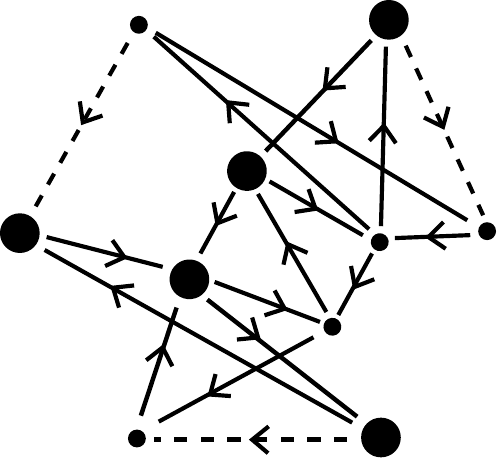}}

\definecolor{inkcol1}{rgb}{0.0,0.0,0.0}
   \put(49.065027,91.857786){\rotatebox{360.0}{\makebox(0,0)[tl]{\strut{}{
    \begin{minipage}[h]{127.655896pt}
\textcolor{inkcol1}{\LARGE{$v_1$}}\\
\end{minipage}}}}}%

\definecolor{inkcol1}{rgb}{0.0,0.0,0.0}
   \put(122.993597,48.286356){\rotatebox{360.0}{\makebox(0,0)[tl]{\strut{}{
    \begin{minipage}[h]{127.655896pt}
\textcolor{inkcol1}{\LARGE{$v_2$}}\\
\end{minipage}}}}}%

\definecolor{inkcol1}{rgb}{0.0,0.0,0.0}
   \put(137.636457,75.429206){\rotatebox{360.0}{\makebox(0,0)[tl]{\strut{}{
    \begin{minipage}[h]{127.655896pt}
\textcolor{inkcol1}{\LARGE{$v_4$}}\\
\end{minipage}}}}}%

\definecolor{inkcol1}{rgb}{0.0,0.0,0.0}
   \put(59.779297,124.000636){\rotatebox{360.0}{\makebox(0,0)[tl]{\strut{}{
    \begin{minipage}[h]{127.655896pt}
\textcolor{inkcol1}{\LARGE{$v_3$}}\\
\end{minipage}}}}}%

 \end{picture}
\endgroup

%% file: figures_gen/MutationExample.tex
\begingroup
 \setlength{\unitlength}{0.8pt}
 \begin{picture}(662.67542,166.07144)
 \put(0,0){\includegraphics{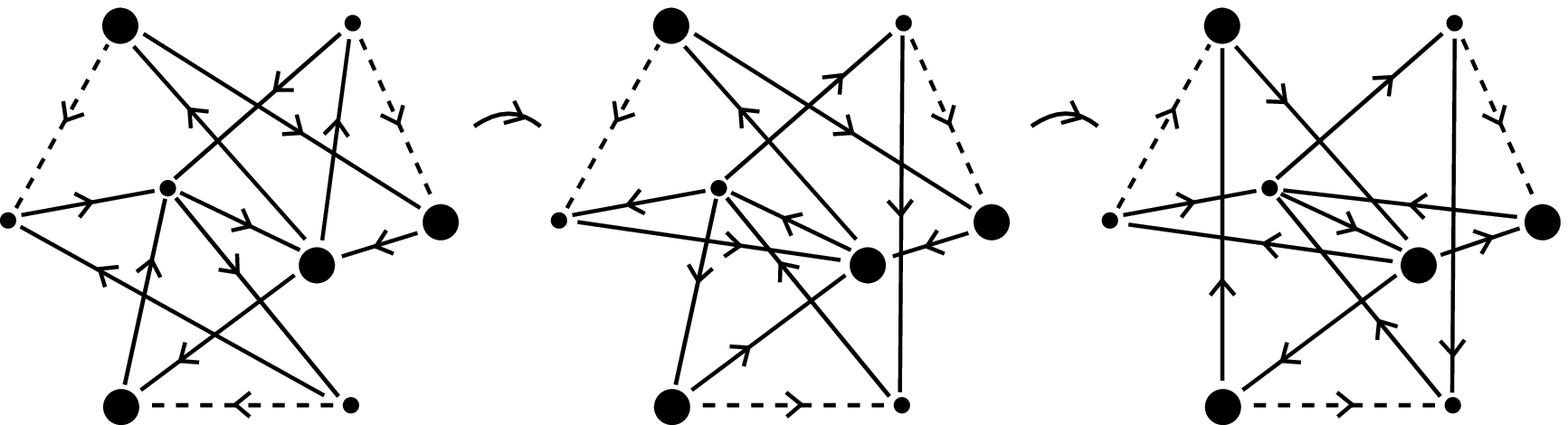}}

\definecolor{inkcol1}{rgb}{0.0,0.0,0.0}
   \put(48.911827,116.678646){\rotatebox{360.0}{\makebox(0,0)[tl]{\strut{}{
    \begin{minipage}[h]{127.655896pt}
\textcolor{inkcol1}{\LARGE{$v_1$}}\\
\end{minipage}}}}}%

\definecolor{inkcol1}{rgb}{0.0,0.0,0.0}
   \put(132.249419,60.110106){\rotatebox{360.0}{\makebox(0,0)[tl]{\strut{}{
    \begin{minipage}[h]{127.655896pt}
\textcolor{inkcol1}{\LARGE{$v_2$}}\\
\end{minipage}}}}}%

\definecolor{inkcol1}{rgb}{0.0,0.0,0.0}
   \put(490.348497,118.193866){\rotatebox{360.0}{\makebox(0,0)[tl]{\strut{}{
    \begin{minipage}[h]{127.655896pt}
\textcolor{inkcol1}{\LARGE{$v_1$}}\\
\end{minipage}}}}}%

\definecolor{inkcol1}{rgb}{0.0,0.0,0.0}
   \put(554.493187,55.059336){\rotatebox{360.0}{\makebox(0,0)[tl]{\strut{}{
    \begin{minipage}[h]{127.655896pt}
\textcolor{inkcol1}{\LARGE{$v_2$}}\\
\end{minipage}}}}}%

\definecolor{inkcol1}{rgb}{0.0,0.0,0.0}
   \put(187.302732,151.023826){\rotatebox{360.0}{\makebox(0,0)[tl]{\strut{}{
    \begin{minipage}[h]{127.655896pt}
\textcolor{inkcol1}{\LARGE{$\mu_{v_1}$}}\\
\end{minipage}}}}}%

\definecolor{inkcol1}{rgb}{0.0,0.0,0.0}
   \put(404.837067,150.162026){\rotatebox{360.0}{\makebox(0,0)[tl]{\strut{}{
    \begin{minipage}[h]{127.655896pt}
\textcolor{inkcol1}{\LARGE{$\mu_{v_2}$}}\\
\end{minipage}}}}}%

\definecolor{inkcol1}{rgb}{0.0,0.0,0.0}
   \put(334.035027,54.430736){\rotatebox{360.0}{\makebox(0,0)[tl]{\strut{}{
    \begin{minipage}[h]{127.655896pt}
\textcolor{inkcol1}{\LARGE{$v_2$}}\\
\end{minipage}}}}}%

\definecolor{inkcol1}{rgb}{0.0,0.0,0.0}
   \put(265.849732,115.039886){\rotatebox{360.0}{\makebox(0,0)[tl]{\strut{}{
    \begin{minipage}[h]{127.655896pt}
\textcolor{inkcol1}{\LARGE{$v_1$}}\\
\end{minipage}}}}}%

 \end{picture}
\endgroup

%% file: figures_gen/MapConf4.tex
\begingroup
 \setlength{\unitlength}{0.8pt}
 \begin{picture}(615.25867,145.03758)
 \put(0,0){\includegraphics{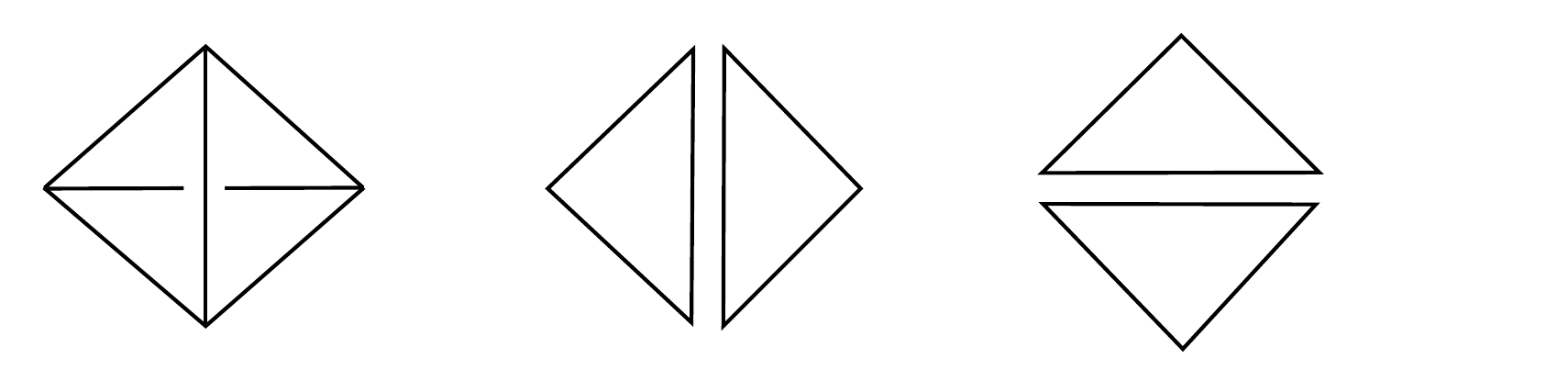}}

\definecolor{inkcol1}{rgb}{0.0,0.0,0.0}
   \put(65.151797,18.580736){\rotatebox{360.0}{\makebox(0,0)[tl]{\strut{}{
    \begin{minipage}[h]{127.655896pt}
\textcolor{inkcol1}{\Large{$g_0N$}}\\
\end{minipage}}}}}%

\definecolor{inkcol1}{rgb}{0.0,0.0,0.0}
   \put(-8.419633,92.866446){\rotatebox{360.0}{\makebox(0,0)[tl]{\strut{}{
    \begin{minipage}[h]{127.655896pt}
\textcolor{inkcol1}{\Large{$g_1N$}}\\
\end{minipage}}}}}%

\definecolor{inkcol1}{rgb}{0.0,0.0,0.0}
   \put(64.080377,148.223586){\rotatebox{360.0}{\makebox(0,0)[tl]{\strut{}{
    \begin{minipage}[h]{127.655896pt}
\textcolor{inkcol1}{\Large{$g_2N$}}\\
\end{minipage}}}}}%

\definecolor{inkcol1}{rgb}{0.0,0.0,0.0}
   \put(135.151807,94.652176){\rotatebox{360.0}{\makebox(0,0)[tl]{\strut{}{
    \begin{minipage}[h]{127.655896pt}
\textcolor{inkcol1}{\Large{$g_3N$}}\\
\end{minipage}}}}}%

\definecolor{inkcol1}{rgb}{0.0,0.0,0.0}
   \put(187.651817,92.152156){\rotatebox{360.0}{\makebox(0,0)[tl]{\strut{}{
    \begin{minipage}[h]{127.655896pt}
\textcolor{inkcol1}{\Large{$g_1N$}}\\
\end{minipage}}}}}%

\definecolor{inkcol1}{rgb}{0.0,0.0,0.0}
   \put(248.008947,144.652156){\rotatebox{360.0}{\makebox(0,0)[tl]{\strut{}{
    \begin{minipage}[h]{127.655896pt}
\textcolor{inkcol1}{\Large{$g_2N$}}\\
\end{minipage}}}}}%

\definecolor{inkcol1}{rgb}{0.0,0.0,0.0}
   \put(289.080377,146.437876){\rotatebox{360.0}{\makebox(0,0)[tl]{\strut{}{
    \begin{minipage}[h]{127.655896pt}
\textcolor{inkcol1}{\Large{$g_2N$}}\\
\end{minipage}}}}}%

\definecolor{inkcol1}{rgb}{0.0,0.0,0.0}
   \put(226.223227,13.937876){\rotatebox{360.0}{\makebox(0,0)[tl]{\strut{}{
    \begin{minipage}[h]{127.655896pt}
\textcolor{inkcol1}{\Large{$g_0s_GN$}}\\
\end{minipage}}}}}%

\definecolor{inkcol1}{rgb}{0.0,0.0,0.0}
   \put(286.937517,18.223586){\rotatebox{360.0}{\makebox(0,0)[tl]{\strut{}{
    \begin{minipage}[h]{127.655896pt}
\textcolor{inkcol1}{\Large{$g_0N$}}\\
\end{minipage}}}}}%

\definecolor{inkcol1}{rgb}{0.0,0.0,0.0}
   \put(382.423847,98.891376){\rotatebox{360.0}{\makebox(0,0)[tl]{\strut{}{
    \begin{minipage}[h]{127.655896pt}
\textcolor{inkcol1}{\Large{$g_1N$}}\\
\end{minipage}}}}}%

\definecolor{inkcol1}{rgb}{0.0,0.0,0.0}
   \put(362.780987,55.319936){\rotatebox{360.0}{\makebox(0,0)[tl]{\strut{}{
    \begin{minipage}[h]{127.655896pt}
\textcolor{inkcol1}{\Large{$g_1s_GN$}}\\
\end{minipage}}}}}%

\definecolor{inkcol1}{rgb}{0.0,0.0,0.0}
   \put(514.209557,97.819946){\rotatebox{360.0}{\makebox(0,0)[tl]{\strut{}{
    \begin{minipage}[h]{127.655896pt}
\textcolor{inkcol1}{\Large{$g_3N$}}\\
\end{minipage}}}}}%

\definecolor{inkcol1}{rgb}{0.0,0.0,0.0}
   \put(514.923837,66.034226){\rotatebox{360.0}{\makebox(0,0)[tl]{\strut{}{
    \begin{minipage}[h]{127.655896pt}
\textcolor{inkcol1}{\Large{$g_3N$}}\\
\end{minipage}}}}}%

\definecolor{inkcol1}{rgb}{0.0,0.0,0.0}
   \put(467.780987,14.962786){\rotatebox{360.0}{\makebox(0,0)[tl]{\strut{}{
    \begin{minipage}[h]{127.655896pt}
\textcolor{inkcol1}{\Large{$g_2N$}}\\
\end{minipage}}}}}%

\definecolor{inkcol1}{rgb}{0.0,0.0,0.0}
   \put(468.138127,151.034226){\rotatebox{360.0}{\makebox(0,0)[tl]{\strut{}{
    \begin{minipage}[h]{127.655896pt}
\textcolor{inkcol1}{\Large{$g_2N$}}\\
\end{minipage}}}}}%

\definecolor{inkcol1}{rgb}{0.0,0.0,0.0}
   \put(339.908497,90.324396){\rotatebox{360.0}{\makebox(0,0)[tl]{\strut{}{
    \begin{minipage}[h]{127.655896pt}
\textcolor{inkcol1}{\Large{$g_3N$}}\\
\end{minipage}}}}}%

 \end{picture}
\endgroup

%% file: figures_gen/QC202.tex
\begingroup
 \setlength{\unitlength}{0.8pt}
 \begin{picture}(453.86667,261.04987)
 \put(0,0){\includegraphics{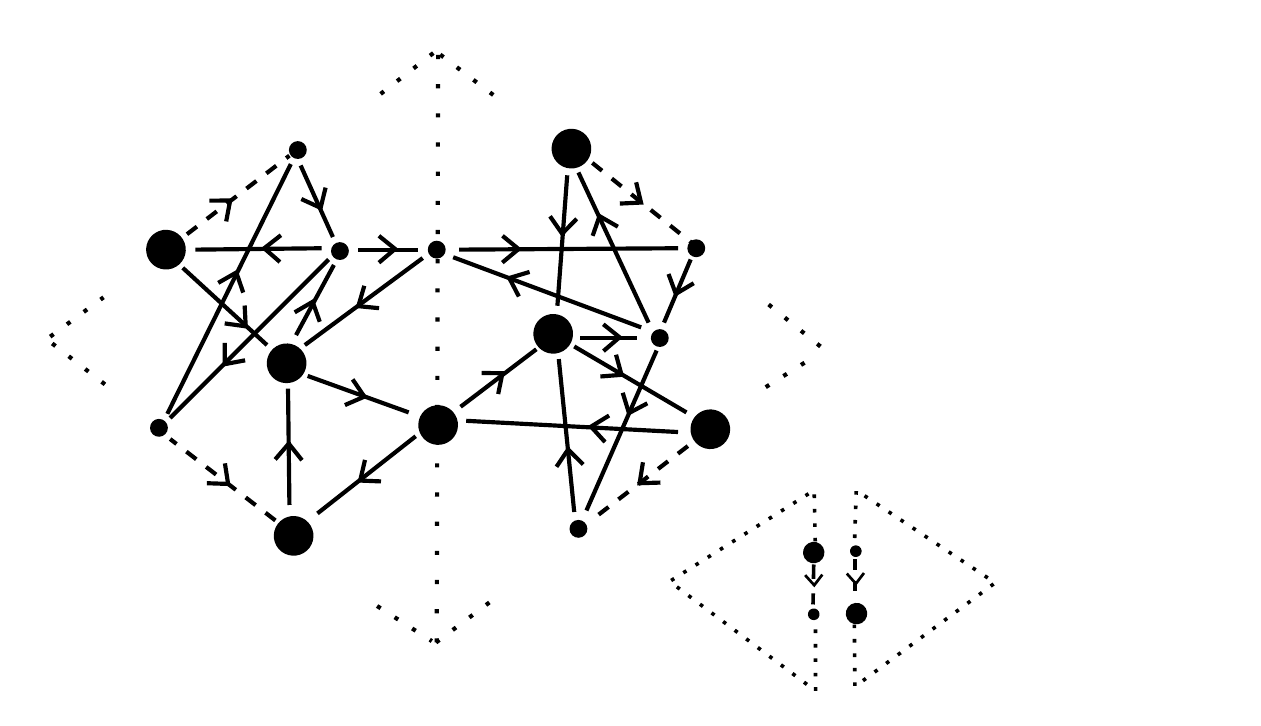}}

\definecolor{inkcol1}{rgb}{0.0,0.0,0.0}
   \put(-0.205343,152.588636){\rotatebox{360.0}{\makebox(0,0)[tl]{\strut{}{
    \begin{minipage}[h]{127.655896pt}
\textcolor{inkcol1}{\Large{$1$}}\\
\end{minipage}}}}}%

\definecolor{inkcol1}{rgb}{0.0,0.0,0.0}
   \put(146.981057,26.038576){\rotatebox{360.0}{\makebox(0,0)[tl]{\strut{}{
    \begin{minipage}[h]{127.655896pt}
\textcolor{inkcol1}{\Large{$0$}}\\
\end{minipage}}}}}%

\definecolor{inkcol1}{rgb}{0.0,0.0,0.0}
   \put(301.263907,148.042956){\rotatebox{360.0}{\makebox(0,0)[tl]{\strut{}{
    \begin{minipage}[h]{127.655896pt}
\textcolor{inkcol1}{\Large{$3$}}\\
\end{minipage}}}}}%

\definecolor{inkcol1}{rgb}{0.0,0.0,0.0}
   \put(282.576087,99.887386){\rotatebox{360.0}{\makebox(0,0)[tl]{\strut{}{
    \begin{minipage}[h]{127.655896pt}
\textcolor{inkcol1}{\Large{$2$}}\\
\end{minipage}}}}}%

\definecolor{inkcol1}{rgb}{0.0,0.0,0.0}
   \put(274.296157,16.391366){\rotatebox{360.0}{\makebox(0,0)[tl]{\strut{}{
    \begin{minipage}[h]{127.655896pt}
\textcolor{inkcol1}{\Large{$0$}}\\
\end{minipage}}}}}%

\definecolor{inkcol1}{rgb}{0.0,0.0,0.0}
   \put(225.563747,58.139366){\rotatebox{360.0}{\makebox(0,0)[tl]{\strut{}{
    \begin{minipage}[h]{127.655896pt}
\textcolor{inkcol1}{\Large{$1$}}\\
\end{minipage}}}}}%

\definecolor{inkcol1}{rgb}{0.0,0.0,0.0}
   \put(357.067407,57.634296){\rotatebox{360.0}{\makebox(0,0)[tl]{\strut{}{
    \begin{minipage}[h]{127.655896pt}
\textcolor{inkcol1}{\Large{$2$}}\\
\end{minipage}}}}}%

\definecolor{inkcol1}{rgb}{0.0,0.0,0.0}
   \put(311.390827,16.182156){\rotatebox{360.0}{\makebox(0,0)[tl]{\strut{}{
    \begin{minipage}[h]{127.655896pt}
\textcolor{inkcol1}{\Large{$0$}}\\
\end{minipage}}}}}%

\definecolor{inkcol1}{rgb}{0.0,0.0,0.0}
   \put(310.053977,98.999256){\rotatebox{360.0}{\makebox(0,0)[tl]{\strut{}{
    \begin{minipage}[h]{127.655896pt}
\textcolor{inkcol1}{\Large{$1$}}\\
\end{minipage}}}}}%

\definecolor{inkcol1}{rgb}{0.0,0.0,0.0}
   \put(310.039099703,42.290714341){\rotatebox{36.207304568065183}{\makebox(0,0)[tl]{\strut{}{
    \begin{minipage}[h]{127.655896pt}
\textcolor{inkcol1}{\Large{$Q_{C_2}$}}\\
\end{minipage}}}}}%

\definecolor{inkcol1}{rgb}{0.0,0.0,0.0}
   \put(270.942272764,29.8958169786){\rotatebox{96.978648576042929}{\makebox(0,0)[tl]{\strut{}{
    \begin{minipage}[h]{127.655896pt}
\textcolor{inkcol1}{\Large{$Q_{C_2}$}}\\
\end{minipage}}}}}%

\definecolor{inkcol1}{rgb}{0.0,0.0,0.0}
   \put(147.001007,264.541716){\rotatebox{360.0}{\makebox(0,0)[tl]{\strut{}{
    \begin{minipage}[h]{127.655896pt}
\textcolor{inkcol1}{\Large{$2$}}\\
\end{minipage}}}}}%

\definecolor{inkcol1}{rgb}{0.0,0.0,0.0}
   \put(161.046267,190.279136){\rotatebox{360.0}{\makebox(0,0)[tl]{\strut{}{
    \begin{minipage}[h]{127.655896pt}
\textcolor{inkcol1}{\Large{$v_0$}}\\
\end{minipage}}}}}%

\definecolor{inkcol1}{rgb}{0.0,0.0,0.0}
   \put(162.622787,103.641366){\rotatebox{360.0}{\makebox(0,0)[tl]{\strut{}{
    \begin{minipage}[h]{127.655896pt}
\textcolor{inkcol1}{\Large{$v_\infty$}}\\
\end{minipage}}}}}%

\definecolor{inkcol1}{rgb}{0.0,0.0,0.0}
   \put(172.979927,151.855656){\rotatebox{360.0}{\makebox(0,0)[tl]{\strut{}{
    \begin{minipage}[h]{127.655896pt}
\textcolor{inkcol1}{\Large{$v_1$}}\\
\end{minipage}}}}}%

\definecolor{inkcol1}{rgb}{0.0,0.0,0.0}
   \put(241.551347,143.641366){\rotatebox{360.0}{\makebox(0,0)[tl]{\strut{}{
    \begin{minipage}[h]{127.655896pt}
\textcolor{inkcol1}{\Large{$v_2$}}\\
\end{minipage}}}}}%

\definecolor{inkcol1}{rgb}{0.0,0.0,0.0}
   \put(77.979917,123.641356){\rotatebox{360.0}{\makebox(0,0)[tl]{\strut{}{
    \begin{minipage}[h]{127.655896pt}
\textcolor{inkcol1}{\Large{$\bar v_1$}}\\
\end{minipage}}}}}%

\definecolor{inkcol1}{rgb}{0.0,0.0,0.0}
   \put(119.051357,191.141356){\rotatebox{360.0}{\makebox(0,0)[tl]{\strut{}{
    \begin{minipage}[h]{127.655896pt}
\textcolor{inkcol1}{\Large{$\bar v_2$}}\\
\end{minipage}}}}}%

\definecolor{inkcol1}{rgb}{0.0,0.0,0.0}
   \put(30.107587,103.734896){\rotatebox{360.0}{\makebox(0,0)[tl]{\strut{}{
    \begin{minipage}[h]{127.655896pt}
\textcolor{inkcol1}{\Large{$v_{10}$}}\\
\end{minipage}}}}}%

\definecolor{inkcol1}{rgb}{0.0,0.0,0.0}
   \put(72.028907,64.844026){\rotatebox{360.0}{\makebox(0,0)[tl]{\strut{}{
    \begin{minipage}[h]{127.655896pt}
\textcolor{inkcol1}{\Large{$v_{01}$}}\\
\end{minipage}}}}}%

\definecolor{inkcol1}{rgb}{0.0,0.0,0.0}
   \put(30.026357,191.618156){\rotatebox{360.0}{\makebox(0,0)[tl]{\strut{}{
    \begin{minipage}[h]{127.655896pt}
\textcolor{inkcol1}{\Large{$v_{12}$}}\\
\end{minipage}}}}}%

\definecolor{inkcol1}{rgb}{0.0,0.0,0.0}
   \put(84.069527,226.973506){\rotatebox{360.0}{\makebox(0,0)[tl]{\strut{}{
    \begin{minipage}[h]{127.655896pt}
\textcolor{inkcol1}{\Large{$v_{21}$}}\\
\end{minipage}}}}}%

\definecolor{inkcol1}{rgb}{0.0,0.0,0.0}
   \put(209.146407,225.351526){\rotatebox{360.0}{\makebox(0,0)[tl]{\strut{}{
    \begin{minipage}[h]{127.655896pt}
\textcolor{inkcol1}{\Large{$v_{23}$}}\\
\end{minipage}}}}}%

\definecolor{inkcol1}{rgb}{0.0,0.0,0.0}
   \put(252.055977,189.125636){\rotatebox{360.0}{\makebox(0,0)[tl]{\strut{}{
    \begin{minipage}[h]{127.655896pt}
\textcolor{inkcol1}{\Large{$v_{32}$}}\\
\end{minipage}}}}}%

\definecolor{inkcol1}{rgb}{0.0,0.0,0.0}
   \put(194.896117,64.349896){\rotatebox{360.0}{\makebox(0,0)[tl]{\strut{}{
    \begin{minipage}[h]{127.655896pt}
\textcolor{inkcol1}{\Large{$v_{03}$}}\\
\end{minipage}}}}}%

\definecolor{inkcol1}{rgb}{0.0,0.0,0.0}
   \put(249.734537,97.684926){\rotatebox{360.0}{\makebox(0,0)[tl]{\strut{}{
    \begin{minipage}[h]{127.655896pt}
\textcolor{inkcol1}{\Large{$v_{30}$}}\\
\end{minipage}}}}}%

 \end{picture}
\endgroup

%% file: figures_gen/QC213.tex
\begingroup
 \setlength{\unitlength}{0.8pt}
 \begin{picture}(423.30069,263.20374)
 \put(0,0){\includegraphics{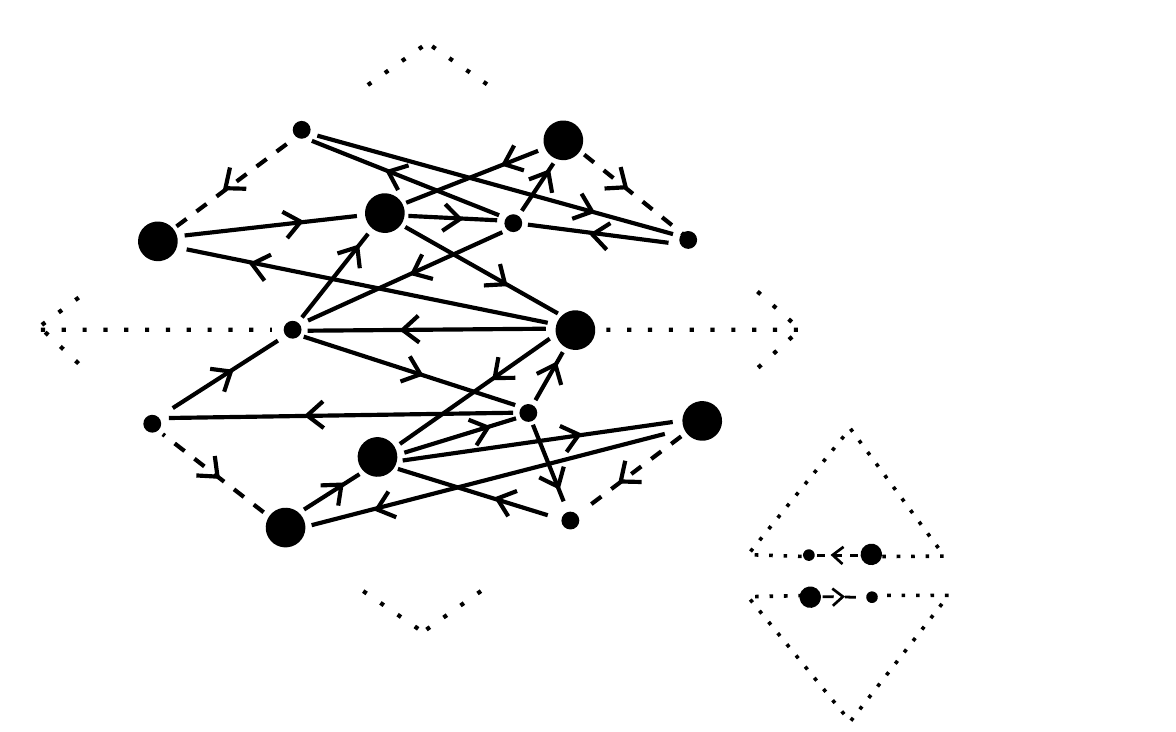}}

\definecolor{inkcol1}{rgb}{0.0,0.0,0.0}
   \put(149.261337,268.175456){\rotatebox{360.0}{\makebox(0,0)[tl]{\strut{}{
    \begin{minipage}[h]{127.655896pt}
\textcolor{inkcol1}{\Large{$2$}}\\
\end{minipage}}}}}%

\definecolor{inkcol1}{rgb}{0.0,0.0,0.0}
   \put(148.792157,36.013696){\rotatebox{360.0}{\makebox(0,0)[tl]{\strut{}{
    \begin{minipage}[h]{127.655896pt}
\textcolor{inkcol1}{\Large{$0$}}\\
\end{minipage}}}}}%

\definecolor{inkcol1}{rgb}{0.0,0.0,0.0}
   \put(1.223227,161.691026){\rotatebox{360.0}{\makebox(0,0)[tl]{\strut{}{
    \begin{minipage}[h]{127.655896pt}
\textcolor{inkcol1}{\Large{$1$}}\\
\end{minipage}}}}}%

\definecolor{inkcol1}{rgb}{0.0,0.0,0.0}
   \put(291.728737,162.823736){\rotatebox{360.0}{\makebox(0,0)[tl]{\strut{}{
    \begin{minipage}[h]{127.655896pt}
\textcolor{inkcol1}{\Large{$3$}}\\
\end{minipage}}}}}%

\definecolor{inkcol1}{rgb}{0.0,0.0,0.0}
   \put(291.129567,86.766436){\rotatebox{360.0}{\makebox(0,0)[tl]{\strut{}{
    \begin{minipage}[h]{127.655896pt}
\textcolor{inkcol1}{\Large{$Q_{C_2}$}}\\
\end{minipage}}}}}%

\definecolor{inkcol1}{rgb}{0.0,0.0,0.0}
   \put(342.323027,75.654746){\rotatebox{360.0}{\makebox(0,0)[tl]{\strut{}{
    \begin{minipage}[h]{127.655896pt}
\textcolor{inkcol1}{\Large{$2$}}\\
\end{minipage}}}}}%

\definecolor{inkcol1}{rgb}{0.0,0.0,0.0}
   \put(254.858167,77.465846){\rotatebox{360.0}{\makebox(0,0)[tl]{\strut{}{
    \begin{minipage}[h]{127.655896pt}
\textcolor{inkcol1}{\Large{$0$}}\\
\end{minipage}}}}}%

\definecolor{inkcol1}{rgb}{0.0,0.0,0.0}
   \put(297.789647,128.825166){\rotatebox{360.0}{\makebox(0,0)[tl]{\strut{}{
    \begin{minipage}[h]{127.655896pt}
\textcolor{inkcol1}{\Large{$1$}}\\
\end{minipage}}}}}%

\definecolor{inkcol1}{rgb}{0.0,0.0,0.0}
   \put(311.558817,13.025296){\rotatebox{360.0}{\makebox(0,0)[tl]{\strut{}{
    \begin{minipage}[h]{127.655896pt}
\textcolor{inkcol1}{\Large{$0$}}\\
\end{minipage}}}}}%

\definecolor{inkcol1}{rgb}{0.0,0.0,0.0}
   \put(251.827707,56.844376){\rotatebox{360.0}{\makebox(0,0)[tl]{\strut{}{
    \begin{minipage}[h]{127.655896pt}
\textcolor{inkcol1}{\Large{$1$}}\\
\end{minipage}}}}}%

\definecolor{inkcol1}{rgb}{0.0,0.0,0.0}
   \put(342.950657,58.359606){\rotatebox{360.0}{\makebox(0,0)[tl]{\strut{}{
    \begin{minipage}[h]{127.655896pt}
\textcolor{inkcol1}{\Large{$2$}}\\
\end{minipage}}}}}%

\definecolor{inkcol1}{rgb}{0.0,0.0,0.0}
   \put(294.754772474,18.5691780593){\rotatebox{57.306553827632627}{\makebox(0,0)[tl]{\strut{}{
    \begin{minipage}[h]{127.655896pt}
\textcolor{inkcol1}{\Large{$Q_{C_2}$}}\\
\end{minipage}}}}}%

\definecolor{inkcol1}{rgb}{0.0,0.0,0.0}
   \put(80.837067,161.498516){\rotatebox{360.0}{\makebox(0,0)[tl]{\strut{}{
    \begin{minipage}[h]{127.655896pt}
\textcolor{inkcol1}{\Large{$v_\infty$}}\\
\end{minipage}}}}}%

\definecolor{inkcol1}{rgb}{0.0,0.0,0.0}
   \put(212.048997,162.631206){\rotatebox{360.0}{\makebox(0,0)[tl]{\strut{}{
    \begin{minipage}[h]{127.655896pt}
\textcolor{inkcol1}{\Large{$v_0$}}\\
\end{minipage}}}}}%

\definecolor{inkcol1}{rgb}{0.0,0.0,0.0}
   \put(113.816767,204.355656){\rotatebox{360.0}{\makebox(0,0)[tl]{\strut{}{
    \begin{minipage}[h]{127.655896pt}
\textcolor{inkcol1}{\Large{$v_1$}}\\
\end{minipage}}}}}%

\definecolor{inkcol1}{rgb}{0.0,0.0,0.0}
   \put(183.998777,178.468036){\rotatebox{360.0}{\makebox(0,0)[tl]{\strut{}{
    \begin{minipage}[h]{127.655896pt}
\textcolor{inkcol1}{\Large{$v_2$}}\\
\end{minipage}}}}}%

\definecolor{inkcol1}{rgb}{0.0,0.0,0.0}
   \put(104.655817,107.421596){\rotatebox{360.0}{\makebox(0,0)[tl]{\strut{}{
    \begin{minipage}[h]{127.655896pt}
\textcolor{inkcol1}{\Large{$\bar v_1$}}\\
\end{minipage}}}}}%

\definecolor{inkcol1}{rgb}{0.0,0.0,0.0}
   \put(198.852537,125.351816){\rotatebox{360.0}{\makebox(0,0)[tl]{\strut{}{
    \begin{minipage}[h]{127.655896pt}
\textcolor{inkcol1}{\Large{$\bar v_2$}}\\
\end{minipage}}}}}%

\definecolor{inkcol1}{rgb}{0.0,0.0,0.0}
   \put(28.763667,199.194306){\rotatebox{360.0}{\makebox(0,0)[tl]{\strut{}{
    \begin{minipage}[h]{127.655896pt}
\textcolor{inkcol1}{\Large{$v_{12}$}}\\
\end{minipage}}}}}%

\definecolor{inkcol1}{rgb}{0.0,0.0,0.0}
   \put(83.816977,235.307256){\rotatebox{360.0}{\makebox(0,0)[tl]{\strut{}{
    \begin{minipage}[h]{127.655896pt}
\textcolor{inkcol1}{\Large{$v_{21}$}}\\
\end{minipage}}}}}%

\definecolor{inkcol1}{rgb}{0.0,0.0,0.0}
   \put(70.684997,68.632096){\rotatebox{360.0}{\makebox(0,0)[tl]{\strut{}{
    \begin{minipage}[h]{127.655896pt}
\textcolor{inkcol1}{\Large{$v_{01}$}}\\
\end{minipage}}}}}%

\definecolor{inkcol1}{rgb}{0.0,0.0,0.0}
   \put(23.712907,110.300886){\rotatebox{360.0}{\makebox(0,0)[tl]{\strut{}{
    \begin{minipage}[h]{127.655896pt}
\textcolor{inkcol1}{\Large{$v_{10}$}}\\
\end{minipage}}}}}%

\definecolor{inkcol1}{rgb}{0.0,0.0,0.0}
   \put(204.025137,72.167626){\rotatebox{360.0}{\makebox(0,0)[tl]{\strut{}{
    \begin{minipage}[h]{127.655896pt}
\textcolor{inkcol1}{\Large{$v_{03}$}}\\
\end{minipage}}}}}%

\definecolor{inkcol1}{rgb}{0.0,0.0,0.0}
   \put(239.885547,100.451886){\rotatebox{360.0}{\makebox(0,0)[tl]{\strut{}{
    \begin{minipage}[h]{127.655896pt}
\textcolor{inkcol1}{\Large{$v_{30}$}}\\
\end{minipage}}}}}%

\definecolor{inkcol1}{rgb}{0.0,0.0,0.0}
   \put(204.277677,234.044566){\rotatebox{360.0}{\makebox(0,0)[tl]{\strut{}{
    \begin{minipage}[h]{127.655896pt}
\textcolor{inkcol1}{\Large{$v_{23}$}}\\
\end{minipage}}}}}%

\definecolor{inkcol1}{rgb}{0.0,0.0,0.0}
   \put(247.966767,193.891006){\rotatebox{360.0}{\makebox(0,0)[tl]{\strut{}{
    \begin{minipage}[h]{127.655896pt}
\textcolor{inkcol1}{\Large{$v_{32}$}}\\
\end{minipage}}}}}%

 \end{picture}
\endgroup

%% file: figures_gen/FlipC2.tex
\begingroup
 \setlength{\unitlength}{0.8pt}
 \begin{picture}(516.79193,198.33563)
 \put(0,0){\includegraphics{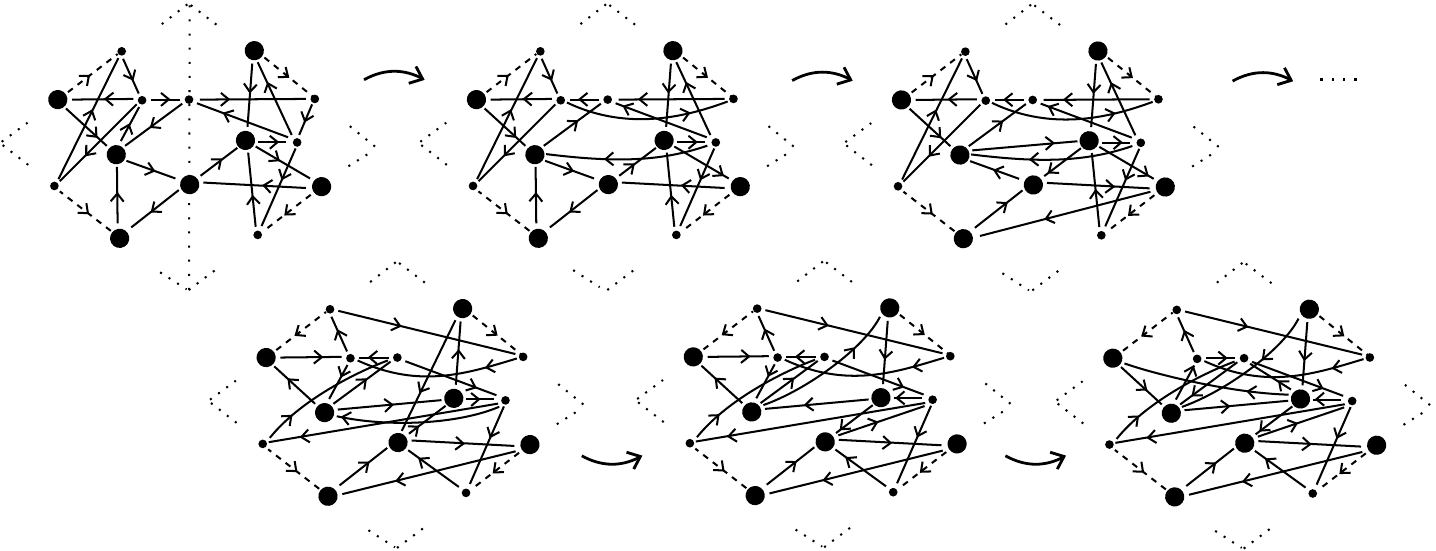}}

\definecolor{inkcol1}{rgb}{0.0,0.0,0.0}
   \put(132.051352,200.324396){\rotatebox{360.0}{\makebox(0,0)[tl]{\strut{}{
    \begin{minipage}[h]{127.655896pt}
\textcolor{inkcol1}{\Large{$\mu_{v_0}$}}\\
\end{minipage}}}}}%

\definecolor{inkcol1}{rgb}{0.0,0.0,0.0}
   \put(280.979922,200.739576){\rotatebox{360.0}{\makebox(0,0)[tl]{\strut{}{
    \begin{minipage}[h]{127.655896pt}
\textcolor{inkcol1}{\Large{$\mu_{v_\infty}$}}\\
\end{minipage}}}}}%

\definecolor{inkcol1}{rgb}{0.0,0.0,0.0}
   \put(205.979932,24.668146){\rotatebox{360.0}{\makebox(0,0)[tl]{\strut{}{
    \begin{minipage}[h]{127.655896pt}
\textcolor{inkcol1}{\Large{$\mu_{v_1}$}}\\
\end{minipage}}}}}%

\definecolor{inkcol1}{rgb}{0.0,0.0,0.0}
   \put(356.337072,25.025286){\rotatebox{360.0}{\makebox(0,0)[tl]{\strut{}{
    \begin{minipage}[h]{127.655896pt}
\textcolor{inkcol1}{\Large{$\mu_{\bar v_1}$}}\\
\end{minipage}}}}}%

\definecolor{inkcol1}{rgb}{0.0,0.0,0.0}
   \put(430.300242,200.427076){\rotatebox{360.0}{\makebox(0,0)[tl]{\strut{}{
    \begin{minipage}[h]{127.655896pt}
\textcolor{inkcol1}{\Large{$\mu_{(1,2,\bar 1,\bar 2,0)}$}}\\
\end{minipage}}}}}%

 \end{picture}
\endgroup

%% file: figures_gen/FundPoly.tex
\begingroup
 \setlength{\unitlength}{0.8pt}
 \begin{picture}(240.42717,188.7)
 \put(0,0){\includegraphics{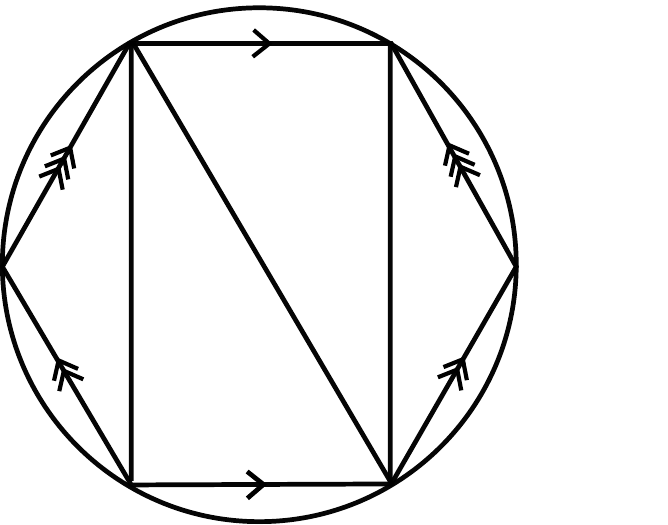}}

\definecolor{inkcol1}{rgb}{0.0,0.0,0.0}
   \put(31.673407,12.971366){\rotatebox{360.0}{\makebox(0,0)[tl]{\strut{}{
    \begin{minipage}[h]{127.655896pt}
\textcolor{inkcol1}{\Large{$4$}}\\
\end{minipage}}}}}%

\definecolor{inkcol1}{rgb}{0.0,0.0,0.0}
   \put(35.959117,189.126496){\rotatebox{360.0}{\makebox(0,0)[tl]{\strut{}{
    \begin{minipage}[h]{127.655896pt}
\textcolor{inkcol1}{\Large{$0$}}\\
\end{minipage}}}}}%

\definecolor{inkcol1}{rgb}{0.0,0.0,0.0}
   \put(140.959117,188.769356){\rotatebox{360.0}{\makebox(0,0)[tl]{\strut{}{
    \begin{minipage}[h]{127.655896pt}
\textcolor{inkcol1}{\Large{$1$}}\\
\end{minipage}}}}}%

\definecolor{inkcol1}{rgb}{0.0,0.0,0.0}
   \put(109.104274273,137.732050047){\rotatebox{304.42299238155897}{\makebox(0,0)[tl]{\strut{}{
    \begin{minipage}[h]{127.655896pt}
\textcolor{inkcol1}{\Large{$Q_G$}}\\
\end{minipage}}}}}%

\definecolor{inkcol1}{rgb}{0.0,0.0,0.0}
   \put(189.394817,103.359766){\rotatebox{360.0}{\makebox(0,0)[tl]{\strut{}{
    \begin{minipage}[h]{127.655896pt}
\textcolor{inkcol1}{\Large{$2$}}\\
\end{minipage}}}}}%

\definecolor{inkcol1}{rgb}{0.0,0.0,0.0}
   \put(140.109117,15.062216){\rotatebox{360.0}{\makebox(0,0)[tl]{\strut{}{
    \begin{minipage}[h]{127.655896pt}
\textcolor{inkcol1}{\Large{$3$}}\\
\end{minipage}}}}}%

\definecolor{inkcol1}{rgb}{0.0,0.0,0.0}
   \put(-16.319463,101.216916){\rotatebox{360.0}{\makebox(0,0)[tl]{\strut{}{
    \begin{minipage}[h]{127.655896pt}
\textcolor{inkcol1}{\Large{$5$}}\\
\end{minipage}}}}}%

\definecolor{inkcol1}{rgb}{0.0,0.0,0.0}
   \put(77.8042145747,81.7734350459){\rotatebox{274.5667478061572}{\makebox(0,0)[tl]{\strut{}{
    \begin{minipage}[h]{127.655896pt}
\textcolor{inkcol1}{\Large{$Q_G$}}\\
\end{minipage}}}}}%

\definecolor{inkcol1}{rgb}{0.0,0.0,0.0}
   \put(168.251720166,107.585147768){\rotatebox{272.13276352358827}{\makebox(0,0)[tl]{\strut{}{
    \begin{minipage}[h]{127.655896pt}
\textcolor{inkcol1}{\Large{$Q_G$}}\\
\end{minipage}}}}}%

\definecolor{inkcol1}{rgb}{0.0,0.0,0.0}
   \put(41.9986503291,113.262308468){\rotatebox{244.05804977772726}{\makebox(0,0)[tl]{\strut{}{
    \begin{minipage}[h]{127.655896pt}
\textcolor{inkcol1}{\Large{$Q_G$}}\\
\end{minipage}}}}}%

 \end{picture}
\endgroup

%% file: figures_gen/FundPolyB2.tex
\begingroup
 \setlength{\unitlength}{0.8pt}
 \begin{picture}(240.5618,188.71875)
 \put(0,0){\includegraphics{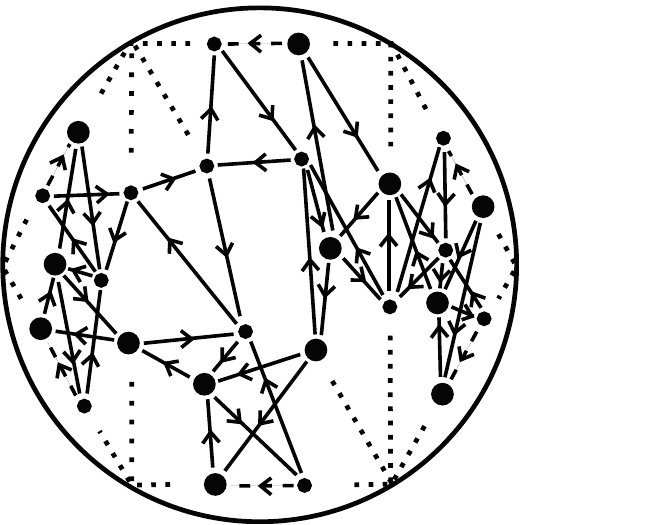}}

\definecolor{inkcol1}{rgb}{0.0,0.0,0.0}
   \put(61.537687,169.169346){\rotatebox{360.0}{\makebox(0,0)[tl]{\strut{}{
    \begin{minipage}[h]{127.655896pt}
\textcolor{inkcol1}{\Large{$a$}}\\
\end{minipage}}}}}%

\definecolor{inkcol1}{rgb}{0.0,0.0,0.0}
   \put(114.930537,171.669356){\rotatebox{360.0}{\makebox(0,0)[tl]{\strut{}{
    \begin{minipage}[h]{127.655896pt}
\textcolor{inkcol1}{\Large{$b$}}\\
\end{minipage}}}}}%

\definecolor{inkcol1}{rgb}{0.0,0.0,0.0}
   \put(57.430537,30.062216){\rotatebox{360.0}{\makebox(0,0)[tl]{\strut{}{
    \begin{minipage}[h]{127.655896pt}
\textcolor{inkcol1}{\Large{$b$}}\\
\end{minipage}}}}}%

\definecolor{inkcol1}{rgb}{0.0,0.0,0.0}
   \put(110.789027,27.919356){\rotatebox{360.0}{\makebox(0,0)[tl]{\strut{}{
    \begin{minipage}[h]{127.655896pt}
\textcolor{inkcol1}{\Large{$a$}}\\
\end{minipage}}}}}%

 \end{picture}
\endgroup

%% file: figures_gen/Epsilonijs.tex
\begingroup
 \setlength{\unitlength}{0.8pt}
 \begin{picture}(704.05853,130.43417)
 \put(0,0){\includegraphics{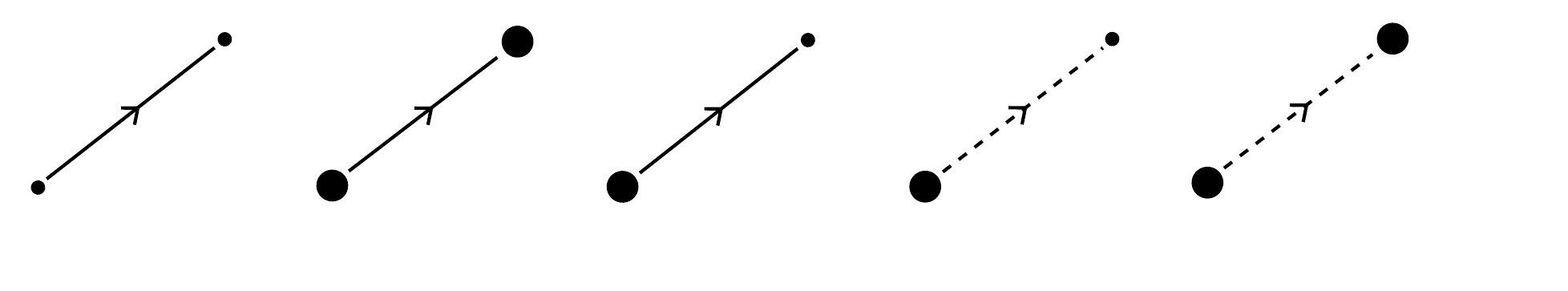}}

\definecolor{inkcol1}{rgb}{0.0,0.0,0.0}
   \put(-0.232683,48.513326){\rotatebox{360.0}{\makebox(0,0)[tl]{\strut{}{
    \begin{minipage}[h]{127.655896pt}
\textcolor{inkcol1}{\LARGE{$i$}}\\
\end{minipage}}}}}%

\definecolor{inkcol1}{rgb}{0.0,0.0,0.0}
   \put(88.476587,136.829896){\rotatebox{360.0}{\makebox(0,0)[tl]{\strut{}{
    \begin{minipage}[h]{127.655896pt}
\textcolor{inkcol1}{\LARGE{$j$}}\\
\end{minipage}}}}}%

\definecolor{inkcol1}{rgb}{0.0,0.0,0.0}
   \put(11.613527,118.209496){\rotatebox{360.0}{\makebox(0,0)[tl]{\strut{}{
    \begin{minipage}[h]{127.655896pt}
\textcolor{inkcol1}{\LARGE{$\varepsilon_{ij}=r$}}\\
\end{minipage}}}}}%

\definecolor{inkcol1}{rgb}{0.0,0.0,0.0}
   \put(39.468877,66.900936){\rotatebox{360.0}{\makebox(0,0)[tl]{\strut{}{
    \begin{minipage}[h]{127.655896pt}
\textcolor{inkcol1}{\LARGE{$\varepsilon_{ji}=-r$}}\\
\end{minipage}}}}}%

\definecolor{inkcol1}{rgb}{0.0,0.0,0.0}
   \put(129.757497,48.092486){\rotatebox{360.0}{\makebox(0,0)[tl]{\strut{}{
    \begin{minipage}[h]{127.655896pt}
\textcolor{inkcol1}{\LARGE{$i$}}\\
\end{minipage}}}}}%

\definecolor{inkcol1}{rgb}{0.0,0.0,0.0}
   \put(220.252461,136.766196){\rotatebox{360.0}{\makebox(0,0)[tl]{\strut{}{
    \begin{minipage}[h]{127.655896pt}
\textcolor{inkcol1}{\LARGE{$j$}}\\
\end{minipage}}}}}%

\definecolor{inkcol1}{rgb}{0.0,0.0,0.0}
   \put(144.103697,118.145806){\rotatebox{360.0}{\makebox(0,0)[tl]{\strut{}{
    \begin{minipage}[h]{127.655896pt}
\textcolor{inkcol1}{\LARGE{$\varepsilon_{ij}=r$}}\\
\end{minipage}}}}}%

\definecolor{inkcol1}{rgb}{0.0,0.0,0.0}
   \put(171.244751,66.837236){\rotatebox{360.0}{\makebox(0,0)[tl]{\strut{}{
    \begin{minipage}[h]{127.655896pt}
\textcolor{inkcol1}{\LARGE{$\varepsilon_{ji}=-r$}}\\
\end{minipage}}}}}%

\definecolor{inkcol1}{rgb}{0.0,0.0,0.0}
   \put(259.5366966,48.153756){\rotatebox{360.0}{\makebox(0,0)[tl]{\strut{}{
    \begin{minipage}[h]{127.655896pt}
\textcolor{inkcol1}{\LARGE{$i$}}\\
\end{minipage}}}}}%

\definecolor{inkcol1}{rgb}{0.0,0.0,0.0}
   \put(350.388824,136.470326){\rotatebox{360.0}{\makebox(0,0)[tl]{\strut{}{
    \begin{minipage}[h]{127.655896pt}
\textcolor{inkcol1}{\LARGE{$j$}}\\
\end{minipage}}}}}%

\definecolor{inkcol1}{rgb}{0.0,0.0,0.0}
   \put(274.5971921,118.207076){\rotatebox{360.0}{\makebox(0,0)[tl]{\strut{}{
    \begin{minipage}[h]{127.655896pt}
\textcolor{inkcol1}{\LARGE{$\varepsilon_{ij}=r$}}\\
\end{minipage}}}}}%

\definecolor{inkcol1}{rgb}{0.0,0.0,0.0}
   \put(301.381115,66.541366){\rotatebox{360.0}{\makebox(0,0)[tl]{\strut{}{
    \begin{minipage}[h]{127.655896pt}
\textcolor{inkcol1}{\LARGE{$\varepsilon_{ji}=-mr$}}\\
\end{minipage}}}}}%

\definecolor{inkcol1}{rgb}{0.0,0.0,0.0}
   \put(394.742637,48.626986){\rotatebox{360.0}{\makebox(0,0)[tl]{\strut{}{
    \begin{minipage}[h]{127.655896pt}
\textcolor{inkcol1}{\LARGE{$i$}}\\
\end{minipage}}}}}%

\definecolor{inkcol1}{rgb}{0.0,0.0,0.0}
   \put(487.023327,136.943556){\rotatebox{360.0}{\makebox(0,0)[tl]{\strut{}{
    \begin{minipage}[h]{127.655896pt}
\textcolor{inkcol1}{\LARGE{$j$}}\\
\end{minipage}}}}}%

\definecolor{inkcol1}{rgb}{0.0,0.0,0.0}
   \put(402.660277,119.037446){\rotatebox{360.0}{\makebox(0,0)[tl]{\strut{}{
    \begin{minipage}[h]{127.655896pt}
\textcolor{inkcol1}{\LARGE{$\varepsilon_{ij}=\frac{r}{2}$}}\\
\end{minipage}}}}}%

\definecolor{inkcol1}{rgb}{0.0,0.0,0.0}
   \put(438.015617,67.014596){\rotatebox{360.0}{\makebox(0,0)[tl]{\strut{}{
    \begin{minipage}[h]{127.655896pt}
\textcolor{inkcol1}{\LARGE{$\varepsilon_{ji}=\frac{-mr}{2}$}}\\
\end{minipage}}}}}%

\definecolor{inkcol1}{rgb}{0.0,0.0,0.0}
   \put(38.283877,92.712146){\rotatebox{360.0}{\makebox(0,0)[tl]{\strut{}{
    \begin{minipage}[h]{127.655896pt}
\textcolor{inkcol1}{\LARGE{$r$}}\\
\end{minipage}}}}}%

\definecolor{inkcol1}{rgb}{0.0,0.0,0.0}
   \put(170.069593,93.426416){\rotatebox{360.0}{\makebox(0,0)[tl]{\strut{}{
    \begin{minipage}[h]{127.655896pt}
\textcolor{inkcol1}{\LARGE{$r$}}\\
\end{minipage}}}}}%

\definecolor{inkcol1}{rgb}{0.0,0.0,0.0}
   \put(301.498161,94.497866){\rotatebox{360.0}{\makebox(0,0)[tl]{\strut{}{
    \begin{minipage}[h]{127.655896pt}
\textcolor{inkcol1}{\LARGE{$r$}}\\
\end{minipage}}}}}%

\definecolor{inkcol1}{rgb}{0.0,0.0,0.0}
   \put(433.504657,93.487706){\rotatebox{360.0}{\makebox(0,0)[tl]{\strut{}{
    \begin{minipage}[h]{127.655896pt}
\textcolor{inkcol1}{\LARGE{$r$}}\\
\end{minipage}}}}}%

\definecolor{inkcol1}{rgb}{0.0,0.0,0.0}
   \put(522.803237,49.413966){\rotatebox{360.0}{\makebox(0,0)[tl]{\strut{}{
    \begin{minipage}[h]{127.655896pt}
\textcolor{inkcol1}{\LARGE{$i$}}\\
\end{minipage}}}}}%

\definecolor{inkcol1}{rgb}{0.0,0.0,0.0}
   \put(609.712447,135.569296){\rotatebox{360.0}{\makebox(0,0)[tl]{\strut{}{
    \begin{minipage}[h]{127.655896pt}
\textcolor{inkcol1}{\LARGE{$j$}}\\
\end{minipage}}}}}%

\definecolor{inkcol1}{rgb}{0.0,0.0,0.0}
   \put(539.234357,120.950436){\rotatebox{360.0}{\makebox(0,0)[tl]{\strut{}{
    \begin{minipage}[h]{127.655896pt}
\textcolor{inkcol1}{\LARGE{$\varepsilon_{ij}=\frac{r}{2}$}}\\
\end{minipage}}}}}%

\definecolor{inkcol1}{rgb}{0.0,0.0,0.0}
   \put(572.091487,67.378996){\rotatebox{360.0}{\makebox(0,0)[tl]{\strut{}{
    \begin{minipage}[h]{127.655896pt}
\textcolor{inkcol1}{\LARGE{$\varepsilon_{ji}=-\frac{r}{2}$}}\\
\end{minipage}}}}}%

\definecolor{inkcol1}{rgb}{0.0,0.0,0.0}
   \put(558.707887,93.649076){\rotatebox{360.0}{\makebox(0,0)[tl]{\strut{}{
    \begin{minipage}[h]{127.655896pt}
\textcolor{inkcol1}{\LARGE{$r$}}\\
\end{minipage}}}}}%

 \end{picture}
\endgroup

%% file: figures_gen/MutationRules.tex
\begingroup
 \setlength{\unitlength}{0.8pt}
 \begin{picture}(748.56567,113.00196)
 \put(0,0){\includegraphics{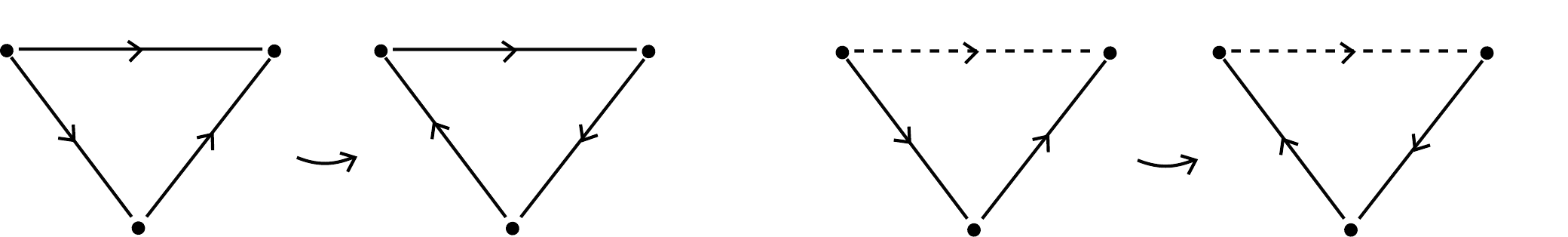}}

\definecolor{inkcol1}{rgb}{0.0,0.0,0.0}
   \put(-10.907383,114.841726){\rotatebox{360.0}{\makebox(0,0)[tl]{\strut{}{
    \begin{minipage}[h]{127.655896pt}
\textcolor{inkcol1}{\LARGE{$d_i$}}\\
\end{minipage}}}}}%

\definecolor{inkcol1}{rgb}{0.0,0.0,0.0}
   \put(59.037737,113.683396){\rotatebox{360.0}{\makebox(0,0)[tl]{\strut{}{
    \begin{minipage}[h]{127.655896pt}
\textcolor{inkcol1}{\LARGE{$t$}}\\
\end{minipage}}}}}%

\definecolor{inkcol1}{rgb}{0.0,0.0,0.0}
   \put(72.279457,12.149776){\rotatebox{360.0}{\makebox(0,0)[tl]{\strut{}{
    \begin{minipage}[h]{127.655896pt}
\textcolor{inkcol1}{\LARGE{$d_k$}}\\
\end{minipage}}}}}%

\definecolor{inkcol1}{rgb}{0.0,0.0,0.0}
   \put(5.499817,59.720646){\rotatebox{360.0}{\makebox(0,0)[tl]{\strut{}{
    \begin{minipage}[h]{127.655896pt}
\textcolor{inkcol1}{\LARGE{$r$}}\\
\end{minipage}}}}}%

\definecolor{inkcol1}{rgb}{0.0,0.0,0.0}
   \put(135.065027,113.310346){\rotatebox{360.0}{\makebox(0,0)[tl]{\strut{}{
    \begin{minipage}[h]{127.655896pt}
\textcolor{inkcol1}{\LARGE{$d_j$}}\\
\end{minipage}}}}}%

\definecolor{inkcol1}{rgb}{0.0,0.0,0.0}
   \put(109.422167,59.720506){\rotatebox{360.0}{\makebox(0,0)[tl]{\strut{}{
    \begin{minipage}[h]{127.655896pt}
\textcolor{inkcol1}{\LARGE{$s$}}\\
\end{minipage}}}}}%

\definecolor{inkcol1}{rgb}{0.0,0.0,0.0}
   \put(167.793107,114.656326){\rotatebox{360.0}{\makebox(0,0)[tl]{\strut{}{
    \begin{minipage}[h]{127.655896pt}
\textcolor{inkcol1}{\LARGE{$d_i$}}\\
\end{minipage}}}}}%

\definecolor{inkcol1}{rgb}{0.0,0.0,0.0}
   \put(203.452507,111.355146){\rotatebox{360.0}{\makebox(0,0)[tl]{\strut{}{
    \begin{minipage}[h]{127.655896pt}
\textcolor{inkcol1}{\LARGE{$t+rs\alpha_{ijk}$}}\\
\end{minipage}}}}}%

\definecolor{inkcol1}{rgb}{0.0,0.0,0.0}
   \put(250.979947,11.964376){\rotatebox{360.0}{\makebox(0,0)[tl]{\strut{}{
    \begin{minipage}[h]{127.655896pt}
\textcolor{inkcol1}{\LARGE{$d_k$}}\\
\end{minipage}}}}}%

\definecolor{inkcol1}{rgb}{0.0,0.0,0.0}
   \put(184.200307,59.535246){\rotatebox{360.0}{\makebox(0,0)[tl]{\strut{}{
    \begin{minipage}[h]{127.655896pt}
\textcolor{inkcol1}{\LARGE{$r$}}\\
\end{minipage}}}}}%

\definecolor{inkcol1}{rgb}{0.0,0.0,0.0}
   \put(313.765517,113.124946){\rotatebox{360.0}{\makebox(0,0)[tl]{\strut{}{
    \begin{minipage}[h]{127.655896pt}
\textcolor{inkcol1}{\LARGE{$d_j$}}\\
\end{minipage}}}}}%

\definecolor{inkcol1}{rgb}{0.0,0.0,0.0}
   \put(288.122657,59.535106){\rotatebox{360.0}{\makebox(0,0)[tl]{\strut{}{
    \begin{minipage}[h]{127.655896pt}
\textcolor{inkcol1}{\LARGE{$s$}}\\
\end{minipage}}}}}%

\definecolor{inkcol1}{rgb}{0.0,0.0,0.0}
   \put(388.150247,113.942036){\rotatebox{360.0}{\makebox(0,0)[tl]{\strut{}{
    \begin{minipage}[h]{127.655896pt}
\textcolor{inkcol1}{\LARGE{$d_i$}}\\
\end{minipage}}}}}%

\definecolor{inkcol1}{rgb}{0.0,0.0,0.0}
   \put(458.095367,112.783706){\rotatebox{360.0}{\makebox(0,0)[tl]{\strut{}{
    \begin{minipage}[h]{127.655896pt}
\textcolor{inkcol1}{\LARGE{$t$}}\\
\end{minipage}}}}}%

\definecolor{inkcol1}{rgb}{0.0,0.0,0.0}
   \put(471.337087,11.250086){\rotatebox{360.0}{\makebox(0,0)[tl]{\strut{}{
    \begin{minipage}[h]{127.655896pt}
\textcolor{inkcol1}{\LARGE{$d_k$}}\\
\end{minipage}}}}}%

\definecolor{inkcol1}{rgb}{0.0,0.0,0.0}
   \put(404.557447,58.820956){\rotatebox{360.0}{\makebox(0,0)[tl]{\strut{}{
    \begin{minipage}[h]{127.655896pt}
\textcolor{inkcol1}{\LARGE{$r$}}\\
\end{minipage}}}}}%

\definecolor{inkcol1}{rgb}{0.0,0.0,0.0}
   \put(534.122657,112.410656){\rotatebox{360.0}{\makebox(0,0)[tl]{\strut{}{
    \begin{minipage}[h]{127.655896pt}
\textcolor{inkcol1}{\LARGE{$d_j$}}\\
\end{minipage}}}}}%

\definecolor{inkcol1}{rgb}{0.0,0.0,0.0}
   \put(508.479797,58.820816){\rotatebox{360.0}{\makebox(0,0)[tl]{\strut{}{
    \begin{minipage}[h]{127.655896pt}
\textcolor{inkcol1}{\LARGE{$s$}}\\
\end{minipage}}}}}%

\definecolor{inkcol1}{rgb}{0.0,0.0,0.0}
   \put(568.150257,113.942036){\rotatebox{360.0}{\makebox(0,0)[tl]{\strut{}{
    \begin{minipage}[h]{127.655896pt}
\textcolor{inkcol1}{\LARGE{$d_i$}}\\
\end{minipage}}}}}%

\definecolor{inkcol1}{rgb}{0.0,0.0,0.0}
   \put(651.337097,11.250086){\rotatebox{360.0}{\makebox(0,0)[tl]{\strut{}{
    \begin{minipage}[h]{127.655896pt}
\textcolor{inkcol1}{\LARGE{$d_k$}}\\
\end{minipage}}}}}%

\definecolor{inkcol1}{rgb}{0.0,0.0,0.0}
   \put(584.557457,58.820956){\rotatebox{360.0}{\makebox(0,0)[tl]{\strut{}{
    \begin{minipage}[h]{127.655896pt}
\textcolor{inkcol1}{\LARGE{$r$}}\\
\end{minipage}}}}}%

\definecolor{inkcol1}{rgb}{0.0,0.0,0.0}
   \put(714.122667,112.410656){\rotatebox{360.0}{\makebox(0,0)[tl]{\strut{}{
    \begin{minipage}[h]{127.655896pt}
\textcolor{inkcol1}{\LARGE{$d_j$}}\\
\end{minipage}}}}}%

\definecolor{inkcol1}{rgb}{0.0,0.0,0.0}
   \put(688.479807,58.820816){\rotatebox{360.0}{\makebox(0,0)[tl]{\strut{}{
    \begin{minipage}[h]{127.655896pt}
\textcolor{inkcol1}{\LARGE{$s$}}\\
\end{minipage}}}}}%

\definecolor{inkcol1}{rgb}{0.0,0.0,0.0}
   \put(137.636457,29.720506){\rotatebox{360.0}{\makebox(0,0)[tl]{\strut{}{
    \begin{minipage}[h]{127.655896pt}
\textcolor{inkcol1}{\LARGE{$\mu_k$}}\\
\end{minipage}}}}}%

\definecolor{inkcol1}{rgb}{0.0,0.0,0.0}
   \put(539.279317,28.643496){\rotatebox{360.0}{\makebox(0,0)[tl]{\strut{}{
    \begin{minipage}[h]{127.655896pt}
\textcolor{inkcol1}{\LARGE{$\mu_k$}}\\
\end{minipage}}}}}%

\definecolor{inkcol1}{rgb}{0.0,0.0,0.0}
   \put(599.894687,110.988926){\rotatebox{360.0}{\makebox(0,0)[tl]{\strut{}{
    \begin{minipage}[h]{127.655896pt}
\textcolor{inkcol1}{\LARGE{$t+2rs\alpha_{ijk}$}}\\
\end{minipage}}}}}%

 \end{picture}
\endgroup

%% file: figures_gen/ExampleB2.tex
\begingroup
 \setlength{\unitlength}{0.8pt}
 \begin{picture}(662.67542,166.07144)
 \put(0,0){\includegraphics{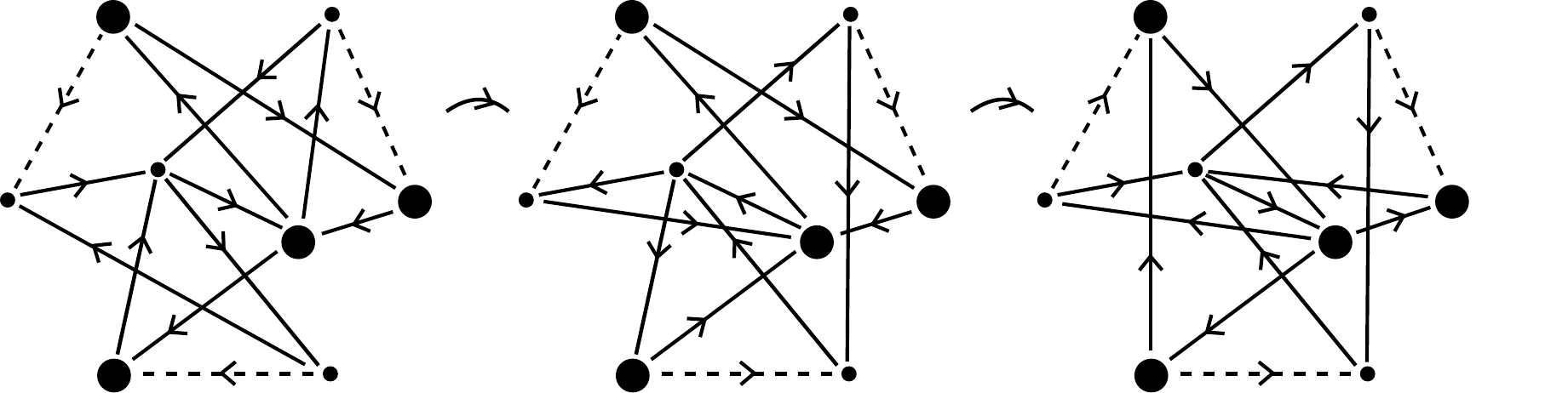}}

\definecolor{inkcol1}{rgb}{0.0,0.0,0.0}
   \put(46.133907,116.678646){\rotatebox{360.0}{\makebox(0,0)[tl]{\strut{}{
    \begin{minipage}[h]{127.655896pt}
\textcolor{inkcol1}{\LARGE{$a_1$}}\\
\end{minipage}}}}}%

\definecolor{inkcol1}{rgb}{0.0,0.0,0.0}
   \put(128.713885,62.130406){\rotatebox{360.0}{\makebox(0,0)[tl]{\strut{}{
    \begin{minipage}[h]{127.655896pt}
\textcolor{inkcol1}{\LARGE{$a_2$}}\\
\end{minipage}}}}}%

\definecolor{inkcol1}{rgb}{0.0,0.0,0.0}
   \put(486.055357,122.992096){\rotatebox{360.0}{\makebox(0,0)[tl]{\strut{}{
    \begin{minipage}[h]{127.655896pt}
\textcolor{inkcol1}{\LARGE{$a_1''$}}\\
\end{minipage}}}}}%

\definecolor{inkcol1}{rgb}{0.0,0.0,0.0}
   \put(550.452577,55.564406){\rotatebox{360.0}{\makebox(0,0)[tl]{\strut{}{
    \begin{minipage}[h]{127.655896pt}
\textcolor{inkcol1}{\LARGE{$a_2''$}}\\
\end{minipage}}}}}%

\definecolor{inkcol1}{rgb}{0.0,0.0,0.0}
   \put(187.302732,151.023826){\rotatebox{360.0}{\makebox(0,0)[tl]{\strut{}{
    \begin{minipage}[h]{127.655896pt}
\textcolor{inkcol1}{\LARGE{$\mu_{v_1}$}}\\
\end{minipage}}}}}%

\definecolor{inkcol1}{rgb}{0.0,0.0,0.0}
   \put(404.837067,150.162026){\rotatebox{360.0}{\makebox(0,0)[tl]{\strut{}{
    \begin{minipage}[h]{127.655896pt}
\textcolor{inkcol1}{\LARGE{$\mu_{v_2}$}}\\
\end{minipage}}}}}%

\definecolor{inkcol1}{rgb}{0.0,0.0,0.0}
   \put(332.519797,55.693426){\rotatebox{360.0}{\makebox(0,0)[tl]{\strut{}{
    \begin{minipage}[h]{127.655896pt}
\textcolor{inkcol1}{\LARGE{$a_2'$}}\\
\end{minipage}}}}}%

\definecolor{inkcol1}{rgb}{0.0,0.0,0.0}
   \put(258.526126,120.595726){\rotatebox{360.0}{\makebox(0,0)[tl]{\strut{}{
    \begin{minipage}[h]{127.655896pt}
\textcolor{inkcol1}{\LARGE{$a_1'$}}\\
\end{minipage}}}}}%

\definecolor{inkcol1}{rgb}{0.0,0.0,0.0}
   \put(-20.023443,79.179476){\rotatebox{360.0}{\makebox(0,0)[tl]{\strut{}{
    \begin{minipage}[h]{127.655896pt}
\textcolor{inkcol1}{\LARGE{$a_{01}$}}\\
\end{minipage}}}}}%

\definecolor{inkcol1}{rgb}{0.0,0.0,0.0}
   \put(9.913117,12.023656){\rotatebox{360.0}{\makebox(0,0)[tl]{\strut{}{
    \begin{minipage}[h]{127.655896pt}
\textcolor{inkcol1}{\LARGE{$a_{02}$}}\\
\end{minipage}}}}}%

\definecolor{inkcol1}{rgb}{0.0,0.0,0.0}
   \put(9.660577,175.415836){\rotatebox{360.0}{\makebox(0,0)[tl]{\strut{}{
    \begin{minipage}[h]{127.655896pt}
\textcolor{inkcol1}{\LARGE{$a_{10}$}}\\
\end{minipage}}}}}%

\definecolor{inkcol1}{rgb}{0.0,0.0,0.0}
   \put(146.536253,175.920916){\rotatebox{360.0}{\makebox(0,0)[tl]{\strut{}{
    \begin{minipage}[h]{127.655896pt}
\textcolor{inkcol1}{\LARGE{$a_{12}$}}\\
\end{minipage}}}}}%

\definecolor{inkcol1}{rgb}{0.0,0.0,0.0}
   \put(163.708846,72.127736){\rotatebox{360.0}{\makebox(0,0)[tl]{\strut{}{
    \begin{minipage}[h]{127.655896pt}
\textcolor{inkcol1}{\LARGE{$a_{21}$}}\\
\end{minipage}}}}}%

\definecolor{inkcol1}{rgb}{0.0,0.0,0.0}
   \put(145.021026,12.023656){\rotatebox{360.0}{\makebox(0,0)[tl]{\strut{}{
    \begin{minipage}[h]{127.655896pt}
\textcolor{inkcol1}{\LARGE{$a_{20}$}}\\
\end{minipage}}}}}%

\definecolor{inkcol1}{rgb}{0.0,0.0,0.0}
   \put(199.3167236,77.683576){\rotatebox{360.0}{\makebox(0,0)[tl]{\strut{}{
    \begin{minipage}[h]{127.655896pt}
\textcolor{inkcol1}{\LARGE{$a_{01}'$}}\\
\end{minipage}}}}}%

\definecolor{inkcol1}{rgb}{0.0,0.0,0.0}
   \put(223.812922,18.842196){\rotatebox{360.0}{\makebox(0,0)[tl]{\strut{}{
    \begin{minipage}[h]{127.655896pt}
\textcolor{inkcol1}{\LARGE{$a_{02}'$}}\\
\end{minipage}}}}}%

\definecolor{inkcol1}{rgb}{0.0,0.0,0.0}
   \put(225.878915,175.668376){\rotatebox{360.0}{\makebox(0,0)[tl]{\strut{}{
    \begin{minipage}[h]{127.655896pt}
\textcolor{inkcol1}{\LARGE{$a_{10}'$}}\\
\end{minipage}}}}}%

\definecolor{inkcol1}{rgb}{0.0,0.0,0.0}
   \put(365.532507,176.678526){\rotatebox{360.0}{\makebox(0,0)[tl]{\strut{}{
    \begin{minipage}[h]{127.655896pt}
\textcolor{inkcol1}{\LARGE{$a_{12}'$}}\\
\end{minipage}}}}}%

\definecolor{inkcol1}{rgb}{0.0,0.0,0.0}
   \put(389.018547,77.178496){\rotatebox{360.0}{\makebox(0,0)[tl]{\strut{}{
    \begin{minipage}[h]{127.655896pt}
\textcolor{inkcol1}{\LARGE{$a_{21}'$}}\\
\end{minipage}}}}}%

\definecolor{inkcol1}{rgb}{0.0,0.0,0.0}
   \put(365.027427,19.852336){\rotatebox{360.0}{\makebox(0,0)[tl]{\strut{}{
    \begin{minipage}[h]{127.655896pt}
\textcolor{inkcol1}{\LARGE{$a_{20}'$}}\\
\end{minipage}}}}}%

\definecolor{inkcol1}{rgb}{0.0,0.0,0.0}
   \put(427.909427,76.168356){\rotatebox{360.0}{\makebox(0,0)[tl]{\strut{}{
    \begin{minipage}[h]{127.655896pt}
\textcolor{inkcol1}{\LARGE{$a_{01}''$}}\\
\end{minipage}}}}}%

\definecolor{inkcol1}{rgb}{0.0,0.0,0.0}
   \put(448.663237,18.084576){\rotatebox{360.0}{\makebox(0,0)[tl]{\strut{}{
    \begin{minipage}[h]{127.655896pt}
\textcolor{inkcol1}{\LARGE{$a_{02}''$}}\\
\end{minipage}}}}}%

\definecolor{inkcol1}{rgb}{0.0,0.0,0.0}
   \put(586.549057,18.842196){\rotatebox{360.0}{\makebox(0,0)[tl]{\strut{}{
    \begin{minipage}[h]{127.655896pt}
\textcolor{inkcol1}{\LARGE{$a_{20}''$}}\\
\end{minipage}}}}}%

\definecolor{inkcol1}{rgb}{0.0,0.0,0.0}
   \put(621.399317,82.481796){\rotatebox{360.0}{\makebox(0,0)[tl]{\strut{}{
    \begin{minipage}[h]{127.655896pt}
\textcolor{inkcol1}{\LARGE{$a_{21}''$}}\\
\end{minipage}}}}}%

\definecolor{inkcol1}{rgb}{0.0,0.0,0.0}
   \put(588.316827,175.415846){\rotatebox{360.0}{\makebox(0,0)[tl]{\strut{}{
    \begin{minipage}[h]{127.655896pt}
\textcolor{inkcol1}{\LARGE{$a_{12}''$}}\\
\end{minipage}}}}}%

\definecolor{inkcol1}{rgb}{0.0,0.0,0.0}
   \put(448.158157,180.466596){\rotatebox{360.0}{\makebox(0,0)[tl]{\strut{}{
    \begin{minipage}[h]{127.655896pt}
\textcolor{inkcol1}{\LARGE{$a_{10}''$}}\\
\end{minipage}}}}}%

 \end{picture}
\endgroup

%% file: figures_gen/EdgeCoordinates.tex
\begingroup
 \setlength{\unitlength}{0.8pt}
 \begin{picture}(382.33902,203.12148)
 \put(0,0){\includegraphics{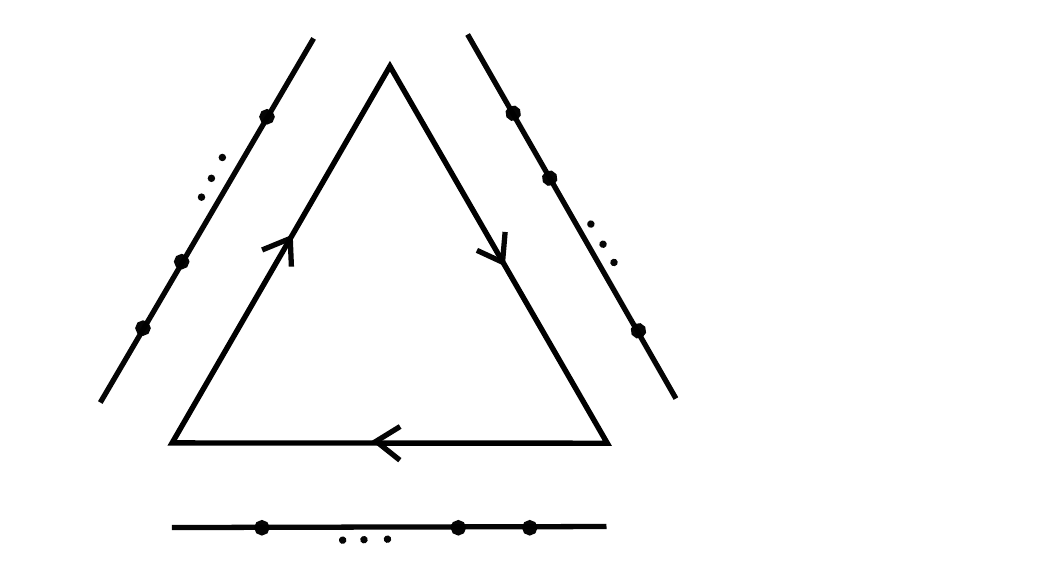}}

\definecolor{inkcol1}{rgb}{0.0,0.0,0.0}
   \put(122.62315,202.464206){\rotatebox{360.0}{\makebox(0,0)[tl]{\strut{}{
    \begin{minipage}[h]{127.655896pt}
\textcolor{inkcol1}{\LARGE{$g_1N$}}\\
\end{minipage}}}}}%

\definecolor{inkcol1}{rgb}{0.0,0.0,0.0}
   \put(23.024288,42.096626){\rotatebox{360.0}{\makebox(0,0)[tl]{\strut{}{
    \begin{minipage}[h]{127.655896pt}
\textcolor{inkcol1}{\LARGE{$g_0N$}}\\
\end{minipage}}}}}%

\definecolor{inkcol1}{rgb}{0.0,0.0,0.0}
   \put(219.09573,48.882356){\rotatebox{360.0}{\makebox(0,0)[tl]{\strut{}{
    \begin{minipage}[h]{127.655896pt}
\textcolor{inkcol1}{\LARGE{$g_2N$}}\\
\end{minipage}}}}}%

\definecolor{inkcol1}{rgb}{0.0,0.0,0.0}
   \put(13.4648587,70.539146){\rotatebox{360.0}{\makebox(0,0)[tl]{\strut{}{
    \begin{minipage}[h]{127.655896pt}
\textcolor{inkcol1}{\Large{$N$}}\\
\end{minipage}}}}}%

\definecolor{inkcol1}{rgb}{0.0,0.0,0.0}
   \put(67.679142,206.949316){\rotatebox{360.0}{\makebox(0,0)[tl]{\strut{}{
    \begin{minipage}[h]{127.655896pt}
\textcolor{inkcol1}{\Large{$\wzero h_1N$}}\\
\end{minipage}}}}}%

\definecolor{inkcol1}{rgb}{0.0,0.0,0.0}
   \put(-18.392285,101.592176){\rotatebox{360.0}{\makebox(0,0)[tl]{\strut{}{
    \begin{minipage}[h]{127.655896pt}
\textcolor{inkcol1}{\Large{$\Delta^{w_1}(h_1)$}}\\
\end{minipage}}}}}%

\definecolor{inkcol1}{rgb}{0.0,0.0,0.0}
   \put(27.861619,177.122296){\rotatebox{360.0}{\makebox(0,0)[tl]{\strut{}{
    \begin{minipage}[h]{127.655896pt}
\textcolor{inkcol1}{\Large{$\Delta^{w_r}(h_1)$}}\\
\end{minipage}}}}}%

\definecolor{inkcol1}{rgb}{0.0,0.0,0.0}
   \put(-3.9240928,126.408006){\rotatebox{360.0}{\makebox(0,0)[tl]{\strut{}{
    \begin{minipage}[h]{127.655896pt}
\textcolor{inkcol1}{\Large{$\Delta^{w_2}(h_1)$}}\\
\end{minipage}}}}}%

\definecolor{inkcol1}{rgb}{0.0,0.0,0.0}
   \put(38.806655,11.219576){\rotatebox{360.0}{\makebox(0,0)[tl]{\strut{}{
    \begin{minipage}[h]{127.655896pt}
\textcolor{inkcol1}{\Large{$N$}}\\
\end{minipage}}}}}%

\definecolor{inkcol1}{rgb}{0.0,0.0,0.0}
   \put(170.94951,206.933876){\rotatebox{360.0}{\makebox(0,0)[tl]{\strut{}{
    \begin{minipage}[h]{127.655896pt}
\textcolor{inkcol1}{\Large{$N$}}\\
\end{minipage}}}}}%

\definecolor{inkcol1}{rgb}{0.0,0.0,0.0}
   \put(219.1638,12.207316){\rotatebox{360.0}{\makebox(0,0)[tl]{\strut{}{
    \begin{minipage}[h]{127.655896pt}
\textcolor{inkcol1}{\Large{$\wzero h_3N$}}\\
\end{minipage}}}}}%

\definecolor{inkcol1}{rgb}{0.0,0.0,0.0}
   \put(244.1638,72.921596){\rotatebox{360.0}{\makebox(0,0)[tl]{\strut{}{
    \begin{minipage}[h]{127.655896pt}
\textcolor{inkcol1}{\Large{$\wzero h_2N$}}\\
\end{minipage}}}}}%

\definecolor{inkcol1}{rgb}{0.0,0.0,0.0}
   \put(186.92119,176.228516){\rotatebox{360.0}{\makebox(0,0)[tl]{\strut{}{
    \begin{minipage}[h]{127.655896pt}
\textcolor{inkcol1}{\Large{$\Delta^{w_1}(h_2)$}}\\
\end{minipage}}}}}%

\definecolor{inkcol1}{rgb}{0.0,0.0,0.0}
   \put(201.56404,154.085656){\rotatebox{360.0}{\makebox(0,0)[tl]{\strut{}{
    \begin{minipage}[h]{127.655896pt}
\textcolor{inkcol1}{\Large{$\Delta^{w_2}(h_2)$}}\\
\end{minipage}}}}}%

\definecolor{inkcol1}{rgb}{0.0,0.0,0.0}
   \put(233.70691,100.157076){\rotatebox{360.0}{\makebox(0,0)[tl]{\strut{}{
    \begin{minipage}[h]{127.655896pt}
\textcolor{inkcol1}{\Large{$\Delta^{w_r}(h_2)$}}\\
\end{minipage}}}}}%

\definecolor{inkcol1}{rgb}{0.0,0.0,0.0}
   \put(159.42118,8.014226){\rotatebox{360.0}{\makebox(0,0)[tl]{\strut{}{
    \begin{minipage}[h]{127.655896pt}
\textcolor{inkcol1}{\Large{$\Delta^{w_1}(h_3)$}}\\
\end{minipage}}}}}%

\definecolor{inkcol1}{rgb}{0.0,0.0,0.0}
   \put(59.064049,7.657096){\rotatebox{360.0}{\makebox(0,0)[tl]{\strut{}{
    \begin{minipage}[h]{127.655896pt}
\textcolor{inkcol1}{\Large{$\Delta^{w_r}(h_3)$}}\\
\end{minipage}}}}}%

\definecolor{inkcol1}{rgb}{0.0,0.0,0.0}
   \put(77.1101773123,126.282223595){\rotatebox{337.54907504206221}{\makebox(0,0)[tl]{\strut{}{
    \begin{minipage}[h]{127.655896pt}
\textcolor{inkcol1}{\Large{$\cong$}}\\
\end{minipage}}}}}%

\definecolor{inkcol1}{rgb}{0.0,0.0,0.0}
   \put(182.118006221,113.128323598){\rotatebox{32.25922341889941}{\makebox(0,0)[tl]{\strut{}{
    \begin{minipage}[h]{127.655896pt}
\textcolor{inkcol1}{\Large{$\cong$}}\\
\end{minipage}}}}}%

\definecolor{inkcol1}{rgb}{0.0,0.0,0.0}
   \put(142.239018461,16.6877536637){\rotatebox{90.582849077661422}{\makebox(0,0)[tl]{\strut{}{
    \begin{minipage}[h]{127.655896pt}
\textcolor{inkcol1}{\Large{$\cong$}}\\
\end{minipage}}}}}%

 \end{picture}
\endgroup

%% file: figures_gen/OrientationChange.tex
\begingroup
 \setlength{\unitlength}{0.8pt}
 \begin{picture}(165.68213,152.27838)
 \put(0,0){\includegraphics{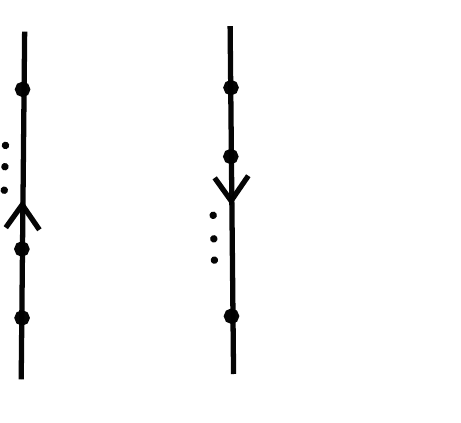}}

\definecolor{inkcol1}{rgb}{0.0,0.0,0.0}
   \put(72.047887,17.853596){\rotatebox{360.0}{\makebox(0,0)[tl]{\strut{}{
    \begin{minipage}[h]{127.655896pt}
\textcolor{inkcol1}{\Large{$g_0s_GN$}}\\
\end{minipage}}}}}%

\definecolor{inkcol1}{rgb}{0.0,0.0,0.0}
   \put(76.088507,163.820656){\rotatebox{360.0}{\makebox(0,0)[tl]{\strut{}{
    \begin{minipage}[h]{127.655896pt}
\textcolor{inkcol1}{\Large{$g_1N$}}\\
\end{minipage}}}}}%

\definecolor{inkcol1}{rgb}{0.0,0.0,0.0}
   \put(-1.693233,17.348536){\rotatebox{360.0}{\makebox(0,0)[tl]{\strut{}{
    \begin{minipage}[h]{127.655896pt}
\textcolor{inkcol1}{\Large{$g_0N$}}\\
\end{minipage}}}}}%

\definecolor{inkcol1}{rgb}{0.0,0.0,0.0}
   \put(0.832147,162.305436){\rotatebox{360.0}{\makebox(0,0)[tl]{\strut{}{
    \begin{minipage}[h]{127.655896pt}
\textcolor{inkcol1}{\Large{$g_1N$}}\\
\end{minipage}}}}}%

\definecolor{inkcol1}{rgb}{0.0,0.0,0.0}
   \put(11.682727,42.729676){\rotatebox{360.0}{\makebox(0,0)[tl]{\strut{}{
    \begin{minipage}[h]{127.655896pt}
\textcolor{inkcol1}{\Large{$x_1$}}\\
\end{minipage}}}}}%

\definecolor{inkcol1}{rgb}{0.0,0.0,0.0}
   \put(11.177647,67.478406){\rotatebox{360.0}{\makebox(0,0)[tl]{\strut{}{
    \begin{minipage}[h]{127.655896pt}
\textcolor{inkcol1}{\Large{$x_2$}}\\
\end{minipage}}}}}%

\definecolor{inkcol1}{rgb}{0.0,0.0,0.0}
   \put(12.152267,126.572336){\rotatebox{360.0}{\makebox(0,0)[tl]{\strut{}{
    \begin{minipage}[h]{127.655896pt}
\textcolor{inkcol1}{\Large{$x_r$}}\\
\end{minipage}}}}}%

\definecolor{inkcol1}{rgb}{0.0,0.0,0.0}
   \put(87.913707,125.057096){\rotatebox{360.0}{\makebox(0,0)[tl]{\strut{}{
    \begin{minipage}[h]{127.655896pt}
\textcolor{inkcol1}{\Large{$x_{\sigma_{G}(1)}$}}\\
\end{minipage}}}}}%

\definecolor{inkcol1}{rgb}{0.0,0.0,0.0}
   \put(87.244647,102.149796){\rotatebox{360.0}{\makebox(0,0)[tl]{\strut{}{
    \begin{minipage}[h]{127.655896pt}
\textcolor{inkcol1}{\Large{$x_{\sigma_{G}(2)}$}}\\
\end{minipage}}}}}%

\definecolor{inkcol1}{rgb}{0.0,0.0,0.0}
   \put(87.749727,42.045716){\rotatebox{360.0}{\makebox(0,0)[tl]{\strut{}{
    \begin{minipage}[h]{127.655896pt}
\textcolor{inkcol1}{\Large{$x_{\sigma_{G}(r)}$}}\\
\end{minipage}}}}}%

 \end{picture}
\endgroup

%% file: figures_gen/MapConf4SL2.tex
\begingroup
 \setlength{\unitlength}{0.8pt}
 \begin{picture}(702.67963,151.01599)
 \put(0,0){\includegraphics{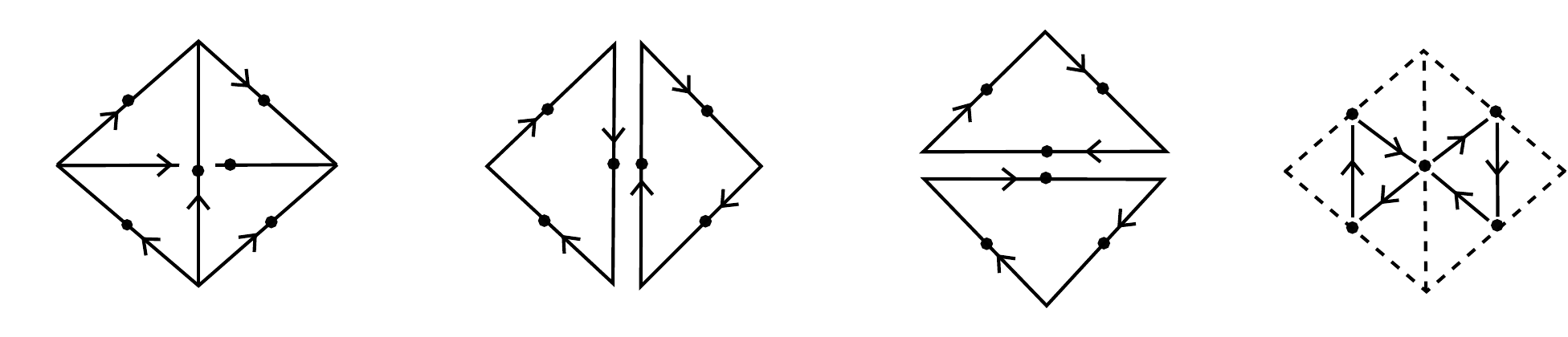}}

\definecolor{inkcol1}{rgb}{0.0,0.0,0.0}
   \put(73.366083,24.248506){\rotatebox{360.0}{\makebox(0,0)[tl]{\strut{}{
    \begin{minipage}[h]{127.655896pt}
\textcolor{inkcol1}{\Large{$g_0N$}}\\
\end{minipage}}}}}%

\definecolor{inkcol1}{rgb}{0.0,0.0,0.0}
   \put(-0.205343,98.534216){\rotatebox{360.0}{\makebox(0,0)[tl]{\strut{}{
    \begin{minipage}[h]{127.655896pt}
\textcolor{inkcol1}{\Large{$g_1N$}}\\
\end{minipage}}}}}%

\definecolor{inkcol1}{rgb}{0.0,0.0,0.0}
   \put(72.294663,153.891356){\rotatebox{360.0}{\makebox(0,0)[tl]{\strut{}{
    \begin{minipage}[h]{127.655896pt}
\textcolor{inkcol1}{\Large{$g_2N$}}\\
\end{minipage}}}}}%

\definecolor{inkcol1}{rgb}{0.0,0.0,0.0}
   \put(143.366093,100.319946){\rotatebox{360.0}{\makebox(0,0)[tl]{\strut{}{
    \begin{minipage}[h]{127.655896pt}
\textcolor{inkcol1}{\Large{$g_3N$}}\\
\end{minipage}}}}}%

\definecolor{inkcol1}{rgb}{0.0,0.0,0.0}
   \put(190.866103,97.462786){\rotatebox{360.0}{\makebox(0,0)[tl]{\strut{}{
    \begin{minipage}[h]{127.655896pt}
\textcolor{inkcol1}{\Large{$g_1N$}}\\
\end{minipage}}}}}%

\definecolor{inkcol1}{rgb}{0.0,0.0,0.0}
   \put(251.223233,149.962786){\rotatebox{360.0}{\makebox(0,0)[tl]{\strut{}{
    \begin{minipage}[h]{127.655896pt}
\textcolor{inkcol1}{\Large{$g_2N$}}\\
\end{minipage}}}}}%

\definecolor{inkcol1}{rgb}{0.0,0.0,0.0}
   \put(292.294663,151.748506){\rotatebox{360.0}{\makebox(0,0)[tl]{\strut{}{
    \begin{minipage}[h]{127.655896pt}
\textcolor{inkcol1}{\Large{$g_2N$}}\\
\end{minipage}}}}}%

\definecolor{inkcol1}{rgb}{0.0,0.0,0.0}
   \put(229.437513,19.248506){\rotatebox{360.0}{\makebox(0,0)[tl]{\strut{}{
    \begin{minipage}[h]{127.655896pt}
\textcolor{inkcol1}{\Large{$g_0s_GN$}}\\
\end{minipage}}}}}%

\definecolor{inkcol1}{rgb}{0.0,0.0,0.0}
   \put(290.151803,23.534216){\rotatebox{360.0}{\makebox(0,0)[tl]{\strut{}{
    \begin{minipage}[h]{127.655896pt}
\textcolor{inkcol1}{\Large{$g_0N$}}\\
\end{minipage}}}}}%

\definecolor{inkcol1}{rgb}{0.0,0.0,0.0}
   \put(387.066703,104.559146){\rotatebox{360.0}{\makebox(0,0)[tl]{\strut{}{
    \begin{minipage}[h]{127.655896pt}
\textcolor{inkcol1}{\Large{$g_1N$}}\\
\end{minipage}}}}}%

\definecolor{inkcol1}{rgb}{0.0,0.0,0.0}
   \put(367.423843,60.987706){\rotatebox{360.0}{\makebox(0,0)[tl]{\strut{}{
    \begin{minipage}[h]{127.655896pt}
\textcolor{inkcol1}{\Large{$g_1s_GN$}}\\
\end{minipage}}}}}%

\definecolor{inkcol1}{rgb}{0.0,0.0,0.0}
   \put(518.852413,103.487716){\rotatebox{360.0}{\makebox(0,0)[tl]{\strut{}{
    \begin{minipage}[h]{127.655896pt}
\textcolor{inkcol1}{\Large{$g_3N$}}\\
\end{minipage}}}}}%

\definecolor{inkcol1}{rgb}{0.0,0.0,0.0}
   \put(519.566693,71.701996){\rotatebox{360.0}{\makebox(0,0)[tl]{\strut{}{
    \begin{minipage}[h]{127.655896pt}
\textcolor{inkcol1}{\Large{$g_3N$}}\\
\end{minipage}}}}}%

\definecolor{inkcol1}{rgb}{0.0,0.0,0.0}
   \put(472.423843,20.630556){\rotatebox{360.0}{\makebox(0,0)[tl]{\strut{}{
    \begin{minipage}[h]{127.655896pt}
\textcolor{inkcol1}{\Large{$g_2N$}}\\
\end{minipage}}}}}%

\definecolor{inkcol1}{rgb}{0.0,0.0,0.0}
   \put(472.780983,156.701996){\rotatebox{360.0}{\makebox(0,0)[tl]{\strut{}{
    \begin{minipage}[h]{127.655896pt}
\textcolor{inkcol1}{\Large{$g_2N$}}\\
\end{minipage}}}}}%

\definecolor{inkcol1}{rgb}{0.0,0.0,0.0}
   \put(343.577093,95.739626){\rotatebox{360.0}{\makebox(0,0)[tl]{\strut{}{
    \begin{minipage}[h]{127.655896pt}
\textcolor{inkcol1}{\Large{$g_3N$}}\\
\end{minipage}}}}}%

\definecolor{inkcol1}{rgb}{0.0,0.0,0.0}
   \put(29.057626,51.026706){\rotatebox{360.0}{\makebox(0,0)[tl]{\strut{}{
    \begin{minipage}[h]{127.655896pt}
\textcolor{inkcol1}{\Large{$c_{01}$}}\\
\end{minipage}}}}}%

\definecolor{inkcol1}{rgb}{0.0,0.0,0.0}
   \put(33.350769,120.222156){\rotatebox{360.0}{\makebox(0,0)[tl]{\strut{}{
    \begin{minipage}[h]{127.655896pt}
\textcolor{inkcol1}{\Large{$c_{12}$}}\\
\end{minipage}}}}}%

\definecolor{inkcol1}{rgb}{0.0,0.0,0.0}
   \put(117.698513,119.969626){\rotatebox{360.0}{\makebox(0,0)[tl]{\strut{}{
    \begin{minipage}[h]{127.655896pt}
\textcolor{inkcol1}{\Large{$c_{23}$}}\\
\end{minipage}}}}}%

\definecolor{inkcol1}{rgb}{0.0,0.0,0.0}
   \put(101.788603,91.685346){\rotatebox{360.0}{\makebox(0,0)[tl]{\strut{}{
    \begin{minipage}[h]{127.655896pt}
\textcolor{inkcol1}{\Large{$c_{13}$}}\\
\end{minipage}}}}}%

\definecolor{inkcol1}{rgb}{0.0,0.0,0.0}
   \put(58.099513,73.250056){\rotatebox{360.0}{\makebox(0,0)[tl]{\strut{}{
    \begin{minipage}[h]{127.655896pt}
\textcolor{inkcol1}{\Large{$c_{02}$}}\\
\end{minipage}}}}}%

\definecolor{inkcol1}{rgb}{0.0,0.0,0.0}
   \put(122.577953,53.400916){\rotatebox{360.0}{\makebox(0,0)[tl]{\strut{}{
    \begin{minipage}[h]{127.655896pt}
\textcolor{inkcol1}{\Large{$c_{03}$}}\\
\end{minipage}}}}}%

\definecolor{inkcol1}{rgb}{0.0,0.0,0.0}
   \put(289.635553,88.835136){\rotatebox{360.0}{\makebox(0,0)[tl]{\strut{}{
    \begin{minipage}[h]{127.655896pt}
\textcolor{inkcol1}{\Large{$c_{02}$}}\\
\end{minipage}}}}}%

\definecolor{inkcol1}{rgb}{0.0,0.0,0.0}
   \put(319.435063,115.856706){\rotatebox{360.0}{\makebox(0,0)[tl]{\strut{}{
    \begin{minipage}[h]{127.655896pt}
\textcolor{inkcol1}{\Large{$c_{23}$}}\\
\end{minipage}}}}}%

\definecolor{inkcol1}{rgb}{0.0,0.0,0.0}
   \put(309.080993,49.944256){\rotatebox{360.0}{\makebox(0,0)[tl]{\strut{}{
    \begin{minipage}[h]{127.655896pt}
\textcolor{inkcol1}{\Large{$-c_{03}$}}\\
\end{minipage}}}}}%

\definecolor{inkcol1}{rgb}{0.0,0.0,0.0}
   \put(205.540363,52.217106){\rotatebox{360.0}{\makebox(0,0)[tl]{\strut{}{
    \begin{minipage}[h]{127.655896pt}
\textcolor{inkcol1}{\Large{$-c_{01}$}}\\
\end{minipage}}}}}%

\definecolor{inkcol1}{rgb}{0.0,0.0,0.0}
   \put(218.924883,119.897326){\rotatebox{360.0}{\makebox(0,0)[tl]{\strut{}{
    \begin{minipage}[h]{127.655896pt}
\textcolor{inkcol1}{\Large{$c_{12}$}}\\
\end{minipage}}}}}%

\definecolor{inkcol1}{rgb}{0.0,0.0,0.0}
   \put(412.621633,124.695556){\rotatebox{360.0}{\makebox(0,0)[tl]{\strut{}{
    \begin{minipage}[h]{127.655896pt}
\textcolor{inkcol1}{\Large{$c_{12}$}}\\
\end{minipage}}}}}%

\definecolor{inkcol1}{rgb}{0.0,0.0,0.0}
   \put(499.494753,122.927786){\rotatebox{360.0}{\makebox(0,0)[tl]{\strut{}{
    \begin{minipage}[h]{127.655896pt}
\textcolor{inkcol1}{\Large{$c_{23}$}}\\
\end{minipage}}}}}%

\definecolor{inkcol1}{rgb}{0.0,0.0,0.0}
   \put(445.704123,103.734886){\rotatebox{360.0}{\makebox(0,0)[tl]{\strut{}{
    \begin{minipage}[h]{127.655896pt}
\textcolor{inkcol1}{\Large{$-c_{13}$}}\\
\end{minipage}}}}}%

\definecolor{inkcol1}{rgb}{0.0,0.0,0.0}
   \put(445.956653,68.127006){\rotatebox{360.0}{\makebox(0,0)[tl]{\strut{}{
    \begin{minipage}[h]{127.655896pt}
\textcolor{inkcol1}{\Large{$-c_{13}$}}\\
\end{minipage}}}}}%

\definecolor{inkcol1}{rgb}{0.0,0.0,0.0}
   \put(498.232073,43.883346){\rotatebox{360.0}{\makebox(0,0)[tl]{\strut{}{
    \begin{minipage}[h]{127.655896pt}
\textcolor{inkcol1}{\Large{$-c_{03}$}}\\
\end{minipage}}}}}%

\definecolor{inkcol1}{rgb}{0.0,0.0,0.0}
   \put(401.509963,46.408726){\rotatebox{360.0}{\makebox(0,0)[tl]{\strut{}{
    \begin{minipage}[h]{127.655896pt}
\textcolor{inkcol1}{\Large{$-c_{01}$}}\\
\end{minipage}}}}}%

\definecolor{inkcol1}{rgb}{0.0,0.0,0.0}
   \put(243.421083,87.572446){\rotatebox{360.0}{\makebox(0,0)[tl]{\strut{}{
    \begin{minipage}[h]{127.655896pt}
\textcolor{inkcol1}{\Large{$c_{02}$}}\\
\end{minipage}}}}}%

 \end{picture}
\endgroup

%% file: figures_gen/NaturalCocycle3d.tex
\begingroup
 \setlength{\unitlength}{0.8pt}
 \begin{picture}(234.20671,190.59085)
 \put(0,0){\includegraphics{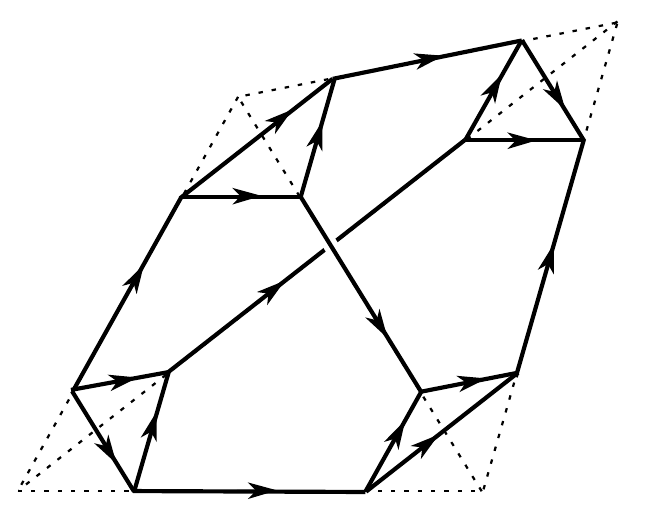}}

\definecolor{inkcol1}{rgb}{0.0,0.0,0.0}
   \put(132.125776,98.110286){\rotatebox{360.0}{\makebox(0,0)[tl]{\strut{}{
    \begin{minipage}[h]{127.655896pt}
\textcolor{inkcol1}{\LARGE{$\alpha_{12}$}}\\
\end{minipage}}}}}%

\definecolor{inkcol1}{rgb}{0.0,0.0,0.0}
   \put(-9.527818,15.444966){\rotatebox{360.0}{\makebox(0,0)[tl]{\strut{}{
    \begin{minipage}[h]{127.655896pt}
\textcolor{inkcol1}{\LARGE{$0$}}\\
\end{minipage}}}}}%

\definecolor{inkcol1}{rgb}{0.0,0.0,0.0}
   \put(70.907446,177.710316){\rotatebox{360.0}{\makebox(0,0)[tl]{\strut{}{
    \begin{minipage}[h]{127.655896pt}
\textcolor{inkcol1}{\LARGE{$1$}}\\
\end{minipage}}}}}%

\definecolor{inkcol1}{rgb}{0.0,0.0,0.0}
   \put(175.984726,20.094396){\rotatebox{360.0}{\makebox(0,0)[tl]{\strut{}{
    \begin{minipage}[h]{127.655896pt}
\textcolor{inkcol1}{\LARGE{$2$}}\\
\end{minipage}}}}}%

\definecolor{inkcol1}{rgb}{0.0,0.0,0.0}
   \put(222.479106,202.817286){\rotatebox{360.0}{\makebox(0,0)[tl]{\strut{}{
    \begin{minipage}[h]{127.655896pt}
\textcolor{inkcol1}{\LARGE{$3$}}\\
\end{minipage}}}}}%

\definecolor{inkcol1}{rgb}{0.0,0.0,0.0}
   \put(116.839056,151.952966){\rotatebox{360.0}{\makebox(0,0)[tl]{\strut{}{
    \begin{minipage}[h]{127.655896pt}
\textcolor{inkcol1}{\LARGE{$\beta^1_{23}$}}\\
\end{minipage}}}}}%

\definecolor{inkcol1}{rgb}{0.0,0.0,0.0}
   \put(201.121796,104.969076){\rotatebox{360.0}{\makebox(0,0)[tl]{\strut{}{
    \begin{minipage}[h]{127.655896pt}
\textcolor{inkcol1}{\LARGE{$\alpha_{23}$}}\\
\end{minipage}}}}}%

\definecolor{inkcol1}{rgb}{0.0,0.0,0.0}
   \put(134.416366,191.854876){\rotatebox{360.0}{\makebox(0,0)[tl]{\strut{}{
    \begin{minipage}[h]{127.655896pt}
\textcolor{inkcol1}{\LARGE{$\alpha_{13}$}}\\
\end{minipage}}}}}%

\definecolor{inkcol1}{rgb}{0.0,0.0,0.0}
   \put(12.735455,103.157526){\rotatebox{360.0}{\makebox(0,0)[tl]{\strut{}{
    \begin{minipage}[h]{127.655896pt}
\textcolor{inkcol1}{\LARGE{$\alpha_{01}$}}\\
\end{minipage}}}}}%

\definecolor{inkcol1}{rgb}{0.0,0.0,0.0}
   \put(76.151786,10.381006){\rotatebox{360.0}{\makebox(0,0)[tl]{\strut{}{
    \begin{minipage}[h]{127.655896pt}
\textcolor{inkcol1}{\LARGE{$\alpha_{02}$}}\\
\end{minipage}}}}}%

\definecolor{inkcol1}{rgb}{0.0,0.0,0.0}
   \put(93.303036,83.426486){\rotatebox{360.0}{\makebox(0,0)[tl]{\strut{}{
    \begin{minipage}[h]{127.655896pt}
\textcolor{inkcol1}{\LARGE{$\alpha_{03}$}}\\
\end{minipage}}}}}%

\definecolor{inkcol1}{rgb}{0.0,0.0,0.0}
   \put(57.538866,51.362936){\rotatebox{360.0}{\makebox(0,0)[tl]{\strut{}{
    \begin{minipage}[h]{127.655896pt}
\textcolor{inkcol1}{\LARGE{$\beta^0_{23}$}}\\
\end{minipage}}}}}%

\definecolor{inkcol1}{rgb}{0.0,0.0,0.0}
   \put(168.015746,136.599256){\rotatebox{360.0}{\makebox(0,0)[tl]{\strut{}{
    \begin{minipage}[h]{127.655896pt}
\textcolor{inkcol1}{\LARGE{$\beta^3_{02}$}}\\
\end{minipage}}}}}%

 \end{picture}
\endgroup

%% file: figures_gen/NaturalCocycle.tex
\begingroup
 \setlength{\unitlength}{0.8pt}
 \begin{picture}(250.69029,200.0932)
 \put(0,0){\includegraphics{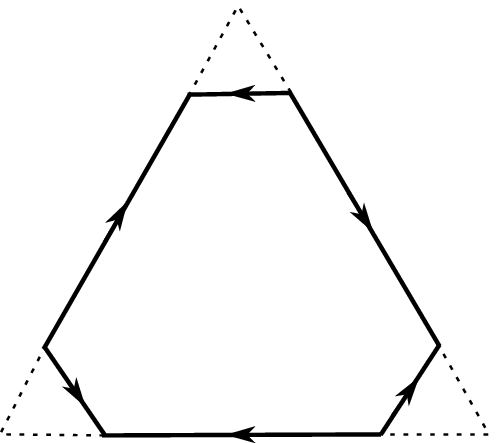}}

\definecolor{inkcol1}{rgb}{0.0,0.0,0.0}
   \put(-0.232683,20.774626){\rotatebox{360.0}{\makebox(0,0)[tl]{\strut{}{
    \begin{minipage}[h]{127.655896pt}
\textcolor{inkcol1}{\LARGE{$0$}}\\
\end{minipage}}}}}%

\definecolor{inkcol1}{rgb}{0.0,0.0,0.0}
   \put(196.195897,28.274626){\rotatebox{360.0}{\makebox(0,0)[tl]{\strut{}{
    \begin{minipage}[h]{127.655896pt}
\textcolor{inkcol1}{\LARGE{$2$}}\\
\end{minipage}}}}}%

\definecolor{inkcol1}{rgb}{0.0,0.0,0.0}
   \put(91.910167,200.060356){\rotatebox{360.0}{\makebox(0,0)[tl]{\strut{}{
    \begin{minipage}[h]{127.655896pt}
\textcolor{inkcol1}{\LARGE{$1$}}\\
\end{minipage}}}}}%

\definecolor{inkcol1}{rgb}{0.0,0.0,0.0}
   \put(49.593997,60.263396){\rotatebox{360.0}{\makebox(0,0)[tl]{\strut{}{
    \begin{minipage}[h]{127.655896pt}
\textcolor{inkcol1}{\LARGE{$\beta^0_{12}$}}\\
\end{minipage}}}}}%

\definecolor{inkcol1}{rgb}{0.0,0.0,0.0}
   \put(93.961797,144.184796){\rotatebox{360.0}{\makebox(0,0)[tl]{\strut{}{
    \begin{minipage}[h]{127.655896pt}
\textcolor{inkcol1}{\LARGE{$\beta^1_{20}$}}\\
\end{minipage}}}}}%

\definecolor{inkcol1}{rgb}{0.0,0.0,0.0}
   \put(130.390357,61.684796){\rotatebox{360.0}{\makebox(0,0)[tl]{\strut{}{
    \begin{minipage}[h]{127.655896pt}
\textcolor{inkcol1}{\LARGE{$\beta^2_{01}$}}\\
\end{minipage}}}}}%

\definecolor{inkcol1}{rgb}{0.0,0.0,0.0}
   \put(30.390357,121.684796){\rotatebox{360.0}{\makebox(0,0)[tl]{\strut{}{
    \begin{minipage}[h]{127.655896pt}
\textcolor{inkcol1}{\LARGE{$\alpha_{01}$}}\\
\end{minipage}}}}}%

\definecolor{inkcol1}{rgb}{0.0,0.0,0.0}
   \put(150.247497,119.040816){\rotatebox{360.0}{\makebox(0,0)[tl]{\strut{}{
    \begin{minipage}[h]{127.655896pt}
\textcolor{inkcol1}{\LARGE{$\alpha_{12}$}}\\
\end{minipage}}}}}%

\definecolor{inkcol1}{rgb}{0.0,0.0,0.0}
   \put(96.676077,20.256216){\rotatebox{360.0}{\makebox(0,0)[tl]{\strut{}{
    \begin{minipage}[h]{127.655896pt}
\textcolor{inkcol1}{\LARGE{$\alpha_{20}$}}\\
\end{minipage}}}}}%

 \end{picture}
\endgroup

%% file: figures_gen/QuiverG2Minors.tex
\begingroup
 \setlength{\unitlength}{0.8pt}
 \begin{picture}(690.20831,164.82143)
 \put(0,0){\includegraphics{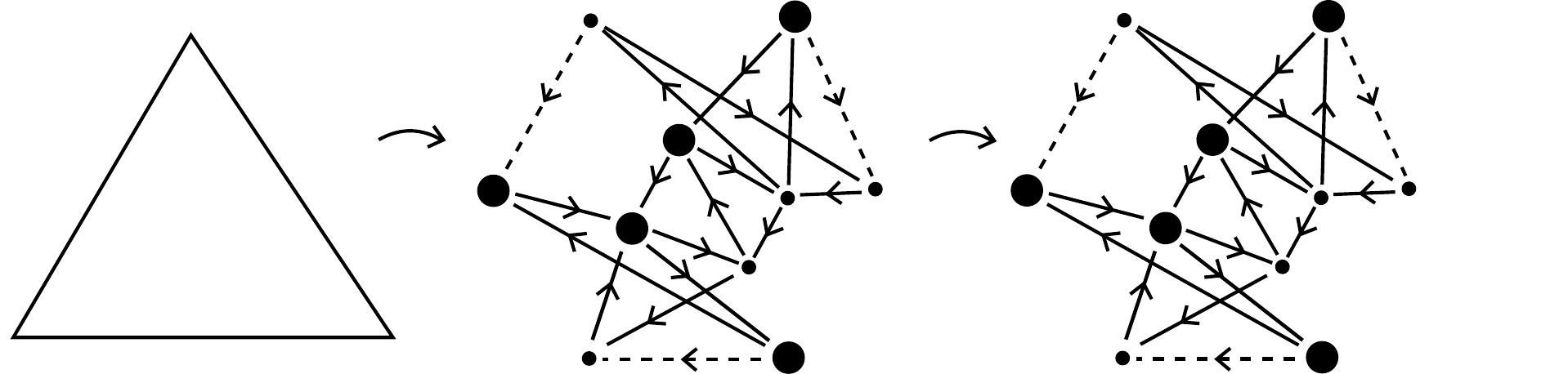}}

\definecolor{inkcol1}{rgb}{0.0,0.0,0.0}
   \put(259.207881,91.857796){\rotatebox{360.0}{\makebox(0,0)[tl]{\strut{}{
    \begin{minipage}[h]{127.655896pt}
\textcolor{inkcol1}{\LARGE{$u_3$}}\\
\end{minipage}}}}}%

\definecolor{inkcol1}{rgb}{0.0,0.0,0.0}
   \put(333.136451,48.286366){\rotatebox{360.0}{\makebox(0,0)[tl]{\strut{}{
    \begin{minipage}[h]{127.655896pt}
\textcolor{inkcol1}{\LARGE{$u_4$}}\\
\end{minipage}}}}}%

\definecolor{inkcol1}{rgb}{0.0,0.0,0.0}
   \put(347.779311,75.429216){\rotatebox{360.0}{\makebox(0,0)[tl]{\strut{}{
    \begin{minipage}[h]{127.655896pt}
\textcolor{inkcol1}{\LARGE{$u_6$}}\\
\end{minipage}}}}}%

\definecolor{inkcol1}{rgb}{0.0,0.0,0.0}
   \put(269.565011,122.214936){\rotatebox{360.0}{\makebox(0,0)[tl]{\strut{}{
    \begin{minipage}[h]{127.655896pt}
\textcolor{inkcol1}{\LARGE{$u_5$}}\\
\end{minipage}}}}}%

\definecolor{inkcol1}{rgb}{0.0,0.0,0.0}
   \put(494.142451,92.036366){\rotatebox{360.0}{\makebox(0,0)[tl]{\strut{}{
    \begin{minipage}[h]{127.655896pt}
\textcolor{inkcol1}{\LARGE{$a_1$}}\\
\end{minipage}}}}}%

\definecolor{inkcol1}{rgb}{0.0,0.0,0.0}
   \put(568.071021,48.464936){\rotatebox{360.0}{\makebox(0,0)[tl]{\strut{}{
    \begin{minipage}[h]{127.655896pt}
\textcolor{inkcol1}{\LARGE{$a_2$}}\\
\end{minipage}}}}}%

\definecolor{inkcol1}{rgb}{0.0,0.0,0.0}
   \put(582.713881,75.607786){\rotatebox{360.0}{\makebox(0,0)[tl]{\strut{}{
    \begin{minipage}[h]{127.655896pt}
\textcolor{inkcol1}{\LARGE{$a_4$}}\\
\end{minipage}}}}}%

\definecolor{inkcol1}{rgb}{0.0,0.0,0.0}
   \put(504.856721,124.179216){\rotatebox{360.0}{\makebox(0,0)[tl]{\strut{}{
    \begin{minipage}[h]{127.655896pt}
\textcolor{inkcol1}{\LARGE{$a_3$}}\\
\end{minipage}}}}}%

\definecolor{inkcol1}{rgb}{0.0,0.0,0.0}
   \put(-18.14926,15.096066){\rotatebox{360.0}{\makebox(0,0)[tl]{\strut{}{
    \begin{minipage}[h]{127.655896pt}
\textcolor{inkcol1}{\LARGE{$g_0N$}}\\
\end{minipage}}}}}%

\definecolor{inkcol1}{rgb}{0.0,0.0,0.0}
   \put(66.493596,170.096066){\rotatebox{360.0}{\makebox(0,0)[tl]{\strut{}{
    \begin{minipage}[h]{127.655896pt}
\textcolor{inkcol1}{\LARGE{$g_1N$}}\\
\end{minipage}}}}}%

\definecolor{inkcol1}{rgb}{0.0,0.0,0.0}
   \put(171.850741,19.738926){\rotatebox{360.0}{\makebox(0,0)[tl]{\strut{}{
    \begin{minipage}[h]{127.655896pt}
\textcolor{inkcol1}{\LARGE{$g_2N$}}\\
\end{minipage}}}}}%

\definecolor{inkcol1}{rgb}{0.0,0.0,0.0}
   \put(166.636451,134.024636){\rotatebox{360.0}{\makebox(0,0)[tl]{\strut{}{
    \begin{minipage}[h]{127.655896pt}
\textcolor{inkcol1}{\LARGE{$\Delta$}}\\
\end{minipage}}}}}%

\definecolor{inkcol1}{rgb}{0.0,0.0,0.0}
   \put(189.350731,70.453216){\rotatebox{360.0}{\makebox(0,0)[tl]{\strut{}{
    \begin{minipage}[h]{127.655896pt}
\textcolor{inkcol1}{\LARGE{$a_{01}$}}\\
\end{minipage}}}}}%

\definecolor{inkcol1}{rgb}{0.0,0.0,0.0}
   \put(232.922171,176.167496){\rotatebox{360.0}{\makebox(0,0)[tl]{\strut{}{
    \begin{minipage}[h]{127.655896pt}
\textcolor{inkcol1}{\LARGE{$a_{12}$}}\\
\end{minipage}}}}}%

\definecolor{inkcol1}{rgb}{0.0,0.0,0.0}
   \put(361.493591,170.096066){\rotatebox{360.0}{\makebox(0,0)[tl]{\strut{}{
    \begin{minipage}[h]{127.655896pt}
\textcolor{inkcol1}{\LARGE{$a_{12}$}}\\
\end{minipage}}}}}%

\definecolor{inkcol1}{rgb}{0.0,0.0,0.0}
   \put(391.850741,96.881786){\rotatebox{360.0}{\makebox(0,0)[tl]{\strut{}{
    \begin{minipage}[h]{127.655896pt}
\textcolor{inkcol1}{\LARGE{$a_{21}$}}\\
\end{minipage}}}}}%

\definecolor{inkcol1}{rgb}{0.0,0.0,0.0}
   \put(354.350741,6.881786){\rotatebox{360.0}{\makebox(0,0)[tl]{\strut{}{
    \begin{minipage}[h]{127.655896pt}
\textcolor{inkcol1}{\LARGE{$a_{20}$}}\\
\end{minipage}}}}}%

\definecolor{inkcol1}{rgb}{0.0,0.0,0.0}
   \put(235.850741,4.714926){\rotatebox{360.0}{\makebox(0,0)[tl]{\strut{}{
    \begin{minipage}[h]{127.655896pt}
\textcolor{inkcol1}{\LARGE{$a_{02}$}}\\
\end{minipage}}}}}%

\definecolor{inkcol1}{rgb}{0.0,0.0,0.0}
   \put(584.422161,2.572066){\rotatebox{360.0}{\makebox(0,0)[tl]{\strut{}{
    \begin{minipage}[h]{127.655896pt}
\textcolor{inkcol1}{\LARGE{$a_{20}$}}\\
\end{minipage}}}}}%

\definecolor{inkcol1}{rgb}{0.0,0.0,0.0}
   \put(472.993601,3.286366){\rotatebox{360.0}{\makebox(0,0)[tl]{\strut{}{
    \begin{minipage}[h]{127.655896pt}
\textcolor{inkcol1}{\LARGE{$a_{02}$}}\\
\end{minipage}}}}}%

\definecolor{inkcol1}{rgb}{0.0,0.0,0.0}
   \put(429.422171,75.072066){\rotatebox{360.0}{\makebox(0,0)[tl]{\strut{}{
    \begin{minipage}[h]{127.655896pt}
\textcolor{inkcol1}{\LARGE{$a_{01}$}}\\
\end{minipage}}}}}%

\definecolor{inkcol1}{rgb}{0.0,0.0,0.0}
   \put(470.136461,177.929216){\rotatebox{360.0}{\makebox(0,0)[tl]{\strut{}{
    \begin{minipage}[h]{127.655896pt}
\textcolor{inkcol1}{\LARGE{$a_{10}$}}\\
\end{minipage}}}}}%

\definecolor{inkcol1}{rgb}{0.0,0.0,0.0}
   \put(596.565021,171.857786){\rotatebox{360.0}{\makebox(0,0)[tl]{\strut{}{
    \begin{minipage}[h]{127.655896pt}
\textcolor{inkcol1}{\LARGE{$a_{12}$}}\\
\end{minipage}}}}}%

\definecolor{inkcol1}{rgb}{0.0,0.0,0.0}
   \put(627.993601,95.786366){\rotatebox{360.0}{\makebox(0,0)[tl]{\strut{}{
    \begin{minipage}[h]{127.655896pt}
\textcolor{inkcol1}{\LARGE{$a_{21}$}}\\
\end{minipage}}}}}%

\definecolor{inkcol1}{rgb}{0.0,0.0,0.0}
   \put(638.636461,165.453216){\rotatebox{360.0}{\makebox(0,0)[tl]{\strut{}{
    \begin{minipage}[h]{127.655896pt}
\textcolor{inkcol1}{\LARGE{$a_1=u_1a_{12}$}}\\
\end{minipage}}}}}%

\definecolor{inkcol1}{rgb}{0.0,0.0,0.0}
   \put(638.747501,145.994506){\rotatebox{360.0}{\makebox(0,0)[tl]{\strut{}{
    \begin{minipage}[h]{127.655896pt}
\textcolor{inkcol1}{\LARGE{$a_2=u_2\frac{a_{21}a_{01}^3}{a_{10}}$}}\\
\end{minipage}}}}}%

\definecolor{inkcol1}{rgb}{0.0,0.0,0.0}
   \put(639.966101,65.453216){\rotatebox{360.0}{\makebox(0,0)[tl]{\strut{}{
    \begin{minipage}[h]{127.655896pt}
\textcolor{inkcol1}{\LARGE{$a_3=u_3a_{01}a_{12}$}}\\
\end{minipage}}}}}%

\definecolor{inkcol1}{rgb}{0.0,0.0,0.0}
   \put(639.608981,42.596056){\rotatebox{360.0}{\makebox(0,0)[tl]{\strut{}{
    \begin{minipage}[h]{127.655896pt}
\textcolor{inkcol1}{\LARGE{$a_4=u_4\frac{a_{21}a_{01}^3}{a_{10}}$}}\\
\end{minipage}}}}}%

\definecolor{inkcol1}{rgb}{0.0,0.0,0.0}
   \put(400.850741,131.149076){\rotatebox{360.0}{\makebox(0,0)[tl]{\strut{}{
    \begin{minipage}[h]{127.655896pt}
\textcolor{inkcol1}{\LARGE{$m_{G_2}$}}\\
\end{minipage}}}}}%

 \end{picture}
\endgroup

%% file: figures_gen/C2Coordinates.tex
\begingroup
 \setlength{\unitlength}{0.8pt}
 \begin{picture}(596.56866,157.89284)
 \put(0,0){\includegraphics{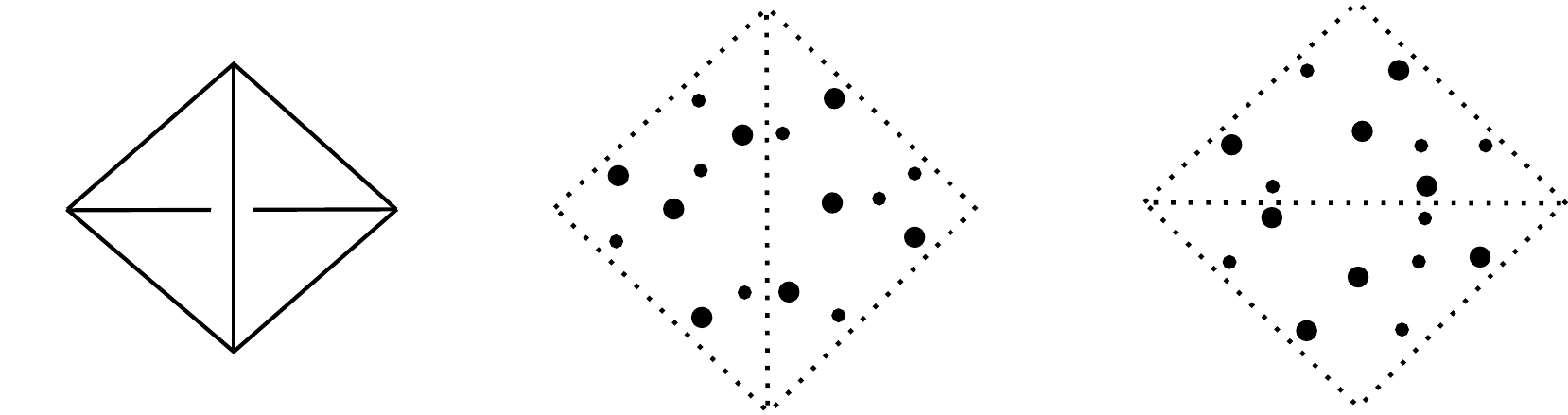}}

\definecolor{inkcol1}{rgb}{0.0,0.0,0.0}
   \put(73.366083,25.355636){\rotatebox{360.0}{\makebox(0,0)[tl]{\strut{}{
    \begin{minipage}[h]{127.655896pt}
\textcolor{inkcol1}{\Large{$g_0N$}}\\
\end{minipage}}}}}%

\definecolor{inkcol1}{rgb}{0.0,0.0,0.0}
   \put(-0.205343,99.641346){\rotatebox{360.0}{\makebox(0,0)[tl]{\strut{}{
    \begin{minipage}[h]{127.655896pt}
\textcolor{inkcol1}{\Large{$g_1N$}}\\
\end{minipage}}}}}%

\definecolor{inkcol1}{rgb}{0.0,0.0,0.0}
   \put(72.294663,154.998486){\rotatebox{360.0}{\makebox(0,0)[tl]{\strut{}{
    \begin{minipage}[h]{127.655896pt}
\textcolor{inkcol1}{\Large{$g_2N$}}\\
\end{minipage}}}}}%

\definecolor{inkcol1}{rgb}{0.0,0.0,0.0}
   \put(143.366097,101.427076){\rotatebox{360.0}{\makebox(0,0)[tl]{\strut{}{
    \begin{minipage}[h]{127.655896pt}
\textcolor{inkcol1}{\Large{$g_3N$}}\\
\end{minipage}}}}}%

\definecolor{inkcol1}{rgb}{0.0,0.0,0.0}
   \put(254.908497,56.431526){\rotatebox{360.0}{\makebox(0,0)[tl]{\strut{}{
    \begin{minipage}[h]{127.655896pt}
\textcolor{inkcol1}{\Large{$a_{20}$}}\\
\end{minipage}}}}}%

\definecolor{inkcol1}{rgb}{0.0,0.0,0.0}
   \put(255.265637,79.288676){\rotatebox{360.0}{\makebox(0,0)[tl]{\strut{}{
    \begin{minipage}[h]{127.655896pt}
\textcolor{inkcol1}{\Large{$f_{1,3}$}}\\
\end{minipage}}}}}%

\definecolor{inkcol1}{rgb}{0.0,0.0,0.0}
   \put(263.837067,96.431526){\rotatebox{360.0}{\makebox(0,0)[tl]{\strut{}{
    \begin{minipage}[h]{127.655896pt}
\textcolor{inkcol1}{\Large{$f_{2,3}$}}\\
\end{minipage}}}}}%

\definecolor{inkcol1}{rgb}{0.0,0.0,0.0}
   \put(239.194197,113.217246){\rotatebox{360.0}{\makebox(0,0)[tl]{\strut{}{
    \begin{minipage}[h]{127.655896pt}
\textcolor{inkcol1}{\Large{$-a_{02}$}}\\
\end{minipage}}}}}%

\definecolor{inkcol1}{rgb}{0.0,0.0,0.0}
   \put(297.979927,116.846686){\rotatebox{360.0}{\makebox(0,0)[tl]{\strut{}{
    \begin{minipage}[h]{127.655896pt}
\textcolor{inkcol1}{\Large{$a_{20}$}}\\
\end{minipage}}}}}%

\definecolor{inkcol1}{rgb}{0.0,0.0,0.0}
   \put(302.622787,55.418136){\rotatebox{360.0}{\makebox(0,0)[tl]{\strut{}{
    \begin{minipage}[h]{127.655896pt}
\textcolor{inkcol1}{\Large{$a_{02}$}}\\
\end{minipage}}}}}%

\definecolor{inkcol1}{rgb}{0.0,0.0,0.0}
   \put(295.051357,79.346706){\rotatebox{360.0}{\makebox(0,0)[tl]{\strut{}{
    \begin{minipage}[h]{127.655896pt}
\textcolor{inkcol1}{\Large{$f_{1,1}$}}\\
\end{minipage}}}}}%

\definecolor{inkcol1}{rgb}{0.0,0.0,0.0}
   \put(317.194217,102.203856){\rotatebox{360.0}{\makebox(0,0)[tl]{\strut{}{
    \begin{minipage}[h]{127.655896pt}
\textcolor{inkcol1}{\Large{$f_{2,1}$}}\\
\end{minipage}}}}}%

\definecolor{inkcol1}{rgb}{0.0,0.0,0.0}
   \put(325.837067,36.811006){\rotatebox{360.0}{\makebox(0,0)[tl]{\strut{}{
    \begin{minipage}[h]{127.655896pt}
\textcolor{inkcol1}{\Large{$a_{03}$}}\\
\end{minipage}}}}}%

\definecolor{inkcol1}{rgb}{0.0,0.0,0.0}
   \put(353.694197,64.311016){\rotatebox{360.0}{\makebox(0,0)[tl]{\strut{}{
    \begin{minipage}[h]{127.655896pt}
\textcolor{inkcol1}{\Large{$a_{30}$}}\\
\end{minipage}}}}}%

\definecolor{inkcol1}{rgb}{0.0,0.0,0.0}
   \put(236.551357,34.668136){\rotatebox{360.0}{\makebox(0,0)[tl]{\strut{}{
    \begin{minipage}[h]{127.655896pt}
\textcolor{inkcol1}{\Large{$a_{01}$}}\\
\end{minipage}}}}}%

\definecolor{inkcol1}{rgb}{0.0,0.0,0.0}
   \put(207.979917,63.239576){\rotatebox{360.0}{\makebox(0,0)[tl]{\strut{}{
    \begin{minipage}[h]{127.655896pt}
\textcolor{inkcol1}{\Large{$a_{10}$}}\\
\end{minipage}}}}}%

\definecolor{inkcol1}{rgb}{0.0,0.0,0.0}
   \put(208.337067,112.525286){\rotatebox{360.0}{\makebox(0,0)[tl]{\strut{}{
    \begin{minipage}[h]{127.655896pt}
\textcolor{inkcol1}{\Large{$a_{12}$}}\\
\end{minipage}}}}}%

\definecolor{inkcol1}{rgb}{0.0,0.0,0.0}
   \put(238.694217,137.882436){\rotatebox{360.0}{\makebox(0,0)[tl]{\strut{}{
    \begin{minipage}[h]{127.655896pt}
\textcolor{inkcol1}{\Large{$a_{21}$}}\\
\end{minipage}}}}}%

\definecolor{inkcol1}{rgb}{0.0,0.0,0.0}
   \put(325.479927,135.025296){\rotatebox{360.0}{\makebox(0,0)[tl]{\strut{}{
    \begin{minipage}[h]{127.655896pt}
\textcolor{inkcol1}{\Large{$a_{23}$}}\\
\end{minipage}}}}}%

\definecolor{inkcol1}{rgb}{0.0,0.0,0.0}
   \put(353.694217,107.882426){\rotatebox{360.0}{\makebox(0,0)[tl]{\strut{}{
    \begin{minipage}[h]{127.655896pt}
\textcolor{inkcol1}{\Large{$a_{23}$}}\\
\end{minipage}}}}}%

\definecolor{inkcol1}{rgb}{0.0,0.0,0.0}
   \put(485.051347,118.596716){\rotatebox{360.0}{\makebox(0,0)[tl]{\strut{}{
    \begin{minipage}[h]{127.655896pt}
\textcolor{inkcol1}{\Large{$f_{1,0}$}}\\
\end{minipage}}}}}%

\definecolor{inkcol1}{rgb}{0.0,0.0,0.0}
   \put(530.051357,122.168146){\rotatebox{360.0}{\makebox(0,0)[tl]{\strut{}{
    \begin{minipage}[h]{127.655896pt}
\textcolor{inkcol1}{\Large{$f_{2,0}$}}\\
\end{minipage}}}}}%

\definecolor{inkcol1}{rgb}{0.0,0.0,0.0}
   \put(486.837067,93.953856){\rotatebox{360.0}{\makebox(0,0)[tl]{\strut{}{
    \begin{minipage}[h]{127.655896pt}
\textcolor{inkcol1}{\Large{$a_{13}$}}\\
\end{minipage}}}}}%

\definecolor{inkcol1}{rgb}{0.0,0.0,0.0}
   \put(516.479927,93.953856){\rotatebox{360.0}{\makebox(0,0)[tl]{\strut{}{
    \begin{minipage}[h]{127.655896pt}
\textcolor{inkcol1}{\Large{$a_{31}$}}\\
\end{minipage}}}}}%

\definecolor{inkcol1}{rgb}{0.0,0.0,0.0}
   \put(473.265637,74.668146){\rotatebox{360.0}{\makebox(0,0)[tl]{\strut{}{
    \begin{minipage}[h]{127.655896pt}
\textcolor{inkcol1}{\Large{$-a_{31}$}}\\
\end{minipage}}}}}%

\definecolor{inkcol1}{rgb}{0.0,0.0,0.0}
   \put(516.122777,75.739566){\rotatebox{360.0}{\makebox(0,0)[tl]{\strut{}{
    \begin{minipage}[h]{127.655896pt}
\textcolor{inkcol1}{\Large{$a_{13}$}}\\
\end{minipage}}}}}%

\definecolor{inkcol1}{rgb}{0.0,0.0,0.0}
   \put(496.479917,53.953866){\rotatebox{360.0}{\makebox(0,0)[tl]{\strut{}{
    \begin{minipage}[h]{127.655896pt}
\textcolor{inkcol1}{\Large{$f_{1,2}$}}\\
\end{minipage}}}}}%

\definecolor{inkcol1}{rgb}{0.0,0.0,0.0}
   \put(524.337067,59.311006){\rotatebox{360.0}{\makebox(0,0)[tl]{\strut{}{
    \begin{minipage}[h]{127.655896pt}
\textcolor{inkcol1}{\Large{$f_{2,2}$}}\\
\end{minipage}}}}}%

\definecolor{inkcol1}{rgb}{0.0,0.0,0.0}
   \put(441.908497,122.882436){\rotatebox{360.0}{\makebox(0,0)[tl]{\strut{}{
    \begin{minipage}[h]{127.655896pt}
\textcolor{inkcol1}{\Large{$a_{12}$}}\\
\end{minipage}}}}}%

\definecolor{inkcol1}{rgb}{0.0,0.0,0.0}
   \put(473.694217,150.739586){\rotatebox{360.0}{\makebox(0,0)[tl]{\strut{}{
    \begin{minipage}[h]{127.655896pt}
\textcolor{inkcol1}{\Large{$a_{21}$}}\\
\end{minipage}}}}}%

\definecolor{inkcol1}{rgb}{0.0,0.0,0.0}
   \put(542.265637,145.382426){\rotatebox{360.0}{\makebox(0,0)[tl]{\strut{}{
    \begin{minipage}[h]{127.655896pt}
\textcolor{inkcol1}{\Large{$a_{23}$}}\\
\end{minipage}}}}}%

\definecolor{inkcol1}{rgb}{0.0,0.0,0.0}
   \put(571.551357,117.882436){\rotatebox{360.0}{\makebox(0,0)[tl]{\strut{}{
    \begin{minipage}[h]{127.655896pt}
\textcolor{inkcol1}{\Large{$a_{32}$}}\\
\end{minipage}}}}}%

\definecolor{inkcol1}{rgb}{0.0,0.0,0.0}
   \put(540.837067,28.953856){\rotatebox{360.0}{\makebox(0,0)[tl]{\strut{}{
    \begin{minipage}[h]{127.655896pt}
\textcolor{inkcol1}{\Large{$a_{03}$}}\\
\end{minipage}}}}}%

\definecolor{inkcol1}{rgb}{0.0,0.0,0.0}
   \put(571.551357,58.239566){\rotatebox{360.0}{\makebox(0,0)[tl]{\strut{}{
    \begin{minipage}[h]{127.655896pt}
\textcolor{inkcol1}{\Large{$a_{30}$}}\\
\end{minipage}}}}}%

\definecolor{inkcol1}{rgb}{0.0,0.0,0.0}
   \put(466.194207,30.025296){\rotatebox{360.0}{\makebox(0,0)[tl]{\strut{}{
    \begin{minipage}[h]{127.655896pt}
\textcolor{inkcol1}{\Large{$a_{01}$}}\\
\end{minipage}}}}}%

\definecolor{inkcol1}{rgb}{0.0,0.0,0.0}
   \put(439.765637,54.668156){\rotatebox{360.0}{\makebox(0,0)[tl]{\strut{}{
    \begin{minipage}[h]{127.655896pt}
\textcolor{inkcol1}{\Large{$a_{10}$}}\\
\end{minipage}}}}}%

 \end{picture}
\endgroup

%% file: figures_gen/ObstructionIdentification.tex
\begingroup
 \setlength{\unitlength}{0.8pt}
 \begin{picture}(400.60101,127.05021)
 \put(0,0){\includegraphics{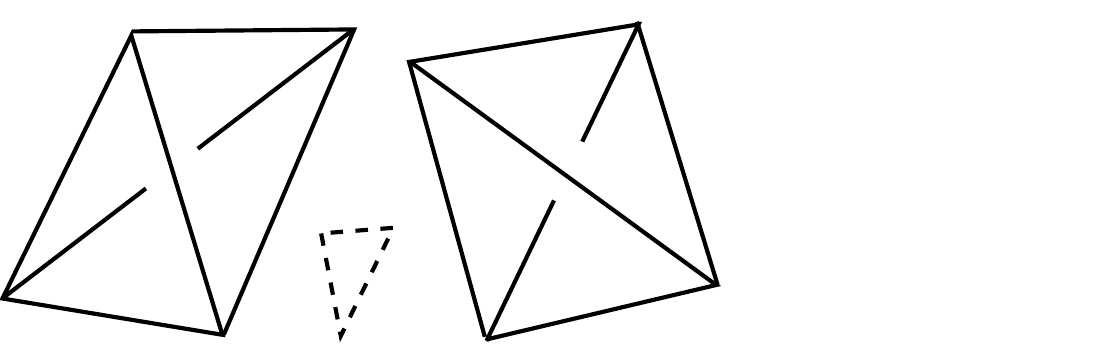}}

\definecolor{inkcol1}{rgb}{0.0,0.0,0.0}
   \put(1.814647,108.734656){\rotatebox{360.0}{\makebox(0,0)[tl]{\strut{}{
    \begin{minipage}[h]{127.655896pt}
\textcolor{inkcol1}{\LARGE{$\Delta_s$}}\\
\end{minipage}}}}}%

\definecolor{inkcol1}{rgb}{0.0,0.0,0.0}
   \put(248.864217,95.954236){\rotatebox{360.0}{\makebox(0,0)[tl]{\strut{}{
    \begin{minipage}[h]{127.655896pt}
\textcolor{inkcol1}{\LARGE{$\Delta_{s'}$}}\\
\end{minipage}}}}}%

\definecolor{inkcol1}{rgb}{0.0,0.0,0.0}
   \put(130.423837,20.950416){\rotatebox{360.0}{\makebox(0,0)[tl]{\strut{}{
    \begin{minipage}[h]{127.655896pt}
\textcolor{inkcol1}{\LARGE{$\Delta^2$}}\\
\end{minipage}}}}}%

\definecolor{inkcol1}{rgb}{0.0,0.0,0.0}
   \put(86.776087,82.734346){\rotatebox{360.0}{\makebox(0,0)[tl]{\strut{}{
    \begin{minipage}[h]{127.655896pt}
\textcolor{inkcol1}{\LARGE{$f$}}\\
\end{minipage}}}}}%

\definecolor{inkcol1}{rgb}{0.0,0.0,0.0}
   \put(186.781187,96.118866){\rotatebox{360.0}{\makebox(0,0)[tl]{\strut{}{
    \begin{minipage}[h]{127.655896pt}
\textcolor{inkcol1}{\LARGE{$f'$}}\\
\end{minipage}}}}}%

\definecolor{inkcol1}{rgb}{0.0,0.0,0.0}
   \put(42.009577,56.371356){\rotatebox{360.0}{\makebox(0,0)[tl]{\strut{}{
    \begin{minipage}[h]{127.655896pt}
\textcolor{inkcol1}{\LARGE{$-$}}\\
\end{minipage}}}}}%

\definecolor{inkcol1}{rgb}{0.0,0.0,0.0}
   \put(75.849677,134.910716){\rotatebox{360.0}{\makebox(0,0)[tl]{\strut{}{
    \begin{minipage}[h]{127.655896pt}
\textcolor{inkcol1}{\LARGE{$-$}}\\
\end{minipage}}}}}%

\definecolor{inkcol1}{rgb}{0.0,0.0,0.0}
   \put(176.864947,65.715266){\rotatebox{360.0}{\makebox(0,0)[tl]{\strut{}{
    \begin{minipage}[h]{127.655896pt}
\textcolor{inkcol1}{\LARGE{$-$}}\\
\end{minipage}}}}}%

\definecolor{inkcol1}{rgb}{0.0,0.0,0.0}
   \put(135.196147,72.533796){\rotatebox{360.0}{\makebox(0,0)[tl]{\strut{}{
    \begin{minipage}[h]{127.655896pt}
\textcolor{inkcol1}{\LARGE{$-$}}\\
\end{minipage}}}}}%

\definecolor{inkcol1}{rgb}{0.0,0.0,0.0}
   \put(132.418227,34.400536){\rotatebox{360.0}{\makebox(0,0)[tl]{\strut{}{
    \begin{minipage}[h]{127.655896pt}
\textcolor{inkcol1}{\Large{$-$}}\\
\end{minipage}}}}}%

\definecolor{inkcol1}{rgb}{0.0,0.0,0.0}
   \put(120.672607,57.791926){\rotatebox{360.0}{\makebox(0,0)[tl]{\strut{}{
    \begin{minipage}[h]{127.655896pt}
\textcolor{inkcol1}{\Large{$-$}}\\
\end{minipage}}}}}%

 \end{picture}
\endgroup

%% file: figures_gen/Kappa.tex
\begingroup
 \setlength{\unitlength}{0.8pt}
 \begin{picture}(306.03638,107.80937)
 \put(0,0){\includegraphics{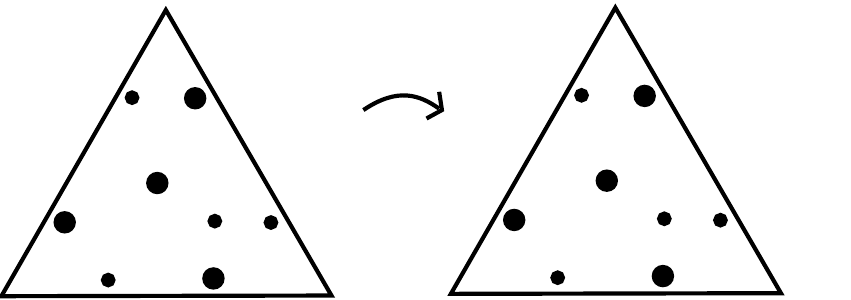}}

\definecolor{inkcol1}{rgb}{0.0,0.0,0.0}
   \put(136.458829,95.767296){\rotatebox{360.0}{\makebox(0,0)[tl]{\strut{}{
    \begin{minipage}[h]{127.655896pt}
\textcolor{inkcol1}{\LARGE{$\kappa_2$}}\\
\end{minipage}}}}}%

\definecolor{inkcol1}{rgb}{0.0,0.0,0.0}
   \put(-15.569123,-0.197194){\rotatebox{360.0}{\makebox(0,0)[tl]{\strut{}{
    \begin{minipage}[h]{127.655896pt}
\textcolor{inkcol1}{\LARGE{$g_0N$}}\\
\end{minipage}}}}}%

\definecolor{inkcol1}{rgb}{0.0,0.0,0.0}
   \put(54.131399,128.597256){\rotatebox{360.0}{\makebox(0,0)[tl]{\strut{}{
    \begin{minipage}[h]{127.655896pt}
\textcolor{inkcol1}{\LARGE{$g_1N$}}\\
\end{minipage}}}}}%

\definecolor{inkcol1}{rgb}{0.0,0.0,0.0}
   \put(104.639029,2.328186){\rotatebox{360.0}{\makebox(0,0)[tl]{\strut{}{
    \begin{minipage}[h]{127.655896pt}
\textcolor{inkcol1}{\LARGE{$g_2N$}}\\
\end{minipage}}}}}%

\definecolor{inkcol1}{rgb}{0.0,0.0,0.0}
   \put(142.772289,1.318036){\rotatebox{360.0}{\makebox(0,0)[tl]{\strut{}{
    \begin{minipage}[h]{127.655896pt}
\textcolor{inkcol1}{\LARGE{$g_0N$}}\\
\end{minipage}}}}}%

\definecolor{inkcol1}{rgb}{0.0,0.0,0.0}
   \put(211.715199,129.607406){\rotatebox{360.0}{\makebox(0,0)[tl]{\strut{}{
    \begin{minipage}[h]{127.655896pt}
\textcolor{inkcol1}{\LARGE{$g_1N$}}\\
\end{minipage}}}}}%

\definecolor{inkcol1}{rgb}{0.0,0.0,0.0}
   \put(282.425879,8.894176){\rotatebox{360.0}{\makebox(0,0)[tl]{\strut{}{
    \begin{minipage}[h]{127.655896pt}
\textcolor{inkcol1}{\LARGE{$g_2s_{C_2}N$}}\\
\end{minipage}}}}}%

\definecolor{inkcol1}{rgb}{0.0,0.0,0.0}
   \put(76.354749,87.938616){\rotatebox{360.0}{\makebox(0,0)[tl]{\strut{}{
    \begin{minipage}[h]{127.655896pt}
\textcolor{inkcol1}{\LARGE{$a_{12}$}}\\
\end{minipage}}}}}%

\definecolor{inkcol1}{rgb}{0.0,0.0,0.0}
   \put(45.520989,37.782546){\rotatebox{360.0}{\makebox(0,0)[tl]{\strut{}{
    \begin{minipage}[h]{127.655896pt}
\textcolor{inkcol1}{\LARGE{$a_1$}}\\
\end{minipage}}}}}%

\definecolor{inkcol1}{rgb}{0.0,0.0,0.0}
   \put(69.512119,47.126466){\rotatebox{360.0}{\makebox(0,0)[tl]{\strut{}{
    \begin{minipage}[h]{127.655896pt}
\textcolor{inkcol1}{\LARGE{$a_2$}}\\
\end{minipage}}}}}%

\definecolor{inkcol1}{rgb}{0.0,0.0,0.0}
   \put(192.498189,38.035086){\rotatebox{360.0}{\makebox(0,0)[tl]{\strut{}{
    \begin{minipage}[h]{127.655896pt}
\textcolor{inkcol1}{\LARGE{$-a_1$}}\\
\end{minipage}}}}}%

\definecolor{inkcol1}{rgb}{0.0,0.0,0.0}
   \put(226.338309,48.894226){\rotatebox{360.0}{\makebox(0,0)[tl]{\strut{}{
    \begin{minipage}[h]{127.655896pt}
\textcolor{inkcol1}{\LARGE{$a_2$}}\\
\end{minipage}}}}}%

\definecolor{inkcol1}{rgb}{0.0,0.0,0.0}
   \put(25.065409,0.154376){\rotatebox{360.0}{\makebox(0,0)[tl]{\strut{}{
    \begin{minipage}[h]{127.655896pt}
\textcolor{inkcol1}{\LARGE{$a_{02}$}}\\
\end{minipage}}}}}%

\definecolor{inkcol1}{rgb}{0.0,0.0,0.0}
   \put(62.441059,0.154376){\rotatebox{360.0}{\makebox(0,0)[tl]{\strut{}{
    \begin{minipage}[h]{127.655896pt}
\textcolor{inkcol1}{\LARGE{$a_{20}$}}\\
\end{minipage}}}}}%

\definecolor{inkcol1}{rgb}{0.0,0.0,0.0}
   \put(-9.027236,47.631536){\rotatebox{360.0}{\makebox(0,0)[tl]{\strut{}{
    \begin{minipage}[h]{127.655896pt}
\textcolor{inkcol1}{\LARGE{$a_{01}$}}\\
\end{minipage}}}}}%

\definecolor{inkcol1}{rgb}{0.0,0.0,0.0}
   \put(12.943584,87.280016){\rotatebox{360.0}{\makebox(0,0)[tl]{\strut{}{
    \begin{minipage}[h]{127.655896pt}
\textcolor{inkcol1}{\LARGE{$a_{10}$}}\\
\end{minipage}}}}}%

\definecolor{inkcol1}{rgb}{0.0,0.0,0.0}
   \put(101.331929,45.358686){\rotatebox{360.0}{\makebox(0,0)[tl]{\strut{}{
    \begin{minipage}[h]{127.655896pt}
\textcolor{inkcol1}{\LARGE{$a_{21}$}}\\
\end{minipage}}}}}%

\definecolor{inkcol1}{rgb}{0.0,0.0,0.0}
   \put(187.447429,0.659446){\rotatebox{360.0}{\makebox(0,0)[tl]{\strut{}{
    \begin{minipage}[h]{127.655896pt}
\textcolor{inkcol1}{\LARGE{$a_{02}$}}\\
\end{minipage}}}}}%

\definecolor{inkcol1}{rgb}{0.0,0.0,0.0}
   \put(221.792619,1.922136){\rotatebox{360.0}{\makebox(0,0)[tl]{\strut{}{
    \begin{minipage}[h]{127.655896pt}
\textcolor{inkcol1}{\LARGE{$-a_{20}$}}\\
\end{minipage}}}}}%

\definecolor{inkcol1}{rgb}{0.0,0.0,0.0}
   \put(152.849709,47.884076){\rotatebox{360.0}{\makebox(0,0)[tl]{\strut{}{
    \begin{minipage}[h]{127.655896pt}
\textcolor{inkcol1}{\LARGE{$a_{01}$}}\\
\end{minipage}}}}}%

\definecolor{inkcol1}{rgb}{0.0,0.0,0.0}
   \put(176.588289,90.815556){\rotatebox{360.0}{\makebox(0,0)[tl]{\strut{}{
    \begin{minipage}[h]{127.655896pt}
\textcolor{inkcol1}{\LARGE{$a_{10}$}}\\
\end{minipage}}}}}%

\definecolor{inkcol1}{rgb}{0.0,0.0,0.0}
   \put(240.480449,86.522416){\rotatebox{360.0}{\makebox(0,0)[tl]{\strut{}{
    \begin{minipage}[h]{127.655896pt}
\textcolor{inkcol1}{\LARGE{$-a_{12}$}}\\
\end{minipage}}}}}%

\definecolor{inkcol1}{rgb}{0.0,0.0,0.0}
   \put(264.724099,44.601086){\rotatebox{360.0}{\makebox(0,0)[tl]{\strut{}{
    \begin{minipage}[h]{127.655896pt}
\textcolor{inkcol1}{\LARGE{$a_{21}$}}\\
\end{minipage}}}}}%

 \end{picture}
\endgroup

%% file: figures_gen/Fig8Trig.tex
\begingroup
 \setlength{\unitlength}{0.8pt}
 \begin{picture}(432.0108,148.31714)
 \put(0,0){\includegraphics{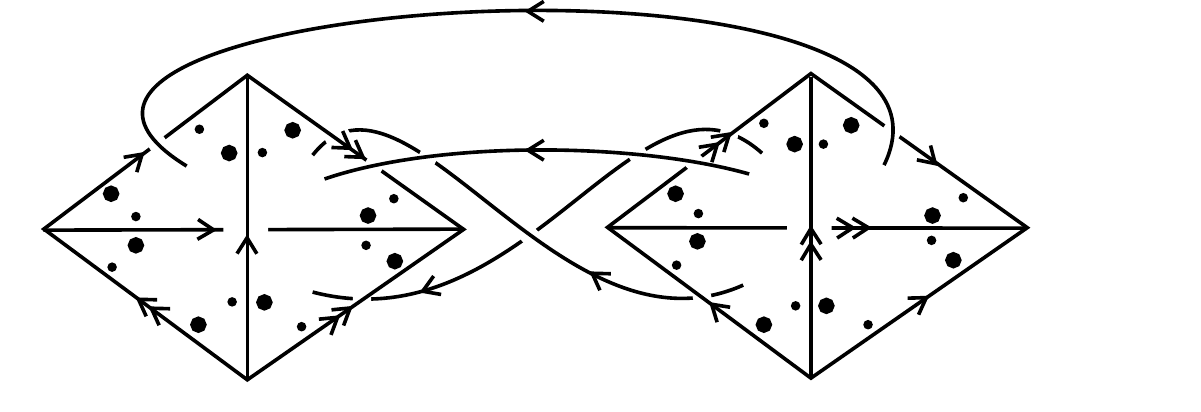}}

\definecolor{inkcol1}{rgb}{0.0,0.0,0.0}
   \put(74.683467,12.732756){\rotatebox{360.0}{\makebox(0,0)[tl]{\strut{}{
    \begin{minipage}[h]{127.655896pt}
\textcolor{inkcol1}{\large{$0$}}\\
\end{minipage}}}}}%

\definecolor{inkcol1}{rgb}{0.0,0.0,0.0}
   \put(0.733227,75.401716){\rotatebox{360.0}{\makebox(0,0)[tl]{\strut{}{
    \begin{minipage}[h]{127.655896pt}
\textcolor{inkcol1}{\large{$1$}}\\
\end{minipage}}}}}%

\definecolor{inkcol1}{rgb}{0.0,0.0,0.0}
   \put(91.037177,130.758856){\rotatebox{360.0}{\makebox(0,0)[tl]{\strut{}{
    \begin{minipage}[h]{127.655896pt}
\textcolor{inkcol1}{\large{$2$}}\\
\end{minipage}}}}}%

\definecolor{inkcol1}{rgb}{0.0,0.0,0.0}
   \put(165.571057,73.977986){\rotatebox{360.0}{\makebox(0,0)[tl]{\strut{}{
    \begin{minipage}[h]{127.655896pt}
\textcolor{inkcol1}{\large{$3$}}\\
\end{minipage}}}}}%

\definecolor{inkcol1}{rgb}{0.0,0.0,0.0}
   \put(276.331027,13.597006){\rotatebox{360.0}{\makebox(0,0)[tl]{\strut{}{
    \begin{minipage}[h]{127.655896pt}
\textcolor{inkcol1}{\large{$0$}}\\
\end{minipage}}}}}%

\definecolor{inkcol1}{rgb}{0.0,0.0,0.0}
   \put(205.616737,74.956856){\rotatebox{360.0}{\makebox(0,0)[tl]{\strut{}{
    \begin{minipage}[h]{127.655896pt}
\textcolor{inkcol1}{\large{$1$}}\\
\end{minipage}}}}}%

\definecolor{inkcol1}{rgb}{0.0,0.0,0.0}
   \put(275.764687,132.633256){\rotatebox{360.0}{\makebox(0,0)[tl]{\strut{}{
    \begin{minipage}[h]{127.655896pt}
\textcolor{inkcol1}{\large{$2$}}\\
\end{minipage}}}}}%

\definecolor{inkcol1}{rgb}{0.0,0.0,0.0}
   \put(372.352327,74.856456){\rotatebox{360.0}{\makebox(0,0)[tl]{\strut{}{
    \begin{minipage}[h]{127.655896pt}
\textcolor{inkcol1}{\large{$3$}}\\
\end{minipage}}}}}%

\definecolor{inkcol1}{rgb}{0.0,0.0,0.0}
   \put(176.351537,142.171576){\rotatebox{360.0}{\makebox(0,0)[tl]{\strut{}{
    \begin{minipage}[h]{127.655896pt}
\textcolor{inkcol1}{\large{$a$}}\\
\end{minipage}}}}}%

\definecolor{inkcol1}{rgb}{0.0,0.0,0.0}
   \put(194.094627,111.605626){\rotatebox{360.0}{\makebox(0,0)[tl]{\strut{}{
    \begin{minipage}[h]{127.655896pt}
\textcolor{inkcol1}{\large{$b$}}\\
\end{minipage}}}}}%

\definecolor{inkcol1}{rgb}{0.0,0.0,0.0}
   \put(215.863567,42.507446){\rotatebox{360.0}{\makebox(0,0)[tl]{\strut{}{
    \begin{minipage}[h]{127.655896pt}
\textcolor{inkcol1}{\large{$c$}}\\
\end{minipage}}}}}%

\definecolor{inkcol1}{rgb}{0.0,0.0,0.0}
   \put(145.182097,42.024036){\rotatebox{360.0}{\makebox(0,0)[tl]{\strut{}{
    \begin{minipage}[h]{127.655896pt}
\textcolor{inkcol1}{\large{$d$}}\\
\end{minipage}}}}}%

\definecolor{inkcol1}{rgb}{0.0,0.0,0.0}
   \put(92.652987,52.821256){\rotatebox{360.0}{\makebox(0,0)[tl]{\strut{}{
    \begin{minipage}[h]{127.655896pt}
\textcolor{inkcol1}{\large{$x$}}\\
\end{minipage}}}}}%

\definecolor{inkcol1}{rgb}{0.0,0.0,0.0}
   \put(59.149637,90.804666){\rotatebox{360.0}{\makebox(0,0)[tl]{\strut{}{
    \begin{minipage}[h]{127.655896pt}
\textcolor{inkcol1}{\large{$-x$}}\\
\end{minipage}}}}}%

\definecolor{inkcol1}{rgb}{0.0,0.0,0.0}
   \put(72.043317,53.791796){\rotatebox{360.0}{\makebox(0,0)[tl]{\strut{}{
    \begin{minipage}[h]{127.655896pt}
\textcolor{inkcol1}{\large{$y$}}\\
\end{minipage}}}}}%

\definecolor{inkcol1}{rgb}{0.0,0.0,0.0}
   \put(91.520287,91.266426){\rotatebox{360.0}{\makebox(0,0)[tl]{\strut{}{
    \begin{minipage}[h]{127.655896pt}
\textcolor{inkcol1}{\large{$y$}}\\
\end{minipage}}}}}%

\definecolor{inkcol1}{rgb}{0.0,0.0,0.0}
   \put(51.665257,29.444796){\rotatebox{360.0}{\makebox(0,0)[tl]{\strut{}{
    \begin{minipage}[h]{127.655896pt}
\textcolor{inkcol1}{\large{$z$}}\\
\end{minipage}}}}}%

\definecolor{inkcol1}{rgb}{0.0,0.0,0.0}
   \put(20.476797,50.657996){\rotatebox{360.0}{\makebox(0,0)[tl]{\strut{}{
    \begin{minipage}[h]{127.655896pt}
\textcolor{inkcol1}{\large{$w$}}\\
\end{minipage}}}}}%

\definecolor{inkcol1}{rgb}{0.0,0.0,0.0}
   \put(24.643667,93.589476){\rotatebox{360.0}{\makebox(0,0)[tl]{\strut{}{
    \begin{minipage}[h]{127.655896pt}
\textcolor{inkcol1}{\large{$x$}}\\
\end{minipage}}}}}%

\definecolor{inkcol1}{rgb}{0.0,0.0,0.0}
   \put(57.221097,117.075526){\rotatebox{360.0}{\makebox(0,0)[tl]{\strut{}{
    \begin{minipage}[h]{127.655896pt}
\textcolor{inkcol1}{\large{$y$}}\\
\end{minipage}}}}}%

\definecolor{inkcol1}{rgb}{0.0,0.0,0.0}
   \put(105.329607,117.075536){\rotatebox{360.0}{\makebox(0,0)[tl]{\strut{}{
    \begin{minipage}[h]{127.655896pt}
\textcolor{inkcol1}{\large{$z$}}\\
\end{minipage}}}}}%

\definecolor{inkcol1}{rgb}{0.0,0.0,0.0}
   \put(143.210327,86.533936){\rotatebox{360.0}{\makebox(0,0)[tl]{\strut{}{
    \begin{minipage}[h]{127.655896pt}
\textcolor{inkcol1}{\large{$w$}}\\
\end{minipage}}}}}%

\definecolor{inkcol1}{rgb}{0.0,0.0,0.0}
   \put(105.203337,79.589136){\rotatebox{360.0}{\makebox(0,0)[tl]{\strut{}{
    \begin{minipage}[h]{127.655896pt}
\textcolor{inkcol1}{\large{$-x$}}\\
\end{minipage}}}}}%

\definecolor{inkcol1}{rgb}{0.0,0.0,0.0}
   \put(114.799787,64.436846){\rotatebox{360.0}{\makebox(0,0)[tl]{\strut{}{
    \begin{minipage}[h]{127.655896pt}
\textcolor{inkcol1}{\large{$y$}}\\
\end{minipage}}}}}%

\definecolor{inkcol1}{rgb}{0.0,0.0,0.0}
   \put(50.907637,63.679236){\rotatebox{360.0}{\makebox(0,0)[tl]{\strut{}{
    \begin{minipage}[h]{127.655896pt}
\textcolor{inkcol1}{\large{$x$}}\\
\end{minipage}}}}}%

\definecolor{inkcol1}{rgb}{0.0,0.0,0.0}
   \put(50.655097,78.578986){\rotatebox{360.0}{\makebox(0,0)[tl]{\strut{}{
    \begin{minipage}[h]{127.655896pt}
\textcolor{inkcol1}{\large{$y$}}\\
\end{minipage}}}}}%

\definecolor{inkcol1}{rgb}{0.0,0.0,0.0}
   \put(111.390527,27.061216){\rotatebox{360.0}{\makebox(0,0)[tl]{\strut{}{
    \begin{minipage}[h]{127.655896pt}
\textcolor{inkcol1}{\large{$w$}}\\
\end{minipage}}}}}%

\definecolor{inkcol1}{rgb}{0.0,0.0,0.0}
   \put(109.875297,56.355636){\rotatebox{360.0}{\makebox(0,0)[tl]{\strut{}{
    \begin{minipage}[h]{127.655896pt}
\textcolor{inkcol1}{\large{$-z$}}\\
\end{minipage}}}}}%

\definecolor{inkcol1}{rgb}{0.0,0.0,0.0}
   \put(224.401347,55.092936){\rotatebox{360.0}{\makebox(0,0)[tl]{\strut{}{
    \begin{minipage}[h]{127.655896pt}
\textcolor{inkcol1}{\large{$y$}}\\
\end{minipage}}}}}%

\definecolor{inkcol1}{rgb}{0.0,0.0,0.0}
   \put(254.453367,29.839136){\rotatebox{360.0}{\makebox(0,0)[tl]{\strut{}{
    \begin{minipage}[h]{127.655896pt}
\textcolor{inkcol1}{\large{$x$}}\\
\end{minipage}}}}}%

\definecolor{inkcol1}{rgb}{0.0,0.0,0.0}
   \put(253.948307,63.426696){\rotatebox{360.0}{\makebox(0,0)[tl]{\strut{}{
    \begin{minipage}[h]{127.655896pt}
\textcolor{inkcol1}{\large{$z$}}\\
\end{minipage}}}}}%

\definecolor{inkcol1}{rgb}{0.0,0.0,0.0}
   \put(254.200837,79.336596){\rotatebox{360.0}{\makebox(0,0)[tl]{\strut{}{
    \begin{minipage}[h]{127.655896pt}
\textcolor{inkcol1}{\large{$w$}}\\
\end{minipage}}}}}%

\definecolor{inkcol1}{rgb}{0.0,0.0,0.0}
   \put(314.052377,83.124676){\rotatebox{360.0}{\makebox(0,0)[tl]{\strut{}{
    \begin{minipage}[h]{127.655896pt}
\textcolor{inkcol1}{\large{$-z$}}\\
\end{minipage}}}}}%

\definecolor{inkcol1}{rgb}{0.0,0.0,0.0}
   \put(317.082837,63.679236){\rotatebox{360.0}{\makebox(0,0)[tl]{\strut{}{
    \begin{minipage}[h]{127.655896pt}
\textcolor{inkcol1}{\large{$w$}}\\
\end{minipage}}}}}%

\definecolor{inkcol1}{rgb}{0.0,0.0,0.0}
   \put(344.104407,52.820096){\rotatebox{360.0}{\makebox(0,0)[tl]{\strut{}{
    \begin{minipage}[h]{127.655896pt}
\textcolor{inkcol1}{\large{$-x$}}\\
\end{minipage}}}}}%

\definecolor{inkcol1}{rgb}{0.0,0.0,0.0}
   \put(313.294757,29.334046){\rotatebox{360.0}{\makebox(0,0)[tl]{\strut{}{
    \begin{minipage}[h]{127.655896pt}
\textcolor{inkcol1}{\large{$y$}}\\
\end{minipage}}}}}%

\definecolor{inkcol1}{rgb}{0.0,0.0,0.0}
   \put(300.162777,49.284566){\rotatebox{360.0}{\makebox(0,0)[tl]{\strut{}{
    \begin{minipage}[h]{127.655896pt}
\textcolor{inkcol1}{\large{$z$}}\\
\end{minipage}}}}}%

\definecolor{inkcol1}{rgb}{0.0,0.0,0.0}
   \put(270.110737,49.284566){\rotatebox{360.0}{\makebox(0,0)[tl]{\strut{}{
    \begin{minipage}[h]{127.655896pt}
\textcolor{inkcol1}{\large{$w$}}\\
\end{minipage}}}}}%

\definecolor{inkcol1}{rgb}{0.0,0.0,0.0}
   \put(223.896267,89.690656){\rotatebox{360.0}{\makebox(0,0)[tl]{\strut{}{
    \begin{minipage}[h]{127.655896pt}
\textcolor{inkcol1}{\large{$z$}}\\
\end{minipage}}}}}%

\definecolor{inkcol1}{rgb}{0.0,0.0,0.0}
   \put(259.251607,118.985096){\rotatebox{360.0}{\makebox(0,0)[tl]{\strut{}{
    \begin{minipage}[h]{127.655896pt}
\textcolor{inkcol1}{\large{$w$}}\\
\end{minipage}}}}}%

\definecolor{inkcol1}{rgb}{0.0,0.0,0.0}
   \put(294.354407,93.731266){\rotatebox{360.0}{\makebox(0,0)[tl]{\strut{}{
    \begin{minipage}[h]{127.655896pt}
\textcolor{inkcol1}{\large{$w$}}\\
\end{minipage}}}}}%

\definecolor{inkcol1}{rgb}{0.0,0.0,0.0}
   \put(264.302367,94.993966){\rotatebox{360.0}{\makebox(0,0)[tl]{\strut{}{
    \begin{minipage}[h]{127.655896pt}
\textcolor{inkcol1}{\large{$-z$}}\\
\end{minipage}}}}}%

\definecolor{inkcol1}{rgb}{0.0,0.0,0.0}
   \put(350.922947,92.973666){\rotatebox{360.0}{\makebox(0,0)[tl]{\strut{}{
    \begin{minipage}[h]{127.655896pt}
\textcolor{inkcol1}{\large{$y$}}\\
\end{minipage}}}}}%

\definecolor{inkcol1}{rgb}{0.0,0.0,0.0}
   \put(306.223687,101.054886){\rotatebox{360.0}{\makebox(0,0)[tl]{\strut{}{
    \begin{minipage}[h]{127.655896pt}
\textcolor{inkcol1}{\large{$x$}}\\
\end{minipage}}}}}%

 \end{picture}
\endgroup